\documentclass[a4paper,12pt]{article}
\usepackage[colorlinks=true, pdfstartview=FitV, linkcolor=blue, citecolor=red, urlcolor=blue]{hyperref}
\usepackage{amsfonts,amssymb}
\usepackage{amsmath,amsthm}
\usepackage{graphicx}
\usepackage{graphics}
\usepackage{epsfig}
\usepackage{epstopdf}
\usepackage{geometry}
\usepackage{color}
\usepackage{url}
\usepackage[dvipsnames]{xcolor}
\definecolor{refkey}{rgb}{1,0,0}
\definecolor{labelkey}{rgb}{1,0,0}
\usepackage{caption}

\usepackage{float}
\usepackage{caption,subfig,color}
\usepackage{tikz, tkz-tab}
\usetikzlibrary{calc}
\usetikzlibrary{matrix}
\usetikzlibrary{backgrounds}
\usepackage{pgfplots}
\pgfplotsset{compat=newest}
\pgfplotsset{plot coordinates/math parser=false}
\usetikzlibrary{plotmarks}
\usetikzlibrary{arrows.meta}
\usepgfplotslibrary{patchplots}

\newcommand\figurescale{0.6}

\usepackage[nointegrals]{wasysym}

\newtheorem{theorem}{Theorem}[section]

\newtheorem{remark}[theorem]{Remark}
\theoremstyle{definition}
\newtheorem{definition}[theorem]{Definition}

\newcommand{\logplottext}[5]{
  \begin{tikzpicture}[inner frame sep = 0pt]
    \centering
    \node[inner sep=0] at (0,0){\includegraphics[width=#2]{#1}};
    \begin{scope}[scale={#2/#3},axes/.style={scale=#4}]
      \draw (-1.5,-1) rectangle (1.5,1);
      \draw foreach[count=\n from -1,evaluate={\m=int(\n*2)}] \x in {-1,0.,...,1}{(\x,-1) node[below,axes]{$\m$}};
      \draw foreach[count=\n from -1,evaluate={\m=int(\n*2)}] \y in {-1,0.,...,1}{(-1.5,\y) node[left,axes]{$\m$}};
      \node[axes,white,left] at (1.5,-.75) {$\mathbf{#5}$};      
    \end{scope}
  \end{tikzpicture}
}

\newcommand{\myTDrepr}[3]{
  \begin{tikzpicture}
    \centering
    \node at (0,0){\includegraphics[width=#2]{#1}};
    \begin{scope}[scale={#2/#3},axes/.style={scale=1}]
      \node[axes] at (-0.81,-2.56) {$0$};
      \node[axes] at (0.77,-2.85) {$2$};
      \node[axes] at (2.07,-2.85) {$-2$};
      \node[axes] at (2.67,-2.05) {$0$};
      \node[axes] at (3.25,-1.3) {$2$};
      \node[axes] at (-2.48,-2.3) {$-2$};
      \node[axes] at (-3.23,-1.7) {$0$};
      \node[axes] at (-3.23,-0.97) {$2$};
      \node[axes] at (-3.23,-0.25) {$4$};
      \node[axes] at (-3.23,0.5) {$6$};
      \node[axes] at (-3.23,1.2) {$8$};
      \node at (3.07,-3.08) {$y$};
      \node at (-0.41,-3.2) {$x$};
      \node[rotate=90] at (-3.63,-0.3) {$t$};
    \end{scope}
  \end{tikzpicture}
}

\title{Discrete transparent boundary conditions \\
for the two-dimensional leap-frog scheme}

\author{Christophe Besse\thanks{Institut de Math\'ematiques de Toulouse ; 
UMR5219, Universit\'e de Toulouse ; CNRS, Universit\'e Paul Sabatier, F-31062 Toulouse Cedex 9, 
France. Email: {\tt christophe.besse@math.univ-toulouse.fr}. Research of the author was supported 
by ANR project NABUCO, ANR-17-CE40-0025.}, 
Jean-Fran\c{c}ois Coulombel\thanks{Institut de Math\'ematiques de Toulouse ; 
UMR5219, Universit\'e de Toulouse ; CNRS, Universit\'e Paul Sabatier, F-31062 Toulouse Cedex 9, 
France. Email: {\tt jean-francois.coulombel@math.univ-toulouse.fr}. Research of the author was supported 
by ANR project NABUCO, ANR-17-CE40-0025.}, 
Pascal Noble\thanks{Institut de Math\'ematiques de Toulouse ; 
UMR5219, Universit\'e de Toulouse ; CNRS, INSA, F-31077 Toulouse, France. 
Email: {\tt pascal.noble@math.univ-toulouse.fr}. Research of the author was supported 
by ANR project NABUCO, ANR-17-CE40-0025.}}

\date{\today}
\begin{document}
\maketitle

\begin{abstract}
We develop a general strategy in order to implement (approximate) discrete transparent boundary conditions for finite difference 
approximations of the two-dimensional transport equation. The computational domain is a rectangle equipped with a Cartesian 
grid. For the two-dimensional leap-frog scheme, we explain why our strategy provides with explicit numerical boundary conditions 
on the four sides of the rectangle and why it does not require prescribing any condition at the four corners of the computational 
domain. The stability of the numerical boundary condition on each side of the rectangle is analyzed by means of the so-called 
normal mode analysis. Numerical investigations for the full problem on the rectangle show that strong instabilities may occur 
when coupling stable strategies on each side of the rectangle. Other coupling strategies yield promising results.
\end{abstract}

\section{Introduction}

In this article, we are concerned with the construction and numerical implementation of discrete transparent boundary conditions for 
linear transport equations. This research area has been increasingly active since the pioneering work by Engquist and Majda \cite{engquist-majda} 
and our goal here is to explain why following the same `small frequency approximation' strategy can be successful at the fully discrete level 
when one deals with a rectangular computational domain. Namely, we shall construct (approximate) fully discrete transparent boundary 
conditions for the two-dimensional leap-frog scheme on a rectangle. Our transparent boundary conditions will not be exact since the 
derivation of such conditions on a rectangle seems to be out of reach. Our strategy is to start from the \emph{exact} transparent boundary 
conditions for a \emph{half-space} and, as in \cite{engquist-majda}, to localize them with respect to the tangential variable by means of 
a suitable `small frequency approximation'. This general strategy may be applied to any finite difference scheme. The nice feature of the 
leap-frog approximation is that it does not dissipate high frequency signals -and therefore allows for a precise analysis of wave reflections 
on each side of the rectangle- and, moreover, its stencil exhibits some `dimensional splitting'. The latter feature will be helpful in our 
construction since we shall not have to develop a specific numerical treatment for the four corners of the computational domain. Hence 
we shall focus on the numerical boundary conditions on each of the four sides and on their coupling through the mesh points that are 
closest to the four corners.

Exact Discrete Transparent Boundary Conditions (DTBC in what follows) are analytically available in two space dimensions only when the 
computational domain is a half-space or the product of a half-line by a torus. We refer for instance to \cite{AES,AESS,aabes,jfc} and references 
therein for various aspects of the theory, the main features of which are recalled below. In either case, the geometry is cylindrical with no 
tangential boundary so DTBC can be derived by performing a Fourier transform (or a Fourier series decomposition) with respect to the 
spatial tangential variables, which yields a one-dimensional problem that is parametrized by the tangential frequency\footnote{The same 
strategy applies when the computational domain is a disk by using polar coordinates.}. In that framework, the two-dimensional DTBC thus 
take the form of \emph{tangential pseudo-differential operators} that are intrinsically \emph{nonlocal}. In order to be implemented on a 
rectangular computational domain, one therefore has to make the DTBC \emph{local}, which we achieve here by using a suitable `small 
frequency approximation', which is meaningful (at least formally) for smooth solutions. The frequency cut-off amounts to retaining only 
finitely many terms in the Taylor expansion of the symbol of the DTBC operator and by rewriting the Taylor expansion as a trigonometric 
polynomial (up to higher order terms). This cut-off strategy is implemented here so that the stencil used for the numerical boundary conditions 
is not wider than that of the interior numerical scheme. Higher order strategies could be relevant but would require a specific treatment when 
getting closer to the corner. We do not investigate this issue here.

As in all numerical approximation problems, a crucial feature for the efficiency of the nume\-rical scheme is stability. Following the 
general result of \cite{jfc}, the DTBC for the leap-frog scheme on a half-space do not meet the strongest stability properties because 
the leap-frog scheme exhibits \emph{glancing} wave packets (even in one space dimension). On a half-space, the exact DTBC meet 
some kind of \emph{neutral} stability which calls for special care since a slight modification in the numerical boundary conditions may 
yield \emph{violent instabilities}. In what follows, we show that the zero, first and second order tangential frequency cut-off -as described 
above- maintain neutral stability when implemented on a half-space. In other words, our approximate, local, transparent boundary conditions 
are as stable as the exact ones on a half-space and we therefore feel rather safe to implement them on a rectangular geometry. If instabilities 
occur in the rectangular geometry, they will necessarily stem from the coupling between the various numerical boundary conditions on each 
side of the rectangle. A thorough stability analysis on the rectangle seems to be out of reach at the present time (it is not even fully understood 
at the continuous level !), but the numerical investigations that we report below show that stability in the rectangle should not be taken for 
granted. Some coupling strategies exhibit good stability -and convergence- properties while at least one displays violent instabilities. In a 
future work, we shall extend our strategy to implicit discretizations of dispersive equations.

\section{The transport equation and the leap-frog approximation}

In this article, we consider the linear advection equation in two space dimensions:
\begin{equation}\label{eq2}
\left\{\begin{array}{ll}
\displaystyle
\partial_t u +\mathbf{c} \cdot \nabla  u=0 \, ,\quad t\geq 0 \, ,\quad (x,y) \in\mathbb{R}^2 \, ,\\
\displaystyle
u|_{t=0}=u_0 \, ,
\end{array}\right.
\end{equation}
where the space coordinates are denoted $(x,y)$, the velocity reads $\mathbf{c}=(c_x,c_y)^T$ and $\nabla=(\partial_x,\partial_y)^T$. It is always 
assumed from now on that the velocity $\mathbf{c}$ is nonzero. To fix ideas, we shall always take $c_x$ and $c_y$ to be nonnegative, but this has 
no impact on the analysis below. The solution to \eqref{eq2} reads:
$$
u(t,x,y) \, = \, u_0(x-c_x \, t,y-c_y \, t) \, ,
$$
so if the initial condition $u_0$ is supported, say in the square $[-L,L] \times [-L,L]$, the solution vanishes in the larger square $[-2L,2L] \times [-2L,2L]$ 
after some time $T \sim L/|\mathbf{c}|$. Given a finite difference approximation to \eqref{eq2} associated with a Cartesian grid, our goal in this article is 
to propose a systematic construction of \emph{transparent} numerical boundary conditions\footnote{Sometimes these are also called non-reflecting or 
absorbing numerical boundary conditions.} in order to simulate the advection of the initial condition through the computational domain and eventually 
its exit through the numerical boundary. As a test case, we choose the two-dimensional leap-frog scheme (see \eqref{LF2d} below) which is explicit and 
propagates high frequency signals without dissipation. More simple outflow strategies, such as extrapolation boundary conditions, may be implemented 
for dissipative schemes (e.g., the Lax-Wendroff scheme). Even though extrapolation boundary conditions are much cheaper from a computational point 
of view and are quite efficient in terms of absence of wave reflection (for dissipative schemes), they only give poor intuition for more complex wave 
propagation phenomena modeled for instance by the Schr\"odinger or Airy equations.

In all what follows, the computational domain is a fixed rectangle $[x_\ell,x_r] \times [y_b,y_t]$ ($\ell,r,b,t$ stand for left, right, bottom and top). We consider 
some space steps $\delta x, \delta y >0$ such that the mesh ratios:
\begin{equation*}
\frac{x_r-x_\ell}{\delta x} \, =: \, J+1 \, ,\quad \frac{y_t-y_b}{\delta y} \, =: \, K+1 \, ,
\end{equation*} 
define some integers $J$ and $K$ (that are meant to be large). Eventually, the time step $\delta t>0$ will always be chosen so that the parameters
\begin{equation*}
\mu_x \, := \, c_x \, \dfrac{\delta t}{\delta x} \, ,\quad \mu_y \, := \, c_y \, \dfrac{\delta t}{\delta y} \, ,
\end{equation*}
are fixed. These parameters should satisfy some stability requirement, the so-called Courant-Friedrichs-Lewy condition, as explained below. The grid 
points are defined as $(x_j,y_k)$, with:
$$
\forall \, j=0,\dots,J+1 \, ,\quad \forall \, k =0,\dots,K+1 \, ,\quad x_j \, := \, x_\ell +j \, \delta x \, ,\quad y_k \, := \, y_b +k \, \delta y \, .
$$
The interior points correspond to $1 \le j \le J$ and $1 \le k \le K$. The four sides of the rectangle (the numerical boundary) correspond to 
$j \, \in \{ 0,J+1\}$ and $k \in \{ 0,K+1\}$. The grid is depicted on Figure \ref{fig:maillage} below.

\begin{figure}[htbp]
\begin{center}
\begin{tikzpicture}[scale=2,>=latex]
\draw [thin, dashed] (-3.5,-1.5) grid [step=0.5] (3.5,2);
\draw (-3.5,2) -- (3.5,2);
\draw (-3.5,-1.5) -- (3.5,-1.5);
\draw (-3.5,-1.5) -- (-3.5,2);
\draw (3.5,-1.5) -- (3.5,2);
\draw[black,->] (-4,0) -- (4,0) node[below] {$x$};
\draw[black,->] (0,-2)--(0,2.5) node[right] {$y$};
\draw (-3.8,-1.5) node[right]{$y_0$};
\draw (-3.8,-1) node[right]{$y_1$};
\draw (-3.9,1.5) node[right]{$y_K$};
\draw (-3.8,2.15) node[right]{$y_{K+1}$};
\draw (-3.6,-1.65) node[right]{$x_0$};
\draw (-3.1,-1.65) node[right]{$x_1$};
\draw (-2.6,-1.65) node[right]{$x_2$};
\draw (2.35,-1.65) node[right]{$x_{J-1}$};
\draw (2.85,-1.65) node[right]{$x_J$};
\draw (3.35,-1.65) node[right]{$x_{J+1}$};
\node (centre) at (-3.5,-1.5){$\bullet$};
\node (centre) at (-3.5,2){$\bullet$};
\node (centre) at (3.5,-1.5){$\bullet$};
\node (centre) at (3.5,2){$\bullet$};
\node (centre) at (-3,-1.5){$\times$};
\node (centre) at (-2.5,-1.5){$\times$};
\node (centre) at (3,-1.5){$\times$};
\node (centre) at (2.5,-1.5){$\times$};
\node (centre) at (-3,2){$\times$};
\node (centre) at (-2.5,2){$\times$};
\node (centre) at (3,2){$\times$};
\node (centre) at (2.5,2){$\times$};
\node (centre) at (-3.5,-1){$\times$};
\node (centre) at (-3.5,1.5){$\times$};
\node (centre) at (3.5,-1){$\times$};
\node (centre) at (3.5,1.5){$\times$};
\node (centre) at (0.5,1){${\color{blue}\otimes}$};
\node (centre) at (1.5,1){${\color{blue}\otimes}$};
\node (centre) at (1,0.5){${\color{blue}\otimes}$};
\node (centre) at (1,1.5){${\color{blue}\otimes}$};
\node (centre) at (1,1){${\color{blue}\otimes}$};
\end{tikzpicture}
\caption{The grid on $[x_\ell,x_r] \times [y_b,y_t]$. Interior points correspond to the intersection of dashed lines, the boundary points correspond 
to the full black lines (crosses) and the four corners are labeled with bullets. The interior stencil for the scheme \eqref{LF2d} is marked with (blue) 
circled crosses.}
\label{fig:maillage}
\end{center}
\end{figure}
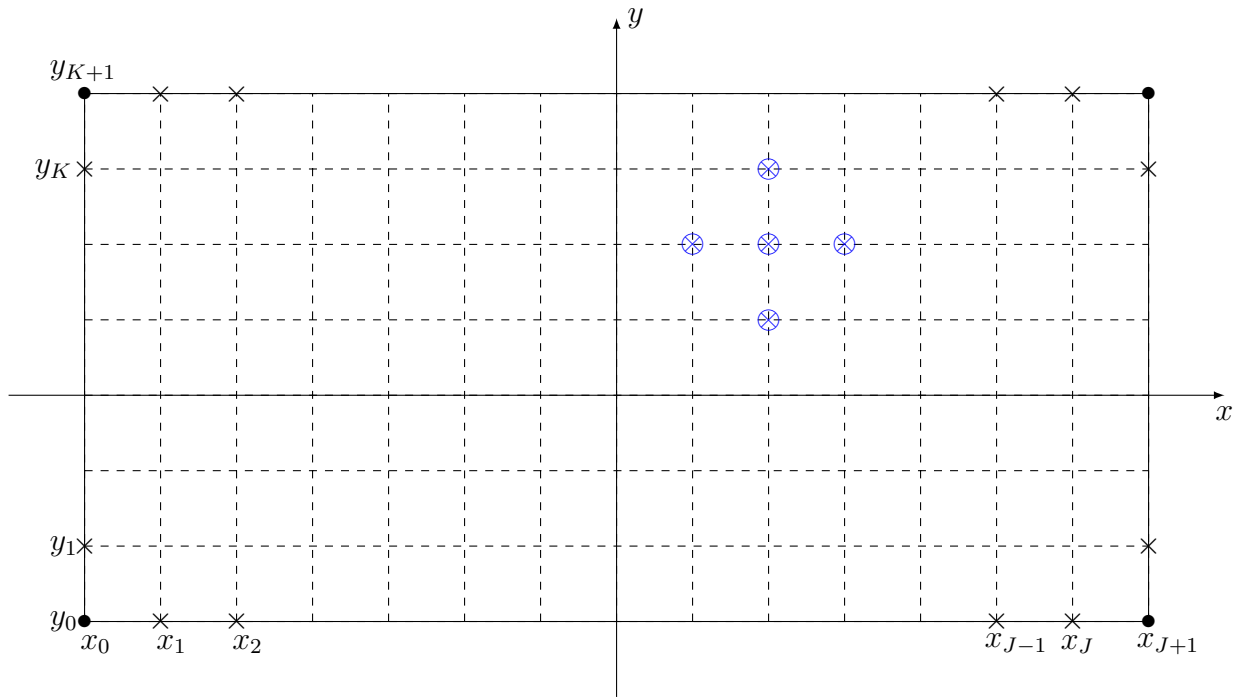

Letting now $u_{j,k}^n$ denote the approximation of the exact solution to \eqref{eq2} at $(x_j,y_k)$ 
and time $n \, \delta t$, the two-dimensional leap-frog scheme reads:
\begin{equation}
\label{LF2d}
\displaystyle
u_{j,k}^{n+2}-u_{j,k}^n +\mu_x \, \big( u_{j+1,k}^{n+1}-u_{j-1,k}^{n+1} \big) +\mu_y \, \big( u_{j,k+1}^{n+1}-u_{j,k-1}^{n+1} \big) \, = \, 0 \, ,
\end{equation}
which holds for $1 \le j \le J$ and $1 \le k \le K$, that is at all interior points. The definition of the numerical scheme requires to prescribe 
numerical boundary conditions for $u_{0,k}^{n+2}$, $u_{J+1,k}^{n+2}$, $k=1,\dots,K$, and $u_{j,0}^{n+2}$, $u_{j,K+1}^{n+2}$, $j=1,\dots,J$. 
A crucial observation for what follows is that the update of the numerical solution in the interior does not require to prescribe any numerical 
boundary condition at the corners (this is because $u_{j,k}^{n+2}$ does not depend on $u_{j \pm 1,k \pm 1}^\sigma$ at any earlier time 
index $\sigma$). As depicted in Figure \ref{fig:maillage}, the stencil of \eqref{LF2d} is the union of the split `1D stencils' which correspond 
to taking either $\mu_x=0$ or $\mu_y=0$. In all what follows, the numerical solution $(u_{j,k}^n)$ will therefore be defined for $(j,k) \in 
\{0,\dots,J+1\} \times \{0,\dots,K+1\} \setminus \{ (0,0), (0,K+1), (J+1,0), (J+1,K+1) \}$. We shall be careful in our definition of the numerical 
boundary conditions to make sure that the four corner values are not involved either.

Let us note that on the whole space $\mathbb{Z}^2$, the $\ell^2$-stability condition for \eqref{LF2d} reads:
\begin{equation}
\label{cfl2D}
\mu_x +\mu_y \le \text{cfl} <1 \, ,
\end{equation}
where cfl is a free parameter to be chosen. (Recall that we consider the case $c_x,c_y \ge 0$, hence $\mu_x,\mu_y \ge 0$, so we do not write 
absolute values in \eqref{cfl2D}.) The inequality \eqref{cfl2D} gives an upper bound for $\delta t$ in terms of $\delta x,\delta y$. The bound 
\eqref{cfl2D}, or its one-dimensional analogue that we shall discuss below, is always assumed to hold from now on.

We wish to construct transparent boundary conditions for \eqref{LF2d} on each side of the rectangle. For ease of reading and also for future 
use on the two-dimensional problem \eqref{LF2d}, we first go back to the one-dimensional problem and briefly recall the construction of DTBC 
for the one-dimensional leap-frog scheme. We discuss the various possible implementations of the DTBC and approximations by means of 
sums of exponentials (this technique was used for instance in \cite{AES}). Once the tools have been set in one space dimension, we shall 
feel free to use them in two space dimensions whenever it makes sense.

\section{Discrete transparent boundary conditions for the one-dimensional leap-frog scheme}

\subsection{The numerical scheme}

We now consider the one-dimensional transport equation:
\begin{equation}
\label{eq1}
\left\{\begin{array}{ll}
\displaystyle
\partial_t u+c_x \, \partial_x u=0 \, ,\\
\displaystyle
u|_{t=0}=u_0 \, ,
\end{array} \right.
\end{equation}
with an initial condition $u_0$ that is supported in an interval $[x_\ell;x_r]$. For concreteness, we shall assume from now on $c_x>0$ so that the 
initial condition is transported towards the right and eventually exits the initial support $[x_\ell;x_r]$. We are interested here in simulating the solution 
to \eqref{eq1} on the fixed interval $[x_\ell;x_r]$ by imposing suitable non-reflecting boundary conditions at the boundary points $x_\ell,x_r$.

As in the previous section, we choose a mesh size $\delta x>0$ such that
\begin{equation*}
\dfrac{x_r-x_\ell}{\delta x} \, =: J+1 \, ,
\end{equation*} 
is an integer. The time step $\delta t>0$ is chosen in such a way that the parameter $\mu_x$ defined by:
\[
\mu_x \, := \, c_x \, \dfrac{\delta t}{\delta x} \, ,
\]
is a fixed (positive) constant. Keeping the notation $x_j :=x_\ell +j \, \delta x$, $j=0,\dots,J+1$, and letting now $u_j^n$ denote the approximation of 
the solution to \eqref{eq1} at $x_j$ and time $n \, \delta t$, the one-dimensional leap-frog scheme reads
\begin{equation}
\label{LF1d}
u_j^{n+2} -u_j^n +\mu_x \, \big( u_{j+1}^{n+1}-u_{j-1}^{n+1} \big) =0 \, ,\quad j=1,\dots,J \, .
\end{equation}
The parameter $\mu_x$ is fixed by imposing:
\begin{equation}
\label{cfl1D}
\mu_x \le \text{cfl} <1 \, ,
\end{equation}
which is a necessary and sufficient stability requirement for \eqref{LF1d} on the whole real line $j \in \mathbb{Z}$, see \cite{gko,RM}. Since $\mu_x$ 
is positive in our framework, we thus fix $\mu_x \in (0,1)$ and let $\delta x,\delta t$ vary accordingly.

The numerical sheme \eqref{LF1d} on the interval $1\leq j \leq J$ requires a definition for the boundary values $u_0^{n+1}$, $u_{J+1}^{n+1}$, and 
a definition for the initial data $(u_j^0)$, $(u_j^1)$. The initial condition $(u_j^0)$ is defined by setting
$$
\forall \, j=0,\dots,J+1 \, ,\quad u_j^0 \, := \, u_0(x_j) \, ,
$$
where $u_0$ is the initial condition for \eqref{eq1}. In the numerical simulations reported below, we have chosen $[x_\ell,x_r]=[-3,3]$ and $u_0(x) 
=\exp (-10\, x^2)$. The first time step value $(u_j^1)$ is defined by imposing the second order Lax-Wendroff scheme in the interior domain, namely:
$$
\forall \, j=1,\dots,J \, ,\quad 
u_j^1 \, := \, u_j^0 -\dfrac{\mu_x}{2} \, (u_{j+1}^0-u_{j-1}^0) +\dfrac{\mu_x^2}{2} \, (u_{j+1}^0-2\, u_j^0 +u_{j-1}^0) \, .
$$
The boundary values $u_0^1$ and $u_{J+1}^1$ are set equal to zero for simplicity (which is almost exact for the Gaussian initial condition that we 
consider). We now recall how to derive the DTBC for the scheme \eqref{LF1d}.

\subsection{Derivation of DTBC}

The goal of this paragraph is to derive the DTBC for the leap-frog scheme \eqref{LF1d}. Following the analysis of \cite{ehrhardt-arnold} and subsequent 
works, we consider the numerical scheme
\begin{equation}
\label{leapfrog1D}
u_j^{n+2} -u_j^n +\mu_x \, \big( u_{j+1}^{n+1}-u_{j-1}^{n+1} \big) =0 \, ,\quad j \in \mathbb{Z} \, ,
\end{equation}
on the whole real line $\mathbb{Z}$ and assume that the initial conditions $(u_j^0)_{j \in \mathbb{Z}}$, $(u_j^1)_{j \in \mathbb{Z}}$ belong to $\ell^2$ 
and vanish outside of the interval $\{ 1,\dots,J \}$. The stability condition \eqref{cfl1D} is assumed to hold. The DTBC are first computed in terms of the 
so-called $\mathcal{Z}$-transform of the sequences $(u_j^n)_{n\in\mathbb{N}}$. Let us briefly recall the definition of the $\mathcal{Z}$-transform (we 
do not pay too much attention here to the convergence of the Laurent series involved in the computations below and rather remain at the level of formal 
series).

\begin{definition}\label{defZT}
Let $(u_n)_{n\in\mathbb{N}} \in \mathbb{C}^\mathbb{N}$ be a sequence. The $\mathcal{Z}$-transform of $(u_n)_{n\in\mathbb{N}}$, which is denoted 
$\widehat{u}$, is defined by:
\begin{equation*}
\widehat{u}(z) \, := \, \sum_{n=0}^\infty \, \dfrac{u_n}{z^n} \, ,
\end{equation*}
the series being well-defined and holomorphic in $\{z \in \mathbb{C} \, , \, |z|>R \}$ where $1/R$ is the convergence radius of $\sum_{n=0}^\infty u_n 
\, w^n$ (here $w$ is a placeholder for $1/z$ and the radius $R$ may be infinite).
\end{definition}

We now let $\hat{u}_j(z)$, $j \in \mathbb{Z}$, denote the $\mathcal{Z}$-transform of the sequence $(u^n_j)_{n\in\mathbb{N}}$:
\begin{equation*}
\hat{u}_j(z) \, := \, \sum_{n=0}^\infty \, u_j^n \, z^{-n} \, .
\end{equation*}
It is shown in \cite{jfc} that the above series converges for $|z|>1$ thanks to the $\ell^2$-stability of the leap-frog scheme\footnote{This is where the 
condition $|\mu_x|<1$ becomes necessary.}. Moreover, the sequence $(\hat{u}_j(z))_{j\in\mathbb{Z}}$ is square integrable. We apply the 
$\mathcal{Z}$-transform to the numerical scheme \eqref{leapfrog1D} and obtain:
\begin{equation}
\label{eq:ZT_LF1d}
\big( z-z^{-1} \big) \, \hat{u}_j +\mu_x \, \big( \hat{u}_{j+1}-\hat{u}_{j-1} \big) \, = \, 0 \, ,
\end{equation}
for $j \le 0$ and $j \ge J+1$. (On the interval $\{1,\dots,J \}$, the initial conditions give rise to a nonzero source term on the right hand side of 
\eqref{eq:ZT_LF1d}.) The equation \eqref{eq:ZT_LF1d} is a recurrence relation and we therefore look for the solutions to the characteristic equation:
\begin{equation}
\label{eq:rec1d}
\kappa^2 +\dfrac{z-z^{-1}}{\mu_x} \, \kappa -1 \, = \, 0 \, .
\end{equation}
The two roots to \eqref{eq:rec1d} are given by
\begin{equation}
\label{eq:sol_rec1d}
\kappa_\pm \, := \, \dfrac{z^{-1}-z \pm \sqrt{z^2+2 \, (2 \, \mu_x^2-1)+z^{-2}}}{2 \, \mu_x} \, ,
\end{equation}
with the standard definition of the square root (the branch cut is $\mathbb{R}^-$). When $|z|\to+\infty$, $|\kappa_-| \to +\infty$ and $|\kappa_+| \to 0$. 
Moreover, $\kappa_\pm$ do not belong to $\mathbb{S}^1$ for $|z|>1$ so we have $|\kappa_+|<1$ and $|\kappa_-|>1$ for $|z|>1$. We thus relabel 
these two roots as a \emph{stable} and \emph{unstable} one and write from now on:
\begin{equation}
\label{eq:stab_root}
\kappa_s^0(z) \, := \, \dfrac{z^{-1}-z+\sqrt{z^2+2 \, (2 \, \mu_x^2-1)+z^{-2}}}{2 \, \mu_x} \, , \, 
\kappa_u^0(z) \, := \, \dfrac{z^{-1}-z-\sqrt{z^2+2 \, (2 \, \mu_x^2-1)+z^{-2}}}{2 \, \mu_x} \, .
\end{equation}
The reason for the superscript $0$ will be clear in the next Section when we deal with the two-dimensional problem. The recurrence relation 
\eqref{eq:ZT_LF1d} and the decay at infinity of $(\hat{u}_j(z))_{j\in\mathbb{Z}}$ implies:
\begin{equation}
\label{eq:dtbc1d}
\hat{u}_{J+1}(z) \, = \, \kappa_s^0(z) \, \hat{u}_J(z) \, ,\quad \hat{u}_0(z) \, = \, \dfrac{1}{\kappa_u^0(z)} \, \hat{u}_1(z) \, .
\end{equation}
It now remains to compute the Laurent series expansion of both $\kappa_s^0$ and $1/\kappa_u^0$ and to perform the inverse $\mathcal{Z}$-transform 
in order to write the numerical boundary conditions \eqref{eq:dtbc1d} in the original discrete time variable $n$. Since we have $\kappa_s^0(z) \, \kappa_u^0 
(z)=-1$, see \eqref{eq:rec1d}, we can equivalently rewrite \eqref{eq:dtbc1d} as:
\begin{equation}
\label{eq:dtbc1d'}
\hat{u}_{J+1}(z) \, = \, \kappa_s^0(z) \, \hat{u}_J(z) \, ,\quad \hat{u}_0(z) \, = \,-\kappa_s^0(z) \, \hat{u}_1(z) \, .
\end{equation}
There are at least two ways to compute the Laurent series expansion of $\kappa_s^0$ for the leap-frog scheme, which we now detail in order to compare 
the possible benefits and drawbacks of each method.
\bigskip

\paragraph{An explicit expansion in terms of Legendre polynomials.} The first method to obtain this expansion relies on the explicit formula \eqref{eq:stab_root} 
and on the expansion
\begin{equation}
\label{eq:Legendre1}
\dfrac{1}{\sqrt{1 -2 \, x \, t+t^2}} \, = \, \sum_{n=0}^\infty \, P_n(x) \, t^n \, ,
\end{equation}
where $P_n$ denotes the $n^{\mathrm{th}}$-Legendre polynomial \cite{szego}. This family of polynomials can be equivalently defined by the first two terms 
$P_0(x):=1$, $P_1(x):=x$ and by the recurrence formula:
\begin{equation}
\label{eq:Leg}
\forall \, n \ge 1 \, ,\quad (n+1) \, P_{n+1}(x) \, := \, (2n+1) \, x \, P_n(x) -n \, P_{n-1}(x) \, .
\end{equation}
Starting from \eqref{eq:stab_root}, we have
\begin{align*}
\sqrt{z^2+2 \, (2\, \mu_x^2-1)+z^{-2}} &= \, z \, \dfrac{1-2\, (1-2\, \mu_x^2) \, z^{-2}+z^{-4}}{\sqrt{1-2 \, (1-2 \, \mu_x^2) \, z^{-2}+z^{-4}}} \\
 &= \, z \, (1-2 \, \alpha_x \, z^{-2}+z^{-4}) \, \sum_{n=0}^\infty \, P_n(\alpha_x) \, z^{-2n} \, ,
\end{align*}
with $\alpha_x:=1-2\, \mu_x^2 \in (-1,1)$. We therefore derive the expansion
\begin{equation}
\label{expressionkappas0}
2 \, \mu_x \, z \, \kappa_s^0 \, (z) \, = \, 1 -z^2 +\big( z^2 -2 \, \alpha_x +z^{-2} \big) \, \sum_{n=0}^\infty \, P_n(\alpha_x) \, z^{-2n} \, .
\end{equation}
We focus below on the right boundary condition at point $x_{J+1} =x_r$. Going back to \eqref{eq:dtbc1d'} and using the above expansion 
\eqref{expressionkappas0} for $\kappa_s^0$, the DTBC \eqref{eq:dtbc1d'} reads
\[
2 \, \mu_x \, z \, \hat{u}_{J+1}(z) \, = \, \Big( 1+\sum_{n=0}^\infty B_n \, z^{-2n} \Big) \, \hat{u}_J(z) \, ,
\]
with
\[
\forall \, n \in \mathbb{N} \, ,\quad B_n \, := \, P_{n+1}(\alpha_x) \, -2 \, \alpha_x \, P_n(\alpha_x) +P_{n-1}(\alpha_x) \, ,
\]
and it is understood that $P_{-1}$ vanishes (in the expression of $B_0$). We compute $B_0=-\alpha_x$, and thanks to the recurrence relation \eqref{eq:Leg}, 
we get the shorter expression
\[
\forall \, n \ge 1 \, ,\quad B_n \, = \, \dfrac{P_{n-1}(\alpha_x)-P_{n+1}(\alpha_x)}{2 \, n+1} \, ,
\]
and finally
\begin{equation}
\label{eq:dtbc1d_2}
2 \, \mu_x \, z \, \hat{u}_{J+1}(z) \, = \, \left( 2 \, \mu_x^2 +\sum_{n=1}^\infty \dfrac{P_{n-1}(\alpha_x)-P_{n+1}(\alpha_x)}{2\, n+1}  \, z^{-2n} \right) 
\, \hat{u}_J(z) \, .
\end{equation}
Recall that in \eqref{eq:dtbc1d_2}, we use the notation $\alpha_x=1-2\, \mu_x^2$.

We now look for an efficient way to implement the coefficients of the above series without using the expression of the Legendre polynomials (which 
would be computationally rather cheap anyway). Let us define the sequence $(s_n^0)_{n\geq 0}$ such that $s_0^0 :=\mu_x$, $s_1^0 :=\mu_x \, 
(1-\mu_x^2)$, and $s_n^0 := (P_{n-1}(\alpha_x)-P_{n+1}(\alpha_x))/((4\, n+2) \, \mu_x)$ for $n \geq 2$. After various simplifications using 
\eqref{eq:Leg}, we obtain
\begin{equation}
\label{eq:sn}
\forall \, n \ge 2 \, ,\quad s_n^0 \, = \, \dfrac{2\, n-1}{n+1} \, (1-2\, \mu_x^2) \, s_{n-1}^0 -\dfrac{n-2}{n+1} \, s_{n-2}^0 \, .
\end{equation}
The DTBC \eqref{eq:dtbc1d_2} therefore reads
\begin{equation}
\label{eq:dtbc1d_3}
\hat{u}_{J+1}(z) \, = \, \sum_{n=0}^\infty \, s_n^0 \, z^{-2\, n-1} \, \hat{u}_J(z) \, .
\end{equation}
Performing the inverse $\mathcal{Z}$-transform, we get\footnote{It is understood in \eqref{eq:dtbc1d_3'} that the summation holds over all integers 
$m$ between $0$ and $(n+1)/2$. In case $n$ is even, then the summation holds over all integers $m$ between $0$ and $n/2$.}:
\begin{equation}
\label{eq:dtbc1d_3'}
u_{J+1}^{n+2} \, = \, \sum_{0 \le m \le (n+1)/2} \, s_m^0 \, u_J^{n+1-2\, m} \, .
\end{equation}
The recurrence relation \eqref{eq:sn} may be efficiently implemented in a computer code, which gives access to the coefficients in the numerical 
boundary condition \eqref{eq:dtbc1d_3'} up to any prescribed final time index $N_f$. Going back to \eqref{eq:dtbc1d'}, we also derive the `left' 
numerical boundary condition:
\begin{equation}
\label{eq:dtbc1d_3''}
u_0^{n+2} \, = \, -\sum_{0 \le m \le (n+1)/2} \, s_m^0 \, u_1^{n+1-2\, m} \, .
\end{equation}

In view of understanding the stability properties of the numerical boundary condition \eqref{eq:dtbc1d_3''}, it might be important to determine the 
asymptotics of the sequence $(s_n^0)_{n \in \mathbb{N}}$. Let us recall that $s_n^0$ equals $(P_{n-1}(\alpha_x)-P_{n+1}(\alpha_x))/((4\, n+2) \, \mu_x)$ 
for $n \geq 2$ (here $\mu_x \in (0,1)$ and $\alpha_x=1-2\, \mu_x^2$). We use the so-called Laplace formula for the Legendre polynomials, see 
\cite[Theorem 8.21.2]{szego}:
\begin{equation}
\label{laplace}
\forall \, \theta \in (0,\pi) \, ,\quad P_n(\cos \theta) \, = \, \sqrt{\dfrac{2}{\pi \, n \, \sin \theta}} \, \cos \left( \big( n+\dfrac{1}{2} \big) \, \theta -\dfrac{\pi}{4} \right) 
+O(n^{-3/2}) \, ,
\end{equation}
where the remainder term is even uniform with respect to $\theta$ on every compact set of the form $[\varepsilon,\pi-\varepsilon]$, $\varepsilon>0$. 
We fix the angle $\theta_x \in (0,\pi)$ such that $\cos \theta_x =\alpha_x=1-2\, \mu_x^2$, and therefore $\sin \theta_x =2\, \mu_x \, \sqrt{1-\mu_x^2}$. 
We then apply the Laplace formula and obtain after a few simplifications:
\begin{equation}
\label{asymptotsn0}
s_n^0 \, = \, \dfrac{(1-\mu_x^2)^{1/4}}{\sqrt{\pi \, \mu_x \, n^3}} \, \sin \left( \big( n+\dfrac{1}{2} \big) \, \theta_x -\dfrac{\pi}{4} \right) +O(n^{-5/2}) \, .
\end{equation}
This asymptotic behavior can be verified on numerical experiments.

\paragraph{Determining inductively the expansion.} As observed in \cite{BNS,jfc}, the Laurent series expansion of $\kappa_s^0$ can also be computed 
inductively by using the equation \eqref{eq:rec1d}. The main benefit is that the method does not rely on an explicit knowledge of a power series 
expansion such as the one involving the Legendre polynomials (which is useful only for three point schemes). The methodology below is therefore 
more flexible in view of being generalized to numerical schemes with larger stencils.

We know that $\kappa_s^0$ tends to zero at infinity so we may plug its Laurent series expansion\footnote{It is proved in \cite{jfc} that the Laurent series 
is convergent for $|z|>1$ which follows from the splitting of the two roots of \eqref{eq:rec1d} for $|z|>1$.} $\sum_{n \geq 1} \sigma_n^0 \, z^{-n}$ into 
\eqref{eq:rec1d} and obtain
\[
\big( z-z^{-1} \big) \, \sum_{n=1}^\infty \sigma_n^0 \, z^{-n} +\mu_x \, \left( \left (\sum_{n=1}^\infty \sigma_n^0 \, z^{-n} \right)^2-1 \right) \, = \, 0 \, ,
\]
which is found to be equivalent to $\sigma_1^0=\mu_x$ and
\begin{equation}
\label{inductionsigman0}
\forall \, n \in \mathbb{N} \, ,\quad \sigma_{n+2}^0 = \sigma_n^0 -\mu_x \, \sum_{m=1}^n \, \sigma_m^0 \, \sigma_{n+1-m}^0 \, ,
\end{equation}
where we use the convention $\sigma_0^0=0$. It is not difficult to see that if $n$ is even, then $\sigma_n^0=0$. Let now $n=2 \, p+1$ be odd, $p\geq 0$. 
An easy computation gives
\[
\sigma_{2 \, p+3}^0 \, = \, \sigma_{2\, p+1}^0 -\mu_x \, \sum_{m=0}^p \sigma_{2\, m+1}^0 \, \sigma_{2(p-m)+1}^0 \, .
\]
This means that the sequence $(s_n^0)_{n \in \mathbb{N}}$ appearing in \eqref{eq:dtbc1d_3'} and that satisfies the recurrence relation \eqref{eq:sn} 
also satisfies
\begin{equation}
\label{recurrencesn0}
s_{n+1}^0 \, = \, s_n^0 -\mu_x \, \sum_{p=0}^n s_p^0 \, s_{n-p}^0 \, ,
\end{equation}
with $s_0^0 =\mu_x$ ($s_n^0$ coincides with $\sigma_{2\, n+1}^0$). Observe that the value $s_1^0 =\mu_x \, (1-\mu_x^2)$ that we have found in the 
previous paragraph is consistent with \eqref{recurrencesn0}. It seems less clear to derive the asymptotic behavior \eqref{asymptotsn0} by starting from 
\eqref{recurrencesn0} rather than from the more explicit representation \eqref{eq:dtbc1d_2}. Both approaches thus seem to have their own interest.
\bigskip

Given the initial conditions $(u^0_j)_{0 \leq j \leq J+1}$ and $(u^1_j)_{0 \leq j \leq J+1}$ as described above ($u^0$ is determined by the initial condition 
for the transport equation and $u^1$ is determined by applying the Lax-Wendroff method), both of which satisfy $u_0^0=u_0^1=u_{J+1}^0=u_{J+1}^1=0$, 
the leap-frog scheme on the interval $\{ 1,\dots,J \}$ with DTBC reads:
\begin{subequations}
\label{LF1d-DTBC}
\begin{align}
u_j^{n+2} &= \, u_j^n-\mu_x \, \left (u_{j+1}^{n+1}-u_{j-1}^{n+1}\right) \, ,\quad \quad 
1 \leq j \leq J \, ,\quad n \in \mathbb{N} \, ,\label{LF1d-DTBC-a} \\
u_0^{n+2} &= \, -\sum_{0 \le m \le (n+1)/2} \, s_m^0 \, u_1^{n+1-2\, m} \, ,\label{LF1d-DTBC-b} \\
u_{J+1}^{n+2} &= \, \sum_{0 \le m \le (n+1)/2} \, s_m^0 \, u_J^{n+1-2\, m} \, ,\label{LF1d-DTBC-c}
\end{align}
\end{subequations}
where the parameter $\mu_x$ in \eqref{LF1d-DTBC} is fixed such that $0 \le \mu_x<1$ and the sequence $(s_n^0)_{n \in \mathbb{N}}$ is defined by 
$s_0^0 =\mu_x$, $s_1^0=\mu_x \, (1-\mu_x^2)$ and the induction relation \eqref{eq:sn} for $n \ge 2$ (or equivalently \eqref{recurrencesn0}, but this 
is slightly less efficient from a numerical point of view). We recall that we consider the case $c_x \ge 0$ in \eqref{eq1}, hence the sign for $\mu_x$, 
but the other case $c_x \le 0$ can be dealt with the same arguments (assuming $-1<\mu_x \le 0$).

\subsection{Stability analysis on a half-line}

In this short paragraph, we explain why the normal mode analysis predicts `neutral stability' for any of the two half-line problems where the 
leap-frog scheme \eqref{LF1d-DTBC-a} is considered on the half-line $\{ j \le J\}$ or $\{ j \ge 1 \}$, in combination with either \eqref{LF1d-DTBC-b} 
or \eqref{LF1d-DTBC-c}. Let us focus on the case $\{ j \ge 1 \}$, since the other case is entirely similar.

The normal mode analysis consists in determining the solutions of the leap-frog scheme \eqref{LF1d-DTBC-a} of the form $u_j^n = z^n \, v_j$, with 
$|z|>1$, and $(v_j) \in \ell^2$. Among such sequences, the final question is to determine whether there exist nonzero ones that satisfy the numerical 
boundary condition \eqref{LF1d-DTBC-b}. Following \cite{gko}, we shall say that the Godunov-Ryabenkii condition is satisfied if there is no nonzero 
sequence $(v_j) \in \ell^2$ such that $u_j^n = z^n \, v_j$ is a solution to the leap-frog scheme that satisfies \eqref{LF1d-DTBC-b}. Otherwise, we shall 
say that $z$ is an unstable eigenvalue. Of course, unstable eigenvalues preclude (in the most violent way) stability estimates hence a desirable feature 
of any numerical boundary condition is to satisfy the Godunov-Ryabenkii condition.

Plugging the ansatz $u_j^n = z^n \, v_j$ in the leap-frog scheme, we find that the sequence $(v_j)$ should belong to $\ell^2$ and satisfy
$$
\forall \, j \ge 1 \, ,\quad \big( z^2-1 \big) \, v_j \, + \, \mu_x \, z \, \big( v_{j+1}-v_{j-1} \big) \, = \, 0 \, .
$$
The computation of the $\ell^2$-solutions to this recurrence relation follows by solving the characteristic equation \eqref{eq:rec1d} and by determining 
among the roots to \eqref{eq:rec1d} which have modulus less than $1$. This classification has already been performed in the previous paragraph, so 
for $|z|>1$, we can conclude that the sequence $(v_j)$ is given by
$$
\forall \, j \ge 1 \, ,\quad v_j \, = \, v_0 \, \kappa_s^0(z)^j \, ,\quad v_0 \in \mathbb{C} \, .
$$
Inserting into \eqref{LF1d-DTBC-b}, we have to determine whether among such sequences (that are parametrized by their initial state $v_0$), 
we can have
$$
v_0 \, = \, -\sum_{0 \le m \le (n+1)/2} \, s_m^0 \, z^{-2\, m-1} \, v_1 \, .
$$
Taking the limit $n \to \infty$, we thus need to determine whether there can hold
$$
v_0 \, = -\kappa_s^0 (z) \, v_1 \, = \, \dfrac{1}{\kappa_u^0(z)} \, v_1 \, = \, \dfrac{\kappa_s^0(z)}{\kappa_u^0(z)} \, v_0 \, ,
$$
where the last equality follows from the expression of $v_1$. In other words, we need to determine whether there can hold $\kappa_s^0(z) 
=\kappa_u^0(z)$ for some complex number $z$ with $|z|>1$. This is clearly impossible since $\kappa_s^0(z)$ has modulus $<1$ and $\kappa_u^0(z)$ 
has modulus $>1$. This means that the so-called Godunov-Ryabenkii condition holds (see \cite{gko}): there does not exist any unstable eigenvalue 
for the half-space problem $\{ j \ge 1 \}$ with the numerical boundary condition \eqref{LF1d-DTBC-b}. Let us observe however that the roots 
$\kappa_s^0$ and $\kappa_u^0$ coincide when $z$ equals one of the four values
$$
\pm i \, \mu_x \pm \sqrt{1-\mu_x^2} \, ,
$$
for which the characteristic equation \eqref{eq:rec1d} has a double root. These values of $z$ correspond to the \emph{glancing} spatial frequencies 
$\kappa$ (for which the associated group velocity vanishes). This is a neutral stability case, whose continuous counterpart has been analyzed in 
details in \cite[chapter 7]{benzoni-serre}, see also \cite{trefethen3} for other examples of this situation.

\subsection{Fast implementation of approximate DTBC with sums of exponentials}

The computation of the DTBC \eqref{LF1d-DTBC-b} and \eqref{LF1d-DTBC-c} at nodes $x_\ell$ and $x_r$ and time $t^{n+2}=(n+2) \, \delta t$ 
requires $n/2$ sums and $n/2$ multiplications. Therefore, the total cost of discrete convolution computations for a full simulation up to time 
$T_f=N_f \, \delta t$ is $\mathcal{O}(N_f^2)$. Since $\delta t$ is related to $\delta x$ by the CFL condition \eqref{cfl1D}, a fine space grid will 
make the time step small and as a consequence the number of time steps $N_f$ large. Since the leap-frog scheme \eqref{LF1d} is explicit, 
each time iteration requires $\mathcal{O}(J)$ operations for interior nodes $x_j$, $j=1,\cdots,J$, which gives $\mathcal{O}(JN_f)$ operations 
for a full simulation. Thus, the total number of operations is $\mathcal{O}(N_f^2+N_fJ)$. If $N_f>J$, which corresponds to \emph{large time} 
simulations, the main part of the computational cost is related to the computation of discrete convolutions. Moreover, it is necessary to store 
in memory the evolution of the solution at the two boundary nodes. It may therefore be useful to reduce this part of the computation by applying 
a more `local' formula.

In \cite{DEDNER2001448,AES}, a fast convolution procedure to compute approximation of discrete convolutions $C_n(v):=\sum_{k=0}^n v_k \, 
\nu_{n-k}$ was introduced. We assume that the convolution coefficients $(\nu_k)_{k \in \mathbb{N}}$ are given and satisfy a decay property similar 
to \eqref{asymptotsn0} that allows to derive their approximation by ``sums of exponentials''. In order to present the idea and the efficiency of the 
method, let us assume for a moment that we are able to compute an approximation $(\tilde{\nu}_k)_{k \in \mathbb{N}}$ of the convolution terms 
$(\nu_k)_{k \in \mathbb{N}}$ given by
\begin{equation}
\label{eq:appr_coef}
\nu_k \, \approx \, \tilde{\nu}_k \, = \, \sum_{m=1}^M b_m \, q_m^{-k} \, ,\quad k\geq 0,
\end{equation}
where the integer $M$, the coefficients $b_m$ and the complex numbers $q_m$, all of them satisfying $|q_m|>1$, have to be defined. Denoting 
$\tilde{C}_n(v) := \sum_{k=0}^n v_k \, \tilde{\nu}_{n-k}$, we compute
\[
\tilde{C}_n(v) \, = \, \sum_{k=0}^n v_k \, \sum_{m=1}^M b_m \, q_m^{-(n-k)} \, = \, \sum_{m=1}^M b_m \, \sum_{k=0}^n v_k \, q_m^{-(n-k)} \, =: \, 
\sum_{m=1}^M C^{(n)}_m(v) \, .
\]
The efficiency of the method relies on the derivation of a recurrence relation to compute the terms $C^{(n)}_m(v)$:
\begin{equation*}
C^{(n)}_m(v) \, = \, b_m \, \sum_{k=0}^{n-1} v_k \, q_m^{-(n-k)} +b_m \, v_n \, = \, b_m \, q_m^{-1} \sum_{k=0}^{n-1} v_k \, q_m^{-(n-1-k)} 
+b_m \, v_n \, ,
\end{equation*}
so we get the recurrence relation
\begin{equation}
\label{eq:expl_coef}
C^{(n)}_m(v) \, = \, q_m^{-1} \, C^{(n-1)}_m(v) +b_m \, v_n \, .
\end{equation}
Assuming $v_0=0$, the initial value for each $C^{(n)}_m(v)$ is given by $C^{(0)}_m(v) =0$. It should be understood that in a practical implementation, 
each value $C^{(n)}_m(v)$ is stored in memory so the recurrence \eqref{eq:expl_coef} represents only three operations at each time iteration.

Accordingly, the original discrete convolution $C_n(v)$ is approximated by
\begin{equation}
\label{eq:dtbc_sumexp}  
C_n(v) \, = \, \sum_{k=0}^n v_k \, \nu_{n-k} \, \approx \, \tilde{C}_n(v) \, = \, \sum_{m=1}^M C^{(n)}_m(v) \, , 
\end{equation}
and we have replace the computation of a nonlocal discrete convolution that involves $n$ operations by the sum of $M$ terms computed by the explicit 
local operations \eqref{eq:expl_coef}. In this case, the total number of operations is $\mathcal{O}(N_f(J+M))$.

It remains to derive the approximation \eqref{eq:appr_coef} of the sequence $(\nu_k)_{k \in \mathbb{N}}$ by $(\tilde{\nu}_k)_{k \in \mathbb{N}}$. The 
sum-of-exponential formula \eqref{eq:appr_coef} is frequent and various techniques are available to compute the coefficients $(b_m)_m$ and $(q_m)_m$ 
(see for example \cite{BEYLKIN200517,BEYLKIN2010131}). We present here a modified version of the algorithm derived in \cite{AES}. Let us consider 
the power series
\[
f(x) \, := \, \sum_{k \in \mathbb{N}} \nu_k \, x^k \, ,\quad \text{for } |x| < 1 \, ,
\]
associated with the original sequence $(\nu_k)_{k \in \mathbb{N}}$ (which we assume to satisfy a polynomial bound of the form $|\nu_k| \lesssim k^\alpha$, 
$\alpha \in \mathbb{R}$, so the power series converges on the unit disk). We then consider its $[N,M]$ Pad\'e approximant, see \cite{baker,bakergravesmorris}, 
which we denote $\tilde{f}(x)=P_N(x)/Q_M(x)$, where $P_N$ and $Q_M$ are respectively polynomials of degree $N$ and $M$. We assume $Q_M$ to be 
unitary in order to normalize the polynomials. Here we shall always consider the case $N<M$. The power series expansion of $\tilde{f}$ at $0$ is given by
\[
\tilde{f}(x) \, = \, \sum_{k \in \mathbb{N}} \tilde{\nu}_k \, x^k \, , 
\]
with $\tilde{\nu}_k=\nu_k$ for $0\leq k\leq N+M$ (by definition of the Pad\'e approximants).

Let us now assume that $Q_M$ has $M$ simple roots $q_m$ with $|q_m|>1$ for $1\leq m \leq M$. So, we have
\[
Q_M(x) \, = \, \prod_{m=1}^M (x-q_m) \, ,\quad \text{and} \quad 
Q_M'(x) \, = \, \sum_{k=1}^M \, \prod_{m \neq k} (x-q_m) \, ,
\]
and thus
\[
Q_M'(q_m)=\prod_{k\neq m} (q_m-q_k) \, .
\]
Let us now define the coefficients
\begin{equation}
\label{eq:bm}
b_m \, := \, - \, \dfrac{P_{N}(q_m)}{q_m \, Q'_M(q_m)} \, ,\quad 1 \leq m \leq M \, .
\end{equation}
So, we obtain
\begin{align*}
\sum_{m=1}^M \dfrac{b_m \, q_m}{q_m-x} \, = \, \sum_{m=1}^M \dfrac{P_{N}(q_m)}{Q'_M(q_m)} \, \dfrac{1}{x-q_m} 
\, &= \, \sum_{m=1}^M \frac{P_N(q_m)}{Q'_M(q_m)} \, \dfrac{\prod_{k \neq m}(x-q_k)}{\prod_{k=1}^M (x-q_k)} \\
&= \, \sum_{m=1}^M \dfrac{P_{N}(q_m)}{Q_M(x)} \, \dfrac{\prod_{k\neq m}(x-q_k)}{Q'_M(q_m)} \, .
\end{align*}
Thereby,
\[
\sum_{m=1}^M \dfrac{b_m \, q_m}{q_m-x} \, = \, 
\dfrac{1}{Q_M(x)} \, \sum_{m=1}^M P_N(q_m) \, \prod_{k\neq m} \dfrac{x-q_k}{q_m-q_k} \, =: \, \dfrac{R_{M-1}(x)}{Q_M(x)} \, .
\]
We are going to show that the polynomial $R_{M-1}$ in the latter equality equals $P_N$. Indeed, we have $R_{M-1}(q_m)=P_{N}(q_m)$ 
for any $1\leq m \leq M$ and $d^\circ(R_{M-1}) \le M-1$. Since $P_N$ is of degree $N \leq M-1$, by uniqueness of the interpolating polynomial, 
we have $R_{M-1}=P_N$. As a partial conclusion, we have
\[
\tilde{f}(x) \, = \, \dfrac{P_N(x)}{Q_M(x)} \, = \, \sum_{m=1}^M \dfrac{b_m \, q_m}{q_m-x} \, ,
\]
where the coefficients $b_m$ are defined in \eqref{eq:bm}. Writing $q_m-x=q_m(1-x/q_m)$, and using $|q_m|>1$, then for $|x| \le 1<|q_m|$, we obtain
\[
\dfrac{1}{q_m-x} \, = \, q_m^{-1} \, \sum_{k=0}^\infty \left( \dfrac{x}{q_m} \right)^k \, .
\]
So,
\begin{equation*}
\tilde{f}(x) \, = \, \sum_{m=1}^M b_m \, \sum_{k=0}^\infty \left( \dfrac{x}{q_m} \right)^k \, = \, 
\sum_{k=0}^\infty \, \sum_{m=1}^M b_m \, q_m^{-k} \, x^k \, .
\end{equation*}
Let us recall that the power series expansion of $\tilde{f}$ at the origin reads $\tilde{f}(x)=\sum_{k\geq 0} \tilde{\nu}_{k} \, x^k$ so proceeding by 
identification, we have
\[
\forall \, k \ge 0 \, ,\quad \tilde{\nu}_k \, = \, \sum_{m=1}^M b_m \, q_m^{-k} \, ,
\]
and the coefficients $\tilde{\nu}_k$ satisfy the property
\[
\tilde{\nu}_k \, = \, \nu_k \quad \text{for } 0 \leq k\leq M+N \, .
\]

The computation of $(\tilde{\nu}_k)_{k \in \mathbb{N}}$ requires to determine the $[N,M]$ Pad\'e approximant of $f$, the roots $q_m$ of $Q_M$ 
(praying for them to be simple) and the coefficients $b_m$ given by \eqref{eq:bm}.

\begin{remark}
In \cite{AES}, the authors suggested to take $N=M-1$. Unfortunately, when we apply the latter approximation strategy to the function $\kappa_s^0(1/x)$, 
we are not always able to verify the assumption $\inf_m |q_m|>1$ on the roots of $Q_M$, even when these roots are computed with (very) high accuracy. 
Our choice to take $N$ `free' and not necessarily equal to $M-1$ allows us to find Pad\'e approximants for which we can verify $\inf_m |q_m|>1$. We 
present below some experiments to emphasize this remark.
\end{remark}

\begin{remark}
In \cite{AES}, the authors justified the convergence of the approximate coefficients $\tilde{\nu}_k$ to ${\nu}_k$ as $M \to \infty$. Their result was based on 
the so-called Baker-Gammel-Wills conjecture, which is unfortunately now known to be false, see \cite{baker-ctrex} and references therein. That does not 
mean however that the above procedure is meaningless, as the numerical results in \cite{AES} and those we present below confirm. The main difficulty 
in applying the above approach is that there does not seem to be any general result available to make sure that the condition $\inf_m |q_m|>1$ holds in 
a systematic way (not mentioning the condition that the roots $q_m$ should be simple...).
\end{remark}

In \cite{AES}, a Maple code is proposed to calculate the Pad\'e approximant of $f$. The use of Maple insures to not be limited to double precision computations. 
This is very important to compute accurately the coefficients of the monomials of $P_N$ and $Q_M$, and accordingly to determine the roots of $Q_M$ with high 
accuracy. This process is however quite limited since it relies on formal calculus softwares and is also restricted to small values of $M$ (typically $M\leq 20$). We 
therefore present below an alternative self-contained procedure which allows us to consider much larger values of $M$ (say, $M \le 200$) provided that $N$ is not 
close to $M$ (say $N \le 50$). When $N$ is close to $M$, the coefficients of $P_N$ and $Q_M$ become extremely large and the numerical computations become 
rather unstable. The $[N,M]$ Pad\'e approximant $\tilde{f}$ of $f$ reads
\[
\tilde{f}(x) \, = \, \dfrac{P_N(x)}{Q_M(x)} \, = \, \dfrac{\sum_{j=0}^N \mathfrak{p}_j \, x^j}{\sum_{j=0}^M \mathfrak{q}_j \, x^j} \, ,
\]
where we fix for instance $\mathfrak{q}_0=1$, and we restrict to the case $N<M$. (The normalization convention for $Q_M$ is not the same as above but this 
does not affect the previous arguments.) We therefore have to compute $N+M+1$ unknowns coefficients, which are solutions of the $N+M+1$ equations given 
by
\begin{equation}
\label{eqpade}
D^{(n)} \, \left[ Q_M(x) \, f(x) \right]_{|_{x=0}} \, = \, 
D^{(n)} \, \left[ P_N(x) \right]_{|_{x=0}} \, , \quad 0\leq n\leq N+M \, ,
\end{equation}
where $D^{(n)}$ denotes the $n^{\text{th}}$ derivative with respect to $x$. Indeed, thanks to the property $Q_M(0) \neq 0$, \eqref{eqpade} can be seen to be 
equivalent to the property
$$
D^{(n)} \, \left[ f(x) \right]_{|_{x=0}} \, = \, 
D^{(n)} \, \left[ \dfrac{P_N(x)}{Q_M(x)} \right]_{|_{x=0}} \, , \quad 0\leq n\leq N+M \, ,
$$
which defines the $[N,M]$ Pad\'e approximant of $f$. The equations \eqref{eqpade} reduce to
\begin{equation}
\label{systlin-1}
\sum_{k=1}^n \mathfrak{q}_k \, \nu_{n-k} -\mathfrak{p}_n \, = \, -\nu_n \, ,\quad \text{for } 0 \leq n \leq N \, ,
\end{equation}
and
\begin{equation}
\label{systlin-2}
\sum_{k=1}^{\min (n,M)} \mathfrak{q}_k \, \nu_{n-k}  \, = \, -\nu_n \, ,\quad \text{for } N<n\leq N+M \, .
\end{equation}
The first equation in \eqref{systlin-1} leads to $\mathfrak{p}_0 =\nu_0 =f(0)$, and the remaining equations in \eqref{systlin-1}, \eqref{systlin-2} can be recast as 
a linear system:
\begin{equation}
\label{systlin}
S \, X \, = \, Y \, ,
\end{equation}
where $S\in \mathbb{M}_{N+M,N+M}(\mathbb{R})$ and $(X,Y)\in \mathbb{R}^{N+M}$. The matrix $S$ in \eqref{systlin} is given by (the first row and 
column below present the indices for the sake of clarity):
\[
\begin{array}{c|cccc|cccc|cccc|}
 & 1 & 2 & & N & N+1 & & & M & M+1 & & & N+M \\ \hline
1 & \nu_0 & & & & & & & & -1 & & & \\
2 & \nu_1 & \ddots & & & & & & & & \ddots & & \\
\vdots & \vdots & \ddots & \ddots & & & & & & & & \ddots & \\
N & \nu_{N-1} & \cdots & \nu_1 & \nu_0 & & & & & & & & -1 \\ \hline
N+1 & \nu_N & \cdots & \cdots & \nu_1 & \nu_0 & & & & 0 & & & \\
\vdots & \vdots & & & \vdots & \nu_1 & \ddots & & & & \ddots & & \\
\vdots & \vdots & & & \vdots & \vdots & \ddots & \ddots & & & & \ddots & \\
M & \nu_{M-1} & \cdots & \cdots & \nu_{M-N} & \nu_{M-N+1} & \cdots & \nu_1 & \nu_0 & & & & 0 \\ \hline
M+1 & \nu_M & \cdots & \cdots & \nu_{M-N+1} & \nu_{M-N} & \cdots & \cdots & \nu_1 & 0 & & & \\
\vdots & \vdots & & & \vdots & \vdots & & & \vdots & & \ddots & & \\
\vdots & \vdots & & & \vdots & \vdots & & & \vdots & & & \ddots & \\
N+M & \nu_{N+M-1} & \cdots & \cdots & \nu_M & \nu_{M-1} & \cdots & \cdots & \nu_N & & & & 0 \\ \hline
\end{array} \, .
\]
The unknown in \eqref{systlin} is $X=(\mathfrak{q}_1,\dots,\mathfrak{q}_M,\mathfrak{p}_1,\dots,\mathfrak{p}_N)^T$, and the right hand side is 
$Y=(-\nu_1,-\nu_2,\cdots,-\nu_{N+M})^T$. The linear system \eqref{systlin} thus displays block matrices
\[
\begin{array}{|ccccc|ccc|} \hline
 &&&&&&& \\
 &&&&&&& \\
 &&A&&&&B& \\
 &&&&&&& \\
 &&&&&&& \\ \hline
 &&&&&&& \\
 &&C&&&&0& \\
 &&&&&&& \\ \hline 
 \end{array} \quad \begin{array}{|c|} \hline
 \\
 \\
\mathfrak{q} \\
 \\
 \\ \hline
 \\
\mathfrak{p} \\
 \\ \hline
\end{array} \quad = \quad \begin{array}{|c|} \hline
 \\
 \\
Y_1 \\
 \\
 \\ \hline
 \\
 Y_2 \\
\\ \hline \end{array} \, ,
\]
with $A\in \mathbb{M}_{M,M}(\mathbb{R})$, $B\in \mathbb{M}_{M,N}(\mathbb{R})$ et $C\in \mathbb{M}_{N,M}(\mathbb{R})$. The matrix $A$ 
is a subtriangular Toeplitz matrix. Its inverse (assuming that $\nu_0$ is nonzero) is also a subtriangular Toeplitz matrix whose first column is 
easy to compute by solving the linear system $A \, x=e_1$ where $e_1$ is the first vector of the canonical basis of $\mathbb{R}^M$. The other 
columns of $A^{-1}$ are deduced from the Toeplitz structure.

Using the Schur complement method, we can compute the unknowns by solving successively
\[
C \, A^{-1} \, B \, \mathfrak{p} \, = \, C \, A^{-1} \, Y_1 -Y_2 \, ,
\]
and then
\[
\mathfrak{q} \, = \, A^{-1} \, (Y_1 -B \, \mathfrak{p}) \, .
\]

In order to get high accuracy both in the computations of $\mathfrak{p}$ and $\mathfrak{q}$, but also of the roots $q_m$ of $Q_M$, we use the 
floating-point arithmetic with arbitrary accuracy Python library \verb+mpmath+. The advantage of using the library \verb+mpmath+ is its ability to 
control the accuracy of floating-point arithmetic. For all the numerical experiments in this paper, we set the decimal accuracy to $80$, meaning 
the library uses $266$ bits to hold an approximation of the numbers that is accurate to $80$ decimal places.
\bigskip

\begin{figure}[htbp]
\centering
\begin{tabular}{cc}
\includegraphics[height=.27\textheight]{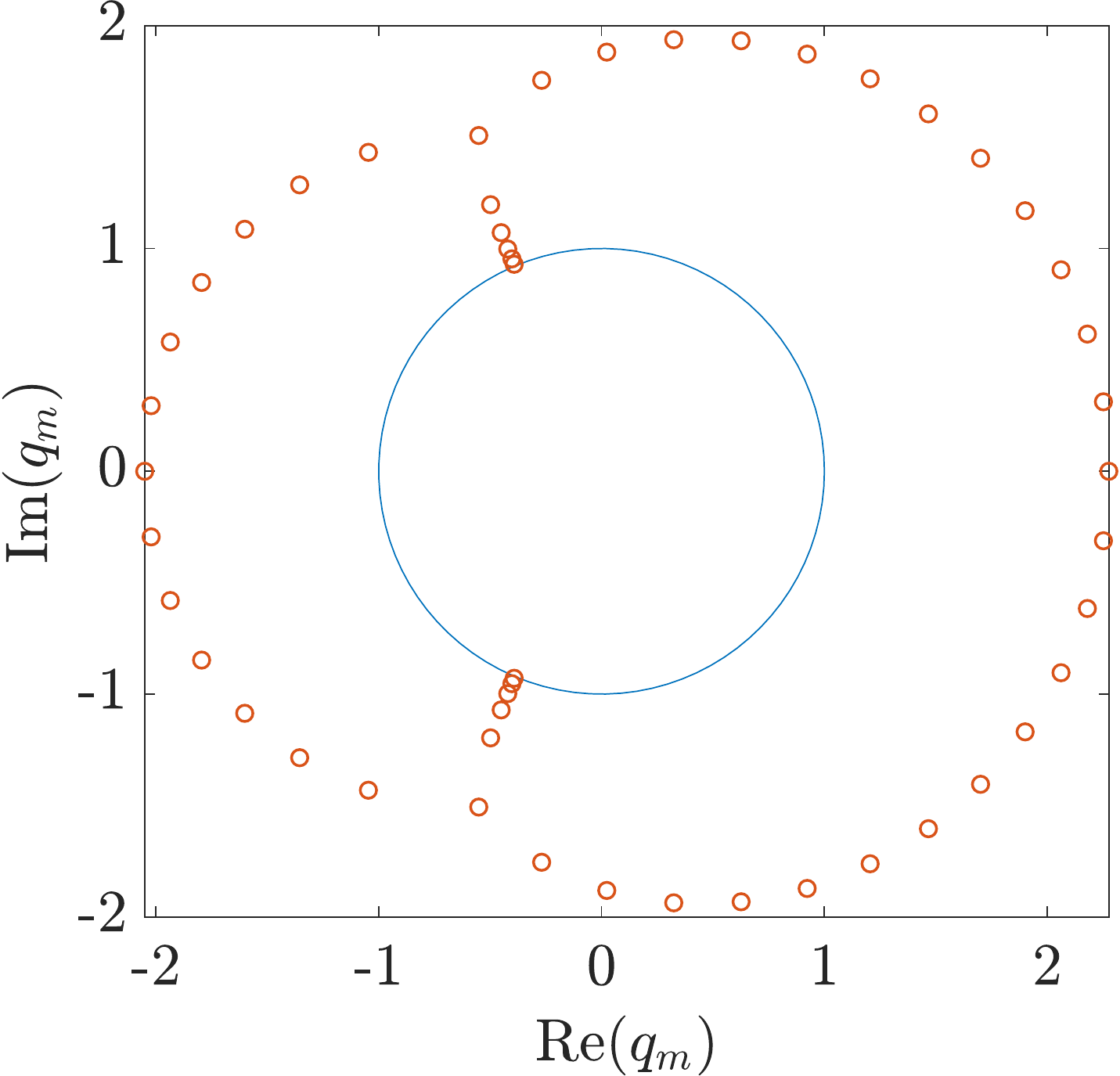} & 
\includegraphics[height=.27\textheight]{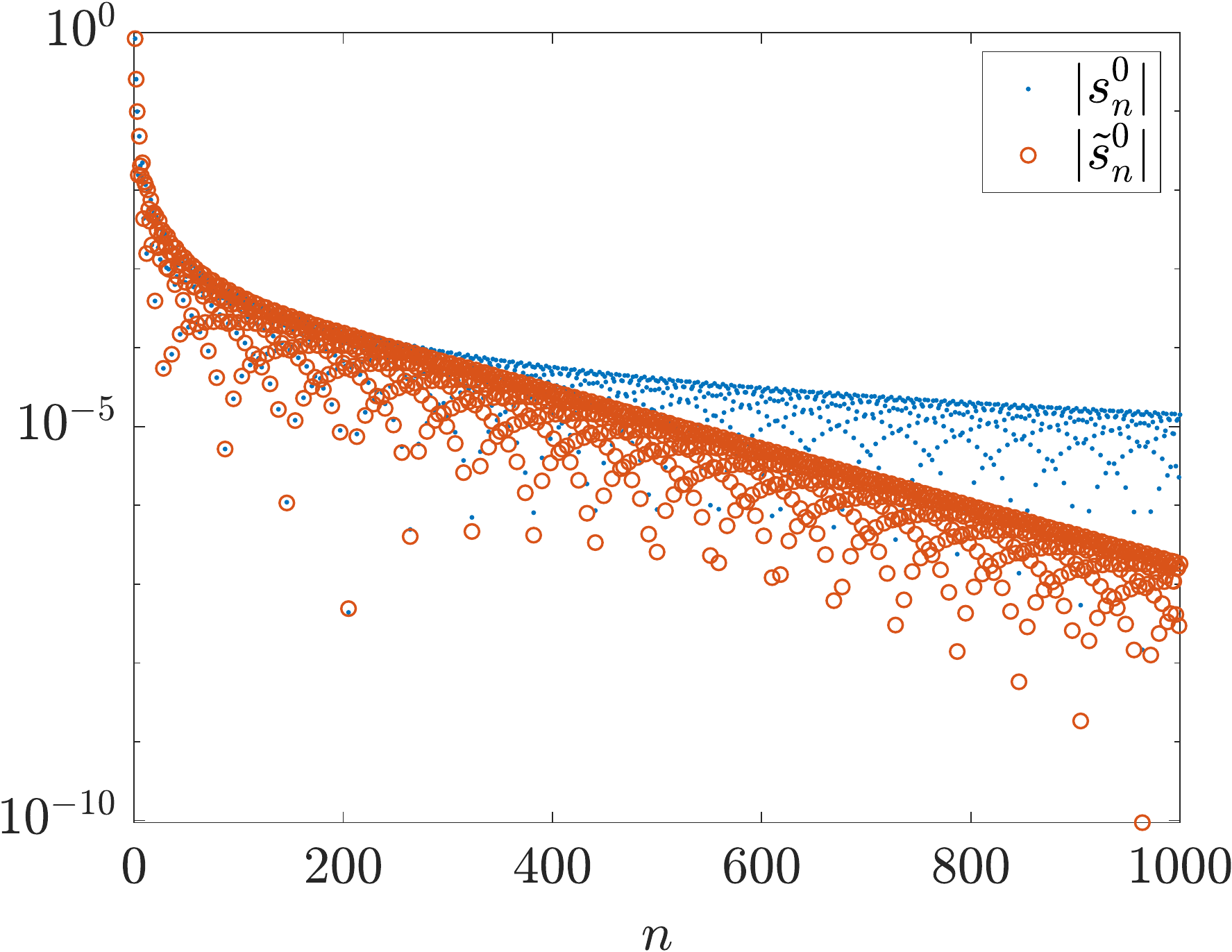} \\
\end{tabular}
\caption{Position of the roots of $Q_M$ (left), and $|s_n^0|$ vs $|\tilde{s}_n^0|$ (right) for $(M,N)=(50,10)$.}
\label{fig:pade01}
\end{figure}

\begin{figure}[htbp]
\centering
\begin{tabular}{cc}
\includegraphics[height=.27\textheight]{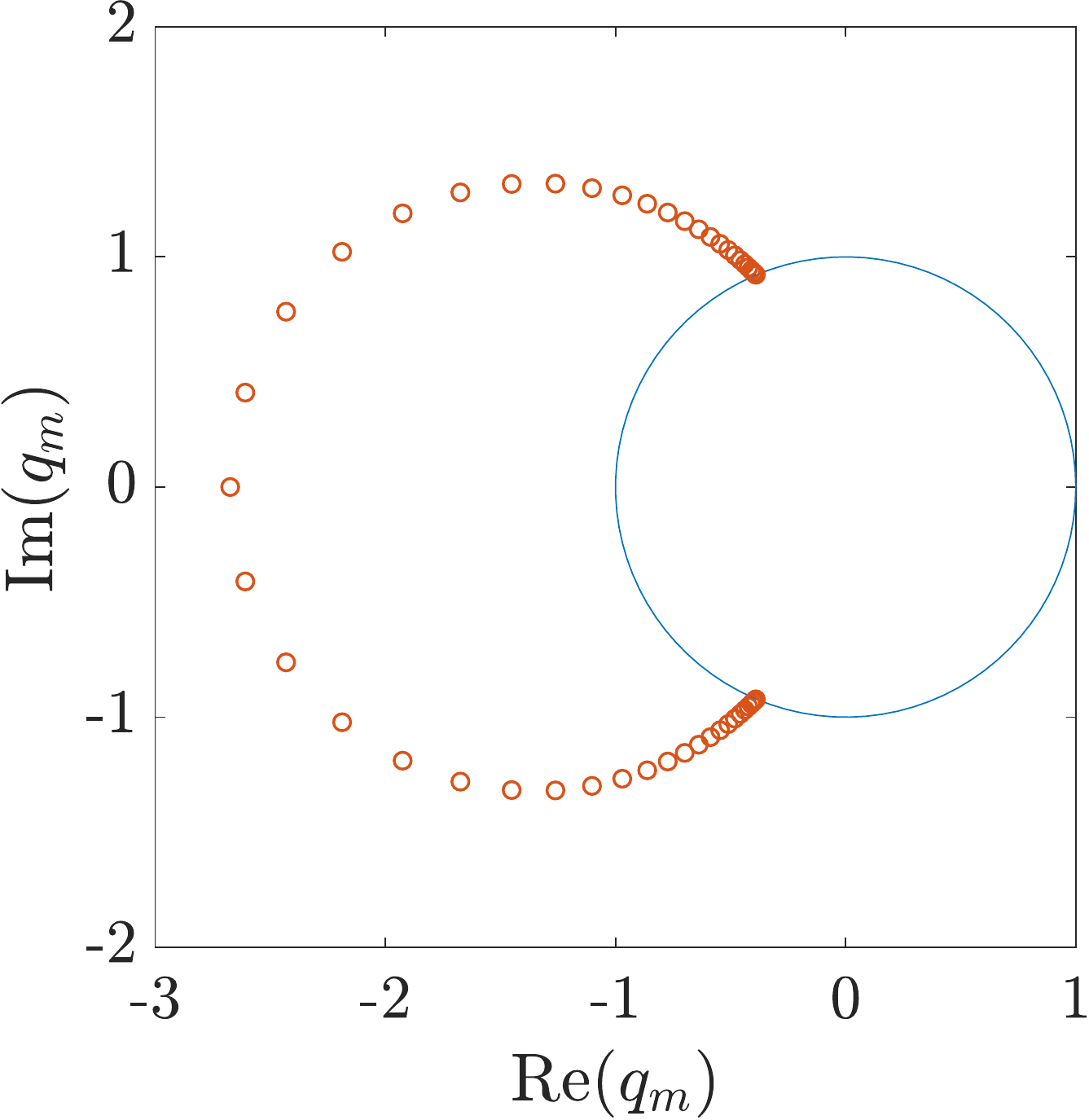} & 
\includegraphics[height=.27\textheight]{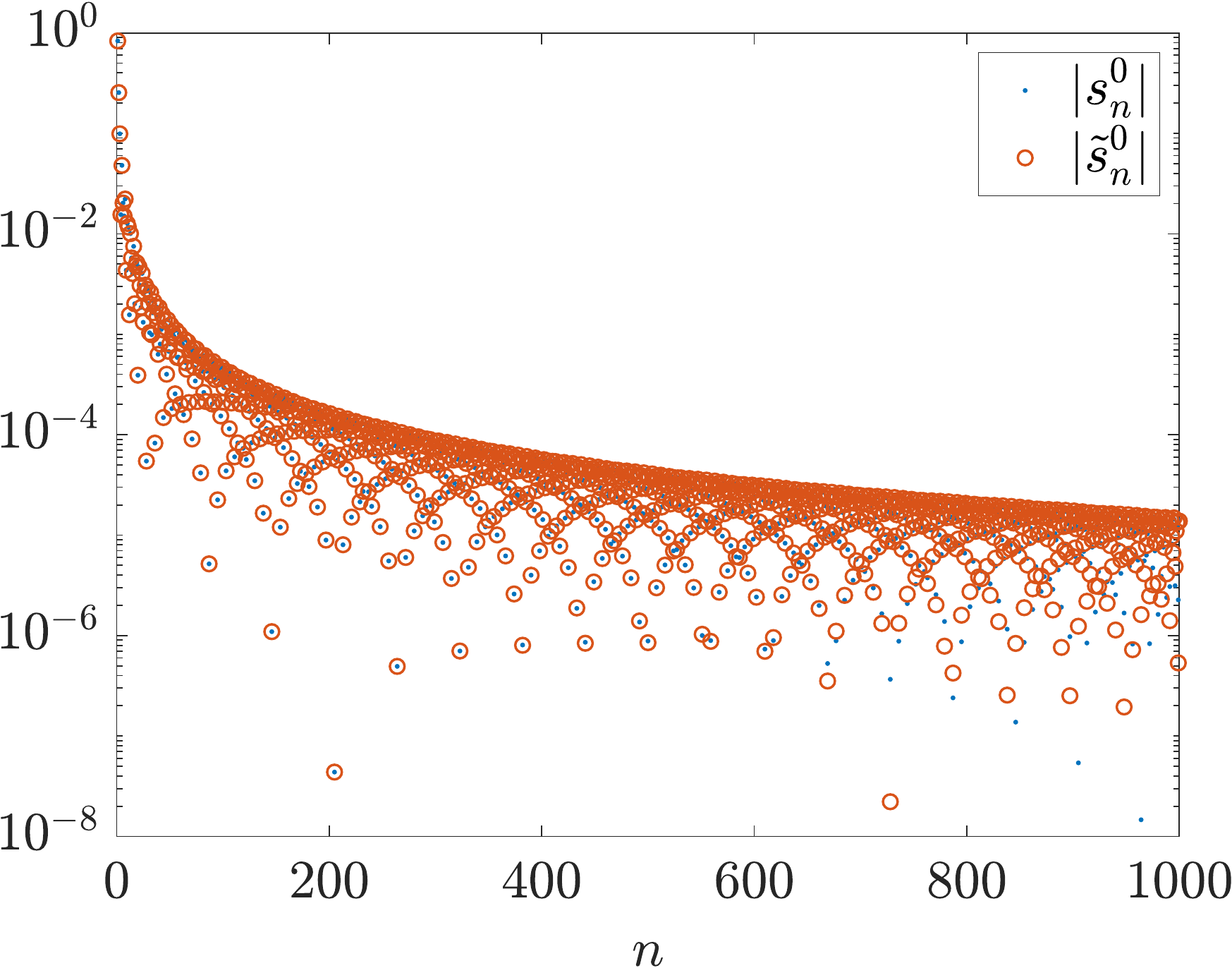} \\
\end{tabular}
\caption{Position of the roots of $Q_M$ (left), and $|s_n^0|$ vs $|\tilde{s}_n^0|$ (right) for $(M,N)=(50,49)$.}
\label{fig:pade02}
\end{figure}

\begin{figure}[htbp]
\centering
\begin{tabular}{cc}
\includegraphics[height=.27\textheight]{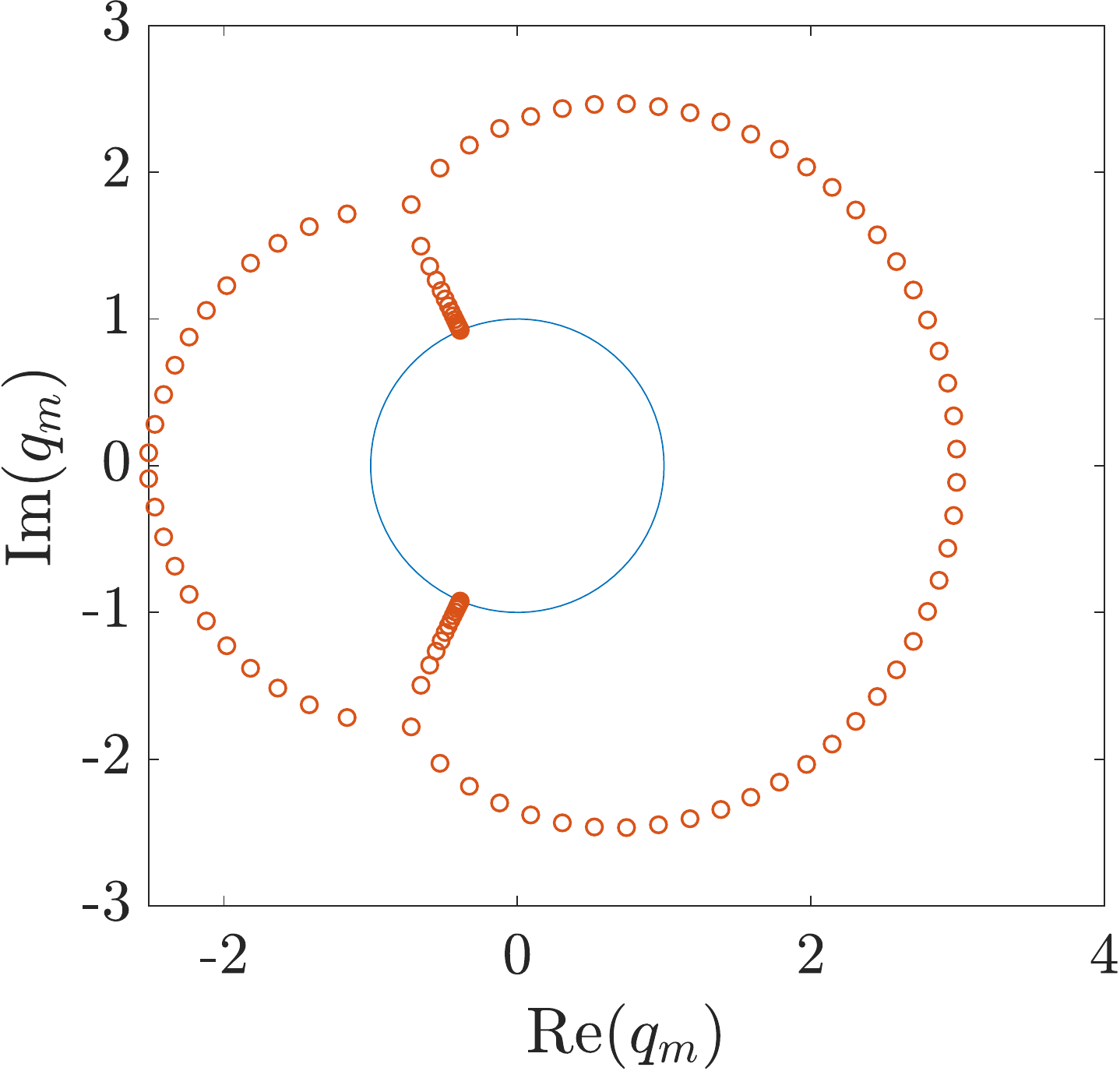} & 
\includegraphics[height=.27\textheight]{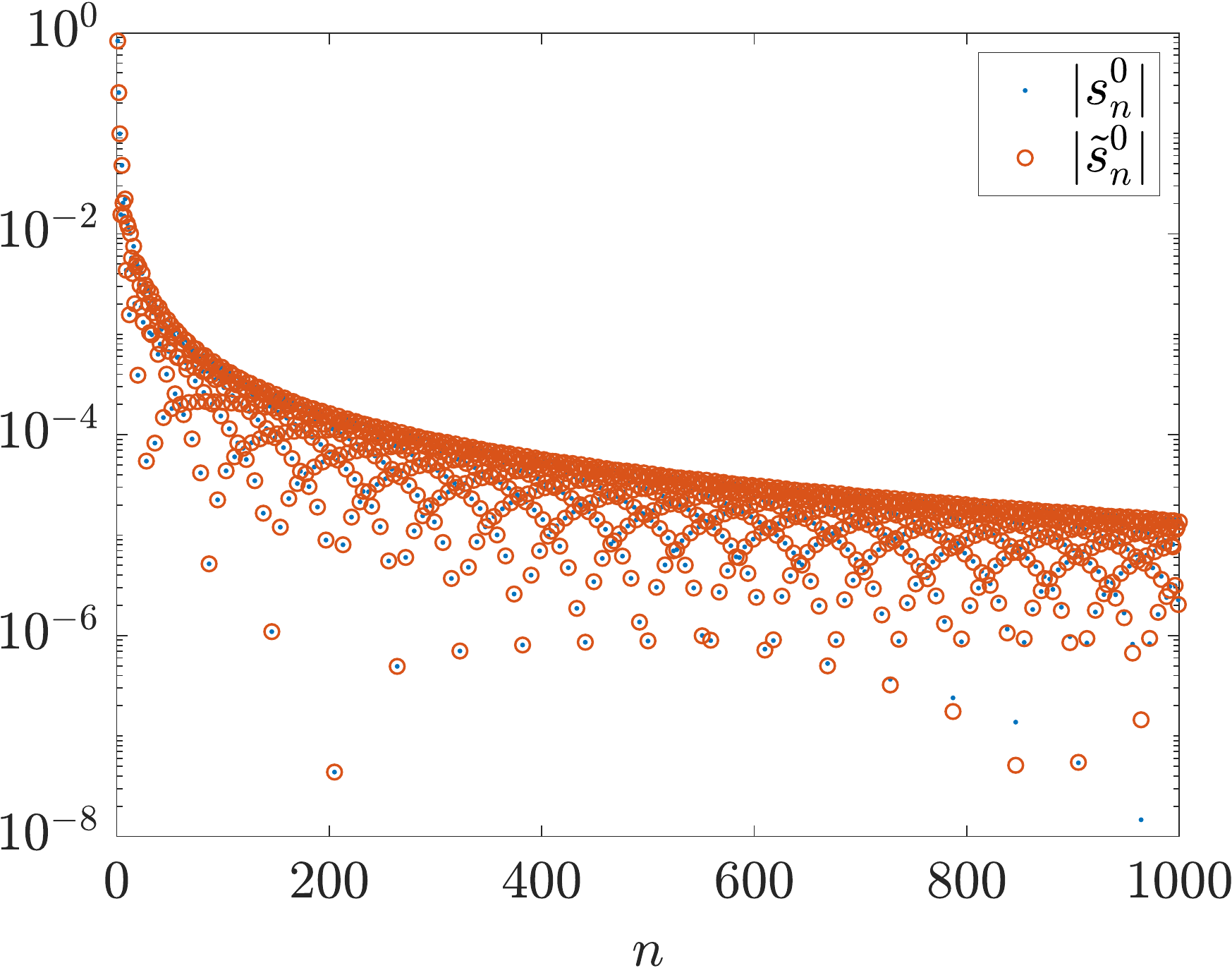} \\
\end{tabular}
\caption{Position of the roots of $Q_M$ (left), and $|s_n^0|$ vs $|\tilde{s}_n^0|$ (right) for $(M,N)=(100,30)$.}
\label{fig:pade03}
\end{figure}

\begin{figure}[htbp]
\centering
\begin{tabular}{cc}
\includegraphics[height=.27\textheight]{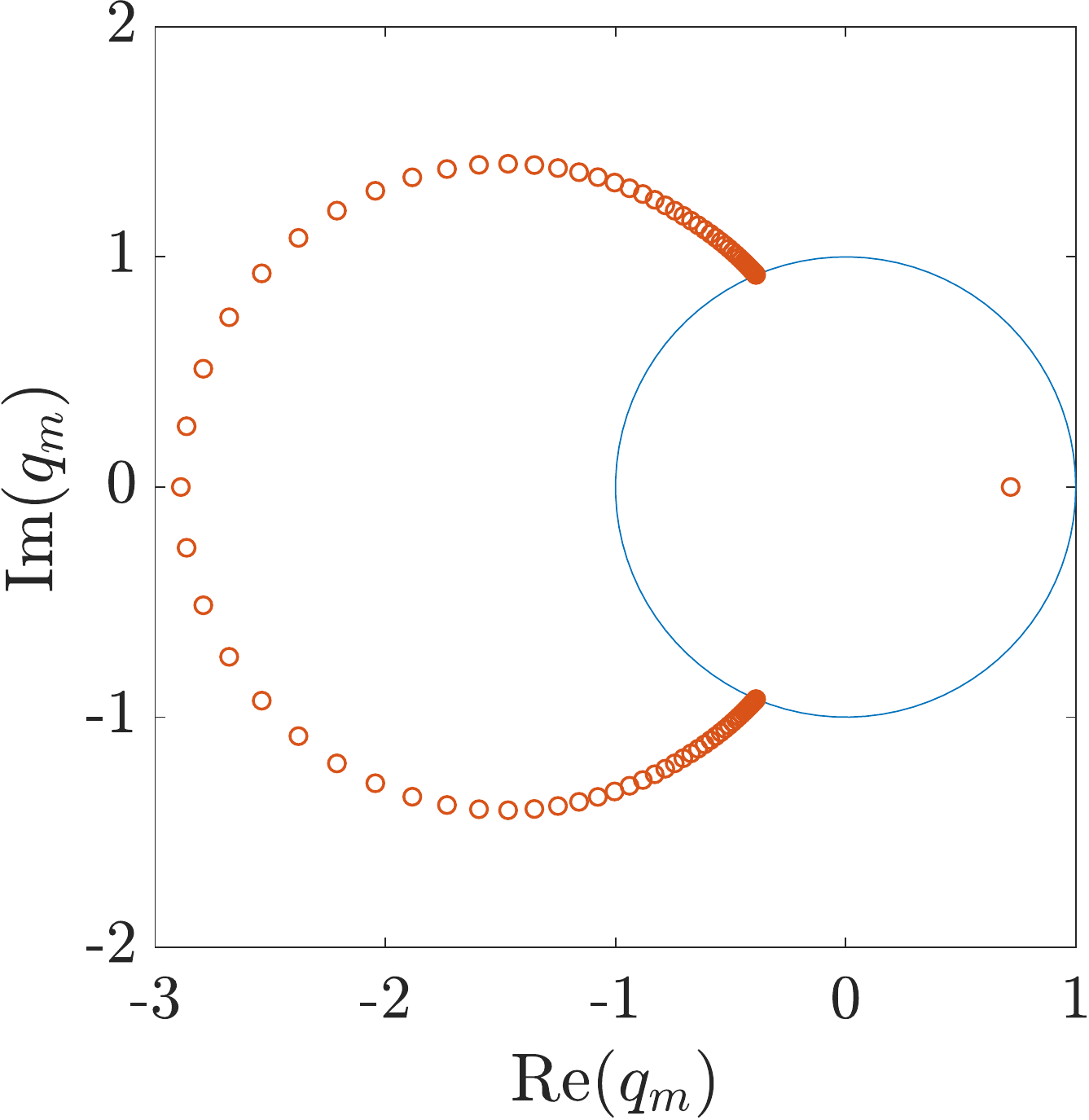} & 
\includegraphics[height=.27\textheight]{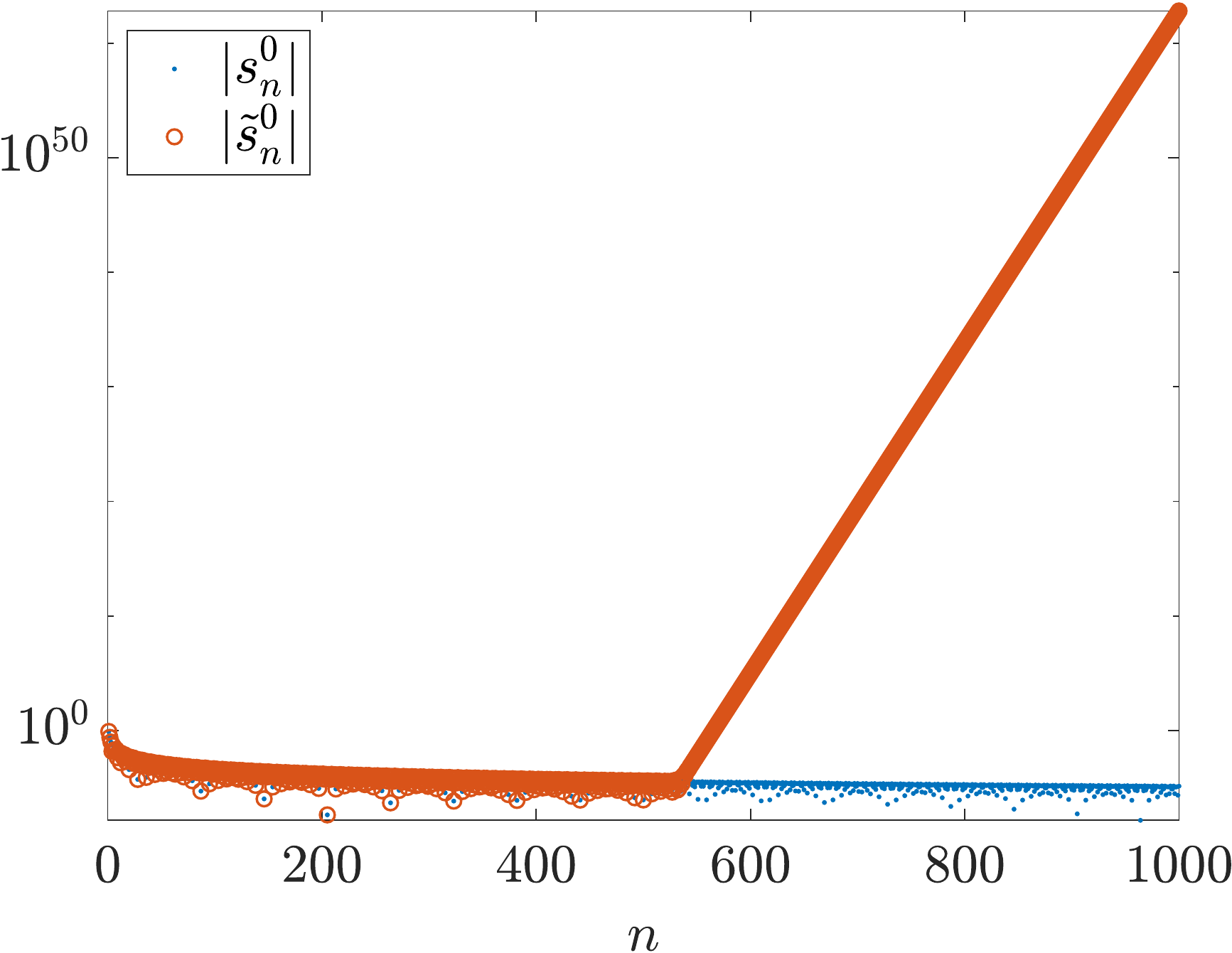} \\
\end{tabular}
\caption{Position of the roots of $Q_M$ (left), and $|s_n^0|$ vs $|\tilde{s}_n^0|$ (right) for $(M,N)=(100,99)$.}
\label{fig:pade04}
\end{figure}

Going back to our numerical scheme \eqref{LF1d-DTBC}, we wish to approximate the sequence $(s_n^0)$ given by \eqref{eq:sn} by a sum of exponentials. 
We apply the above algorithm, which is legitimate since $s_0^0 \neq 0$ ($s_n^0$ plays the above role of $\nu_n$ for all $n$). We present below a comparison 
between the thousand first elements of the sequence $(s_n^0)$ given by \eqref{eq:sn} and their sum of exponentials approximations $(\tilde{s}_n^0)$ for 
various pairs $(N,M)$ with $N<M$. The numerical parameters are chosen to satisfy the CFL condition \eqref{cfl1D} with $\text{cfl}=5/6$. The transport 
velocity is $c=1$. The space and time steps are $\delta x=6/1000$ and $\delta t =5/1000$ (which corresponds to the space interval $[-3,3]$ and choosing 
a thousand grid points). On Figures \ref{fig:pade01} to \ref{fig:pade04}, we plot on the left the position of the roots $(q_m)$ of $Q_M$ in the complex plane, 
compared with the unit circle, and the evolution of $\tilde{s}_n^0$ compared with that of $s_n^0$ 
for $0 \leq n \leq 1000$. We first choose $M=50$ with $N=10$ (see Figure \ref{fig:pade01}), then $M=50$ with $N=49$ (see Figure \ref{fig:pade02}), 
then $M=100$ with $N=30$ (see Figure \ref{fig:pade03}), and eventually $M=100$ with $N=99$ (see Figure \ref{fig:pade04}). On Figures \ref{fig:pade02} 
and \ref{fig:pade04}, one root of the polynomial $Q_M$ is relatively big compared with the other roots and we therefore do not represent it. They 
respectively take the values $-4.6\, 10^{17}$ and $-1.2\, 10^{29}$, which means that some of the coefficients of $Q_M$ are so large that any computation 
involving the roots of $Q_M$ should be taken with much care. We see that the more $N$ is close to $M$, the more $\tilde{s}_n^0$ are close to $s_n^0$ 
for a large interval of values of $n$. The quality of the approximations for $(M,N)=(50,49)$ and $(M,N)=(100,30)$ is equivalent ; the case $(M,N)=(100,30)$ 
is however much more stable since all roots of the polynomial $Q_M$ remain within a `small' ball of the complex plane while $Q_M$ has a very large root 
in the case $(M,N)=(50,49)$. When $(M,N)=(100,99)$, one of the roots that we computed for $Q_M$ has a modulus less than one, which generates blow 
up for sequence $(\tilde{s}_n^0)$ (see Figure \ref{fig:pade04}). This does not prove however that $Q_M$ has a root of modulus less than $1$, because 
the computation is so unstable (due to the very large coefficients in $Q_M$) that the final result should be taken with care even with the decimal accuracy 
set to $80$.

\subsection{Numerical experiments}

We present here some numerical tests with the following parameters:
\[
u_0(x) \, = \, \exp (-10 \, x^2) \, ,\quad x_\ell \, = \, -3 \, ,\quad x_r \, = \, 3,\quad T \, = \, 10 \, ,\quad c \, = \, 1 \, .
\]
The number of grid points is $J+1=1000$, and  the leap-frog scheme \eqref{LF1d} is implemented with a CFL parameter $\mu_x=5/6$ (the 
time step $\delta t$ is computed accordingly). Since the initial condition is roughly concentrated in the interval $[-1,1]$, the exact solution to 
the transport equation more or less vanishes on the interval $[x_\ell,x_r]$ after the time $4$. We measure below the amplitude of the reflected 
wave at the right boundary $x_r$ depending on the choice of numerical boundary conditions.

We present below the evolution of the logarithm $\log_{10} |u_j^n|$ rather than the evolution of the solution $u_j^n$ itself. This gives a much 
more visible representation of wave reflections and their magnitude. Let us observe that in Figures \ref{fig:1d_evol_01} and \ref{fig:1d_evol_02} 
below, we always observe a left going wave emanating from the initial condition, this wave being highly oscillatory (in both space and time) and 
having magnitude $10^{-8}$. This is due to the fact that the leap-frog scheme \eqref{LF1d} supports the wave $(-1)^{j+n}$ which has group 
velocity $-c$ and thus travels backwards \cite{trefethen1}. In Figure \ref{fig:1d_evol_01}, we first compare the leap-frog scheme with exact 
DTBC \eqref{LF1d-DTBC} (left) with the leap-frog scheme implemented with the more classical Neumann type boundary condition 
$u_{J+1}^{n+2}=u_J^{n+1}$, see \cite{trefethen3} (accordingly the numerical boundary condition at $x_\ell$ is $u_0^{n+2}=u_1^{n+1}$). After 
the theoretical exit time for the continuous solution, the DTBC strategy leaves a solution of amplitude $10^{-16}$ in the whole domain while the 
easier (local) Neumann strategy exhibits multiple wave reflections.

We then show on Figure \ref{fig:1d_evol_02} the evolution of the logarithm $\log_{10} |u_j^n|$ when we do not implement the exact DTBC 
\eqref{LF1d-DTBC-b}, \eqref{LF1d-DTBC-c} but rather their approximation by sums of exponentials as explained in the preceding paragraph. 
The numerical simulation with the parameters $(M,N)=(50,6)$ (left of Figure \ref{fig:1d_evol_02}) and $(M,N)=(50,49)$ (right of Figure 
\ref{fig:1d_evol_02}) for the degrees of the Pad\'e approximant. In each of these two cases, the condition that all roots of $Q_M$ are simple 
and lie outside the unit circle has been checked numerically. We observe that in the numerical boundary conditions \eqref{LF1d-DTBC-b} and 
\eqref{LF1d-DTBC-c}, the determination of $u_0^{n+2}$ and $u_{J+1}^{n+2}$ relies on the values $u_1^\sigma$ and $u_J^\sigma$ for either 
only odd or even values of $\sigma$. When we implement the approximate numerical boundary condition \eqref{eq:dtbc_sumexp}, we thus 
need to consider separately the cases where $n$ is either odd or even, and we therefore introduce two sequences $(C^{(2\, p)}_m)$ and 
$(C^{(2\, p+1)}_m)$ for each root $q_m$ of the polynomial $Q_M$. The numerical results presented in Figure \ref{fig:1d_evol_02} exhibit 
very low reflection amplitudes. In the case $(M,N)=(50,49)$, the overall accuracy is comparable with the exact DTBC. Note however that 
the condition $\inf_m |q_m|>1$ for this approximation to make sense is not granted and the implementation of the sum of exponential 
approximation requires first solving for $P_N$ and $Q_M$ and then identifying the roots of $Q_M$ (which is not necessarily easier than 
implementing \eqref{eq:sn} and \eqref{LF1d-DTBC}).

\pgfplotsset{width=7cm}

\pgfplotsset{
     	colormap={matlab}{
     		rgb=(0.000000,0.000000,0.562500)
     		rgb=(0.000000,0.000000,0.625000)
     		rgb=(0.000000,0.000000,0.687500)
     		rgb=(0.000000,0.000000,0.750000)
     		rgb=(0.000000,0.000000,0.812500)
     		rgb=(0.000000,0.000000,0.875000)
     		rgb=(0.000000,0.000000,0.937500)
     		rgb=(0.000000,0.000000,1.000000)
     		rgb=(0.000000,0.062500,1.000000)
     		rgb=(0.000000,0.125000,1.000000)
     		rgb=(0.000000,0.187500,1.000000)
     		rgb=(0.000000,0.250000,1.000000)
     		rgb=(0.000000,0.312500,1.000000)
     		rgb=(0.000000,0.375000,1.000000)
     		rgb=(0.000000,0.437500,1.000000)
     		rgb=(0.000000,0.500000,1.000000)
     		rgb=(0.000000,0.562500,1.000000)
     		rgb=(0.000000,0.625000,1.000000)
     		rgb=(0.000000,0.687500,1.000000)
     		rgb=(0.000000,0.750000,1.000000)
     		rgb=(0.000000,0.812500,1.000000)
     		rgb=(0.000000,0.875000,1.000000)
     		rgb=(0.000000,0.937500,1.000000)
     		rgb=(0.000000,1.000000,1.000000)
     		rgb=(0.062500,1.000000,0.937500)
     		rgb=(0.125000,1.000000,0.875000)
     		rgb=(0.187500,1.000000,0.812500)
     		rgb=(0.250000,1.000000,0.750000)
     		rgb=(0.312500,1.000000,0.687500)
     		rgb=(0.375000,1.000000,0.625000)
     		rgb=(0.437500,1.000000,0.562500)
     		rgb=(0.500000,1.000000,0.500000)
     		rgb=(0.562500,1.000000,0.437500)
     		rgb=(0.625000,1.000000,0.375000)
     		rgb=(0.687500,1.000000,0.312500)
     		rgb=(0.750000,1.000000,0.250000)
     		rgb=(0.812500,1.000000,0.187500)
     		rgb=(0.875000,1.000000,0.125000)
     		rgb=(0.937500,1.000000,0.062500)
     		rgb=(1.000000,1.000000,0.000000)
     		rgb=(1.000000,0.937500,0.000000)
     		rgb=(1.000000,0.875000,0.000000)
     		rgb=(1.000000,0.812500,0.000000)
     		rgb=(1.000000,0.750000,0.000000)
     		rgb=(1.000000,0.687500,0.000000)
     		rgb=(1.000000,0.625000,0.000000)
     		rgb=(1.000000,0.562500,0.000000)
     		rgb=(1.000000,0.500000,0.000000)
     		rgb=(1.000000,0.437500,0.000000)
     		rgb=(1.000000,0.375000,0.000000)
     		rgb=(1.000000,0.312500,0.000000)
     		rgb=(1.000000,0.250000,0.000000)
     		rgb=(1.000000,0.187500,0.000000)
     		rgb=(1.000000,0.125000,0.000000)
     		rgb=(1.000000,0.062500,0.000000)
     		rgb=(1.000000,0.000000,0.000000)
     		rgb=(0.937500,0.000000,0.000000)
     		rgb=(0.875000,0.000000,0.000000)
     		rgb=(0.812500,0.000000,0.000000)
     		rgb=(0.750000,0.000000,0.000000)
     		rgb=(0.687500,0.000000,0.000000)
     		rgb=(0.625000,0.000000,0.000000)
     		rgb=(0.562500,0.000000,0.000000)
     		rgb=(0.500000,0.000000,0.000000)
     	}
      }

      \begin{figure}[htbp]
        \centering
        \begin{tabular}{cc}
          \begin{tikzpicture}
            \begin{axis}[
              xlabel=$x$,ylabel=$t$,
              enlargelimits=false,
              axis on top, width=7.00cm,
              colorbar,
              colorbar style={point meta min=-16,point meta max=0,    ytick pos=right,
                tick label style={font=\footnotesize},
              },
              colorbar/width=3mm,    
              ]
              \addplot graphics [
              xmin=-3,xmax=3,
              ymin=0,ymax=10,
              ] {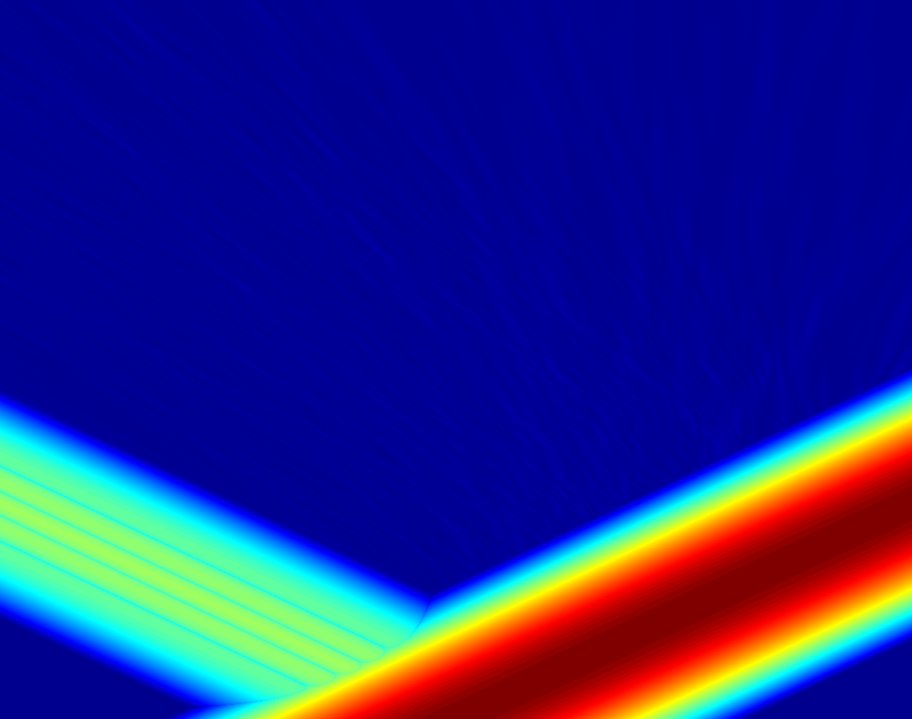};
            \end{axis}
          \end{tikzpicture}
          &
            \noindent \begin{tikzpicture}
              \begin{axis}[
                xlabel=$x$,ylabel=$t$,
                enlargelimits=false,
                axis on top, width=7.00cm,
                colorbar,
                colorbar style={point meta min=-16,point meta max=0,    ytick pos=right,
                  tick label style={font=\footnotesize},
                },
                colorbar/width=3mm,    
                ]
                \addplot graphics [
                xmin=-3,xmax=3,
                ymin=0,ymax=10,
                ] {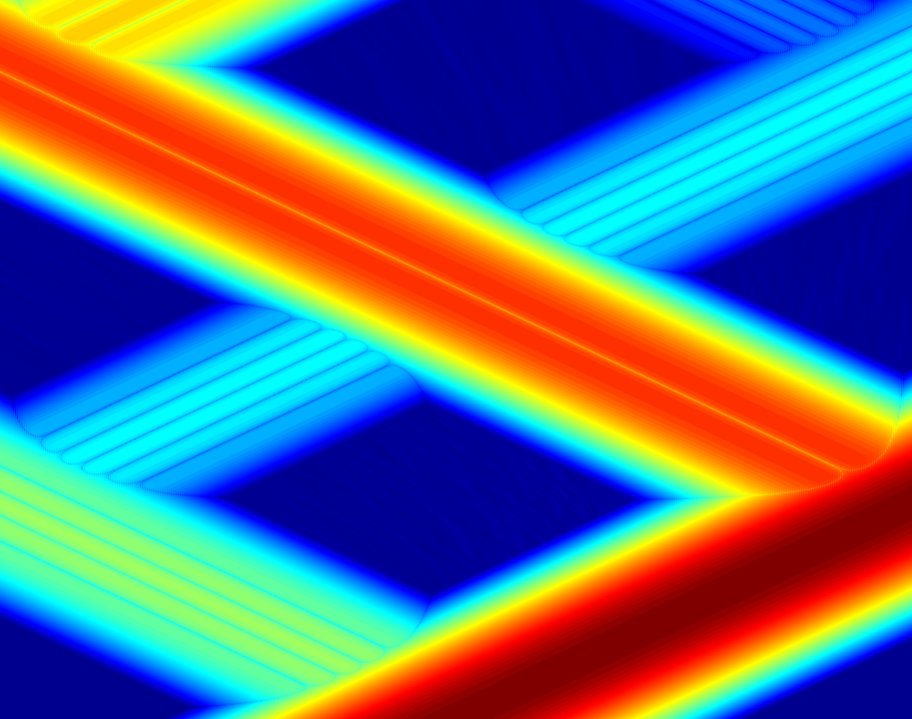};
              \end{axis}
            \end{tikzpicture}
          \\
        \end{tabular}
        \caption{Evolution of $\log_{10} |u_j^n|$ for DTBC \eqref{LF1d-DTBC} (left) and homogeneous Neumann BC $u_0^{n+2}=u_1^{n+1}$, 
        $u_{J+1}^{n+2}=u_J^{n+1}$ (right).}
        \label{fig:1d_evol_01}
      \end{figure}

      \begin{figure}[htbp]
        \centering
        \begin{tabular}{cc}
          \begin{tikzpicture}
            \begin{axis}[
              xlabel=$x$,ylabel=$t$,
              enlargelimits=false,
              axis on top, width=7.00cm,
              colorbar,
              colorbar style={point meta min=-16,point meta max=0,    ytick pos=right,
                tick label style={font=\footnotesize},
              },
              colorbar/width=3mm,    
              ]
              \addplot graphics [
              xmin=-3,xmax=3,
              ymin=0,ymax=10,
              ] {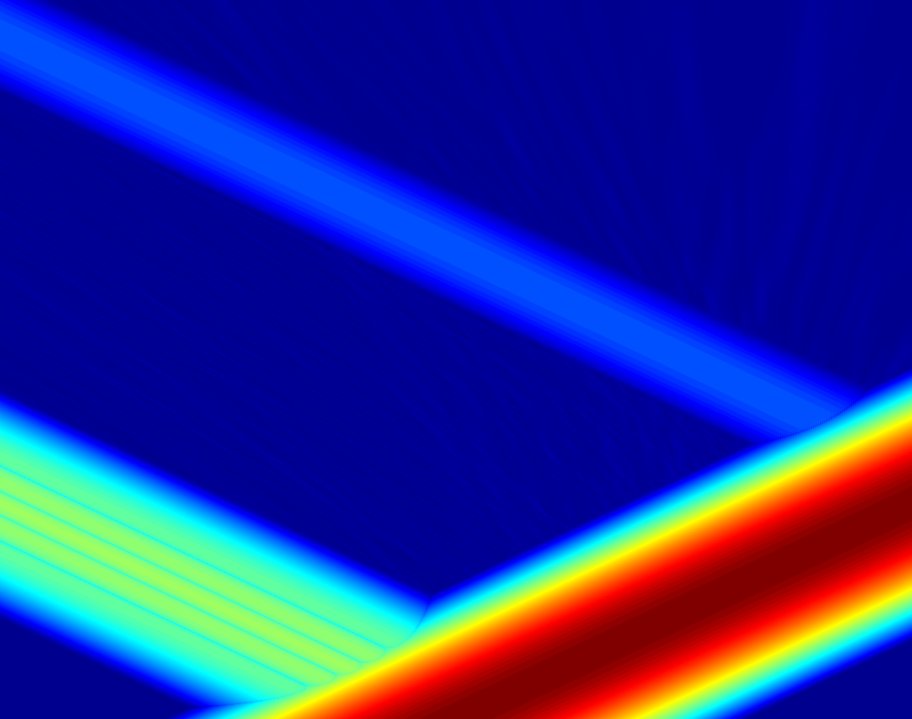};
            \end{axis}
          \end{tikzpicture}
          &
            \noindent \begin{tikzpicture}
              \begin{axis}[
                xlabel=$x$,ylabel=$t$,
                enlargelimits=false,
                axis on top, width=7.00cm,
                colorbar,
                colorbar style={point meta min=-16,point meta max=0,    ytick pos=right,
                  tick label style={font=\footnotesize},
                },
                colorbar/width=3mm,    
                ]
                \addplot graphics [
                xmin=-3,xmax=3,
                ymin=0,ymax=10,
                ] {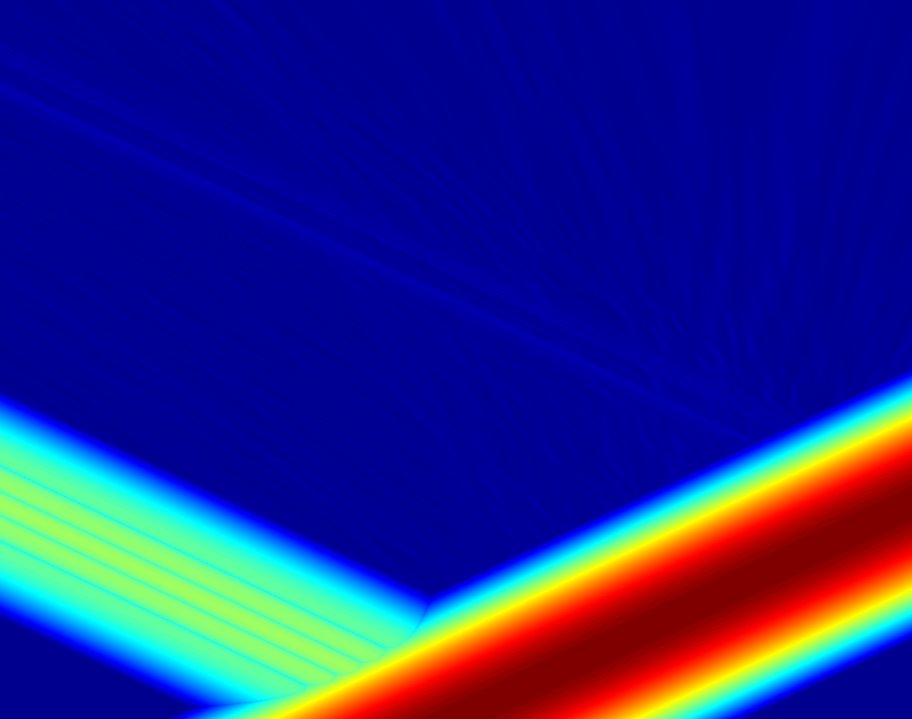};
              \end{axis}
            \end{tikzpicture}
          \\
        \end{tabular}
        \caption{Evolution of $\log_{10} |u_j^n|$ for the sum of exponentials BC \eqref{eq:dtbc_sumexp}: $(M,N)=(50,6)$ (left) and $(M,N)=(50,49)$ (right).}
        \label{fig:1d_evol_02}
      \end{figure}

\section{Exact and approximate DTBC for the two-dimensional leap-frog scheme}

\subsection{Exact DTBC on a half-space and local approximations}

In this paragraph, we go back to the leap-frog scheme \eqref{LF2d} for the linear advection equation \eqref{eq2} in two space dimensions. 
Our first goal is to derive the exact DTBC for \eqref{LF2d} when that scheme is considered on a half-space, be it $\{ j \ge 1 \}$, $\{ j \le J \}$, 
$\{ k \ge 1 \}$ or $\{ k \le K \}$. This means that we shall consider the leap-frog scheme \eqref{LF2d} in the whole plane $(j,k) \in \mathbb{Z}^2$ 
with initial data supported, say, in the half-space $\{ j \ge 1 \}$, and we shall derive the numerical boundary conditions satisfied by the solution 
to \eqref{LF2d} at $j=0$. This is exactly the problem considered in \cite{jfc} in the general framework of multistep finite difference schemes. 
For ease of reading and in order to stick to the notation of the previous Section, we consider from now on the case where the initial data for 
\eqref{LF2d} are supported in $\{ j \ge 1 \}$, and it is understood that in \eqref{LF2d} the tangential variable $k$ lies in $\mathbb{Z}$. We 
consider the step function $u_j^n (\cdot)$ defined by
$$
u_j^n(y) \, := \, u_{j,k}^n \, ,\quad \forall \, y \in [y_k,y_{k+1}) \, ,\quad k \in \mathbb{Z} \, ,
$$
with $y_k=y_b+k \, \delta y$ for all $k \in \mathbb{Z}$ and not only for $k=0,\dots,K+1$. Then \eqref{LF2d} reads
$$
u_j^{n+2}(y)-u_j^n(y) +\mu_x \, \big( u_{j+1}^{n+1}(y)-u_{j-1}^{n+1}(y) \big) +\mu_y \, \big( u_j^{n+1}(y+\delta y)-u_j^{n+1}(y-\delta y) \big) \, 
= \, 0 \, ,
$$
so applying first a partial Fourier transform with respect to $y$ and then the $\mathcal{Z}$-transform in time, we are led (with obvious notation) 
to the recurrence relation
\begin{equation}
\label{eq:LF2d_four_Z}
(z^2-1) \, \hat{u}_j(z,\eta) +\mu_x \, z \, \big( \hat{u}_{j+1}(z,\eta)-\hat{u}_{j-1}(z,\eta) \big) +2\, i \, \mu_y \, \sin(\eta \, \delta y) \, z \, \hat{u}_j(z,\eta) 
\, =0 \, ,
\end{equation}
which holds for all $j \le 0$ (in the half-space $\{ j \ge 1 \}$ where the initial data are supported, the analogous recurrence relation displays a 
nonzero source term on the right hand side).

The solution to the recurrence relation \eqref{eq:LF2d_four_Z} is computed by considering the characteristic equation:
\begin{equation}
\label{eq:rec2d}
(z^2-1) \, \kappa +\mu_x \, z \, (\kappa^2-1) +2\, i \,  \mu_y \, \sin \theta \, z \, \kappa \, = \, 0 \, ,
\end{equation}
where $\theta$ is a placeholder for the rescaled frequency $\eta \, \delta y$. For any real number $\theta$ and $|z|>1$, the roots to 
\eqref{eq:LF2d_four_Z} do not belong to $\mathbb{S}^1$. Therefore one has modulus $<1$ and the other one has modulus $>1$ (see 
the general splitting argument in \cite{jfc}). We label the two roots in such a way that:
\[
|\kappa_s(z,\theta)|<1 \, ,\qquad |\kappa_u(z,\theta)|>1 \, ,
\]
whenever $|z|>1$ and $\theta \in \mathbb{R}$. The discrete boundary condition (on the Fourier-$\mathcal{Z}$-transform side) thus reads
\begin{equation}
\label{eq:dtbc2d_1}
\hat{u}_0 (z,\eta) \, = \, \dfrac{1}{\kappa_u(z,\eta \, \delta y)} \, \hat{u}_1 (z,\eta) \, ,
\end{equation}
which is reminiscent of \eqref{eq:dtbc1d}. The analogous numerical boundary condition at $j=J+1$ (for the other half space problem) is
$$
\hat{u}_{J+1} (z,\eta) \, = \, \kappa_s(z,\eta \, \delta y) \, \hat{u}_J (z,\eta) \, .
$$

If we now apply the inverse Fourier-$\mathcal{Z}$-transform to \eqref{eq:dtbc2d_1}, this would lead to the exact, nonlocal DTBC for \eqref{LF2d} 
on a half space. Nonlocality refers here both to time and to the tangential variable $k$. In other words, determining $u^{n+2}_{0,k}$ requires the 
knowledge of $u_{1,k'}^\sigma$ for all $k' \in \mathbb{Z}$ and at all time steps $\sigma$ up to $n+1$. The nonlocality in time is not so harmful 
from a computational point of view (as we have seen in one space dimension) but the nonlocality with respect to $k$ is much more problematic 
since the set to which $k$ belongs is infinite. We therefore propose to follow the idea developed by Engquist and Majda \cite{engquist-majda} 
for continuous problems and to localize the exact DTBC \eqref{eq:dtbc2d_1} in the $k$-direction. This is performed by means of an asymptotic 
expansion with respect to $\theta=\eta \, \delta y$ in \eqref{eq:dtbc2d_1}. At a formal level, this amounts more or less to assuming that the discrete 
solution displays only bounded tangential frequencies $\eta$ and to performing an asymptotic expansion with respect to $\delta y$ (which is meant 
to be small in practice).

Let us therefore go back to \eqref{eq:rec2d} and observe that again the product of the two roots $\kappa_s$ and $\kappa_u$ equals $-1$. Hence 
\eqref{eq:dtbc2d_1} reads
$$
\hat{u}_0 (z,\eta) \, = \, -\kappa_s(z,\eta \, \delta y) \, \hat{u}_1 (z,\eta) \, ,
$$
and we now expand the stable root $\kappa_s(z,\theta)$ to \eqref{eq:rec2d} with respect to $\theta$. Up to the third order, this expansion reads
\begin{equation}
\label{eq:kappa_expansion}
\kappa_s(z,\theta) \, = \, \kappa_s^0(z) +2\, i \, \sin \theta \, \kappa_s^1(z) -4 \, \sin^2 \dfrac{\theta}{2} \, \kappa_s^2(z) +O(\theta^3) \, ,
\end{equation}
where $\kappa_s^0(z)$ is the stable root for the one-dimensional problem (which corresponds to $\theta=0$ in \eqref{eq:rec2d}) and the choice 
of writing $\sin \theta$ rather than $\theta$, as well as $\sin^2 \, \theta/2$ rather than $\theta^2$, has been made in order to exhibit amplification 
factors that are obviously linked with tangential finite difference operators (the centered first order derivative and the discrete Laplacian). At this 
stage, localizing the exact DTBC \eqref{eq:dtbc2d_1} with respect to $k$ amounts to using in \eqref{eq:dtbc2d_1} finitely many terms of the Taylor 
expansion \eqref{eq:kappa_expansion} (e.g., either the first term, or the two/three first terms). Namely, on the Fourier-$\mathcal{Z}$-transform, our 
localization procedure for \eqref{eq:dtbc2d_1} yields one of the three following choices (in increasing order of approximation):
\begin{subequations}
\label{eq:dtbc2d_2}
\begin{align}
\hat{u}_0 (z,\eta) \, &= \, -\kappa_s^0(z) \, \hat{u}_1 (z,\eta) \, ,\label{DTBC2d-a} \\
\hat{u}_0 (z,\eta) \, &= \, -\kappa_s^0(z) \, \hat{u}_1 (z,\eta) -2\, i \, \sin \, (\eta \, \delta y) \, \kappa_s^1(z) \, \hat{u}_1 (z,\eta) \, ,\label{DTBC2d-b} \\
\hat{u}_0 (z,\eta) \, &= \, -\kappa_s^0(z) \, \hat{u}_1 (z,\eta) -2\, i \, \sin \, (\eta \, \delta y) \, \kappa_s^1(z) \, \hat{u}_1 (z,\eta) 
+4 \, \sin^2 \dfrac{\eta \, \delta y}{2} \, \kappa_s^2(z) \, \hat{u}_1 (z,\eta) \, .\label{DTBC2d-c}
\end{align}
\end{subequations}
After performing an inverse Fourier transform with respect to $\eta$, \eqref{eq:dtbc2d_2} reads:
\begin{subequations}
\label{eq:dtbc2d_2'}
\begin{align}
\hat{u}_{0,k}(z) \, &= \, -\kappa_s^0(z) \, \hat{u}_{1,k}(z) \, ,\label{DTBC2d-k-a} \\
\hat{u}_{0,k}(z) \, &= \, -\kappa_s^0(z) \, \hat{u}_{1,k}(z) -\kappa_s^1(z) \, \big( \hat{u}_{1,k+1}(z)-\hat{u}_{1,k-1}(z) \big) \, ,\label{DTBC2d-k-b} \\
\hat{u}_{0,k}(z) \, &= \, -\kappa_s^0(z) \, \hat{u}_{1,k}(z) -\kappa_s^1(z) \, \big( \hat{u}_{1,k+1}(z)-\hat{u}_{1,k-1}(z) \big) \notag \\
& \qquad \qquad \qquad \qquad -\kappa_s^2(z) \, \big( \hat{u}_{1,k+1}(z)-2 \, \hat{u}_{1,k}(z) +\hat{u}_{1,k-1}(z) \big) \, ,\label{DTBC2d-k-c}
\end{align}
\end{subequations}
where the hat notation refers here to the $\mathcal{Z}$-transform only. In order to perform the inverse $\mathcal{Z}$-transform in \eqref{eq:dtbc2d_2'} 
and write down the approximate non-reflecting boundary conditions in the physical variables, we need to compute the Laurent series expansion of 
the functions $\kappa_s^1,\kappa_s^2$ that appear on the right hand side of \eqref{DTBC2d-k-b}, \eqref{DTBC2d-k-c}. (The Laurent series expansion 
of $\kappa_s^0$ has already been derived in the analysis of the one-dimensional problem.) This is achieved below with either of the two methods that 
we have already used in the one-dimensional case.
\bigskip

\paragraph{Expansions based on Legendre (and Tchebychev) polynomials.} Let us recall that we have already determined the Laurent series expansion 
of the function $\kappa_s^0$ in \eqref{eq:dtbc2d_2'}. We have written it under the form
$$
\kappa_s^0(z) \, = \, \sum_{n=1}^\infty \, \dfrac{\sigma_n^0}{z^n} \, = \, \sum_{n=0}^\infty \, \dfrac{s_n^0}{z^{2\, n+1}} \, ,
$$
where the sequence $(s_n^0)_{n \ge 0}$ is determined by either \eqref{eq:sn} or \eqref{recurrencesn0}. The two first values are $s_0^0=\mu_x$ and 
$s_1^0=\mu_x \, (1-\mu_x^2)$. Let us eventually recall that the expression of $\kappa_s^0(z)$ is given by \eqref{eq:stab_root}.

We are now going to determine the Laurent series expansion of $\kappa_s^1(z)$. The root $\kappa_s(z,\theta)$ to \eqref{eq:rec2d} is simple for 
$\theta \in \mathbb{R}$ and $|z|>1$. It thus depends holomorphically on $z$ and is $\mathcal{C}^\infty$ with respect to $\theta$. Differentiating 
\eqref{eq:rec2d} with respect to $\theta$ at $\theta=0$, and identifying $\kappa_s^1=\partial_\theta \kappa_s (z,0)/(2\, i)$, we obtain
\begin{equation}
\label{eq:eq_kappa_s_1}
\big( z^2-1 +2\, \mu_x \, z \, \kappa_s^0(z) \big) \, \kappa_s^1(z) \, = \, -\mu_y \, z \, \kappa_s^0(z) \, ,
\end{equation}
which yields
\begin{align}
\kappa_s^1(z) \, = \, -\dfrac{\mu_y \, \kappa_s^0}{\sqrt{(z-z^{-1})^2+4 \, \mu_x^2}} 
\, &= \, -\dfrac{\mu_y}{2\, \mu_x} \, \left\{ 1+\dfrac{z^{-2}-1}{\sqrt{(z^{-2}-1)^2 +4\, \mu_x^2 \, z^{-2}}} \right\} \notag \\
&= \, -\dfrac{\mu_y}{2\, \mu_x} \, \sum_{n=1}^\infty \dfrac{P_{n-1}(\alpha_x)-P_n(\alpha_x)}{z^{2\, n}} \, ,
\label{eq:kappa_s_1}
\end{align}
by using \eqref{eq:stab_root} and \eqref{eq:Legendre1}. We can therefore write the Laurent series expansion of $\kappa_s^1$ under the form:
\begin{equation}
\label{Laurentkappas1}
\kappa_s^1(z) \, = \, \sum_{n=0}^\infty \dfrac{s_n^1}{z^{2\, n}} \, , 
\end{equation}
with $s_0^1=0$ and
\begin{equation}
\label{eq:sn1}
\forall \, n \ge 1 \, ,\quad s_n^1 \, := \, \dfrac{\mu_y}{2\, \mu_x} \, \big( P_n(\alpha_x)-P_{n-1}(\alpha_x) \big) \, .
\end{equation}
We can, for instance, compute (recall $s_0^1=0$):
$$
s_1^1 \, = \, - \mu_x \, \mu_y \, ,\quad s_2^1 \, = \, - \mu_x \, \mu_y \, (2-3\, \mu_x^2) \, .
$$
The expression \eqref{eq:sn1} is only seemingly singular with respect to $\mu_x$ for $\mu_x\ll 1$. Indeed, we recall $\alpha_x=1-2\, \mu_x^2$ and 
since all Legendre polynomials satisfy $P_n(1)=1$, the expression $P_n(\alpha_x)-P_{n-1}(\alpha_x)$ can be written as $\mu_x^2 \, Q_n(\mu_x^2)$ 
for some polynomial $Q_n$. We therefore see from \eqref{eq:sn1} that all coefficients $s_n^1$ read $\mu_y \, \mu_x \, \mathcal{Q}_n(\mu_x^2)$ for 
some real polynomial $\mathcal{Q}_n$. Using again the Laplace formula \eqref{laplace}, we also get the asymptotic behavior:
\begin{equation}
\label{asymptotsn1}
s_n^1 \, = \, -\dfrac{\mu_y}{\mu_x} \, \sqrt{\dfrac{\tan (\theta_x/2)}{\pi \, n}} \, \sin \left( n \, \theta_x -\dfrac{\pi}{4} \right) +O(n^{-3/2}) \, ,
\end{equation}
where the angle $\theta_x \in (0,\pi)$ is still defined by $\cos \theta_x =1-2\, \mu_x^2$. The asymptotic behavior \eqref{asymptotsn1} can be verified by 
numerical experiments. Comparing with \eqref{asymptotsn0}, we observe that the sequence $(s_n^1)$ has a slower decay than $(s_n^0)$, but it decays 
to zero nevertheless. This is due to the fact that $\kappa_s^1$ has (four) singularities of the form $(z/z_0-1)^{-1/2}$, $z_0 \in \mathbb{S}^1$, as $z$ 
approaches the unit circle $\mathbb{S}^1$, while $\kappa_s^0$ is continuous up to $\mathbb{S}^1$.
\bigskip

\begin{figure}[htbp]
  \centering
  \begin{tabular}{ccc}
  \includegraphics[width=.32\textwidth]{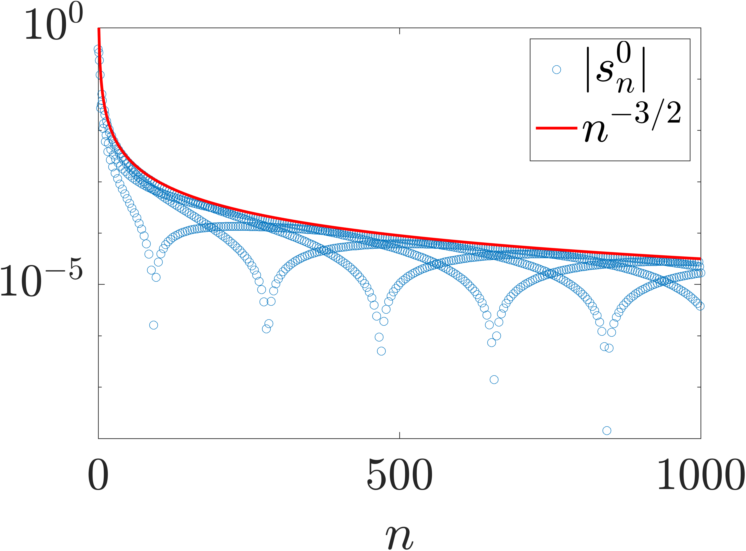}
    &
  \includegraphics[width=.32\textwidth]{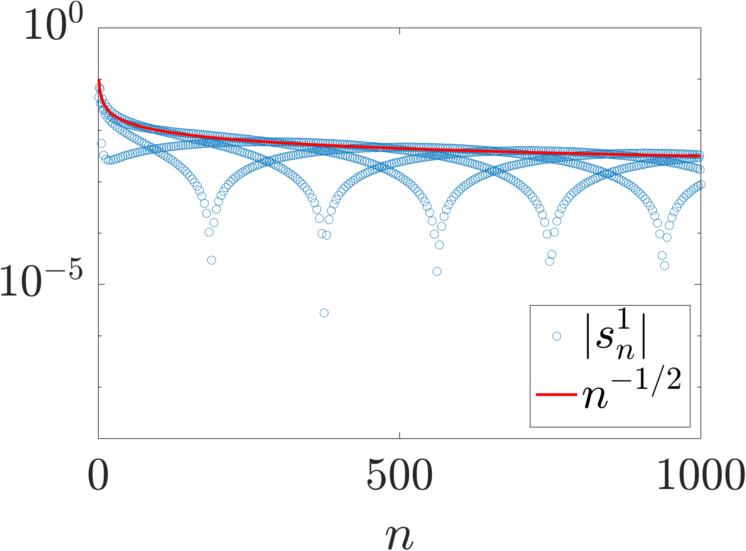}
    &
  \includegraphics[width=.32\textwidth]{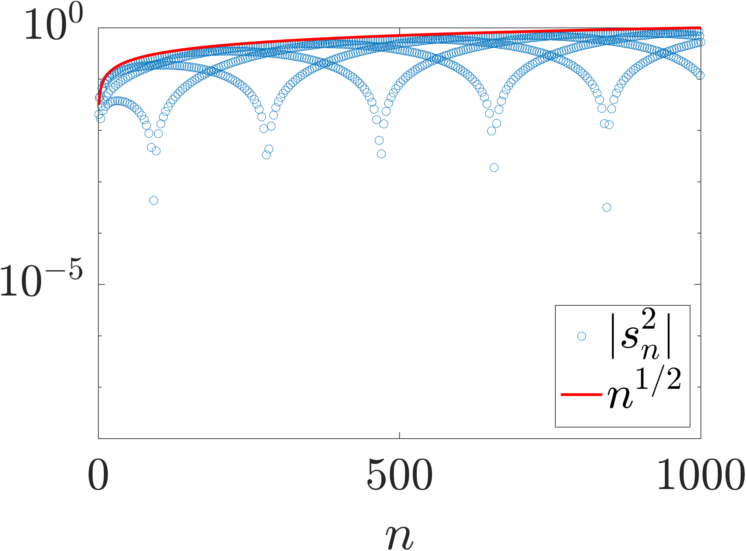}\\
  \end{tabular}
\caption{Asymptotic behaviour of $s_n^0$, $s_n^1$ and $s_n^2$}
\label{lapourifig}
\end{figure}

The Laurent series expansion of $\kappa_s^2$ can be derived by following more or less the same lines. By differentiating twice \eqref{eq:rec2d} with 
respect to $\theta$ and plugging the expression \eqref{eq:eq_kappa_s_1} of $\kappa_s^1$, we first get:
\begin{align}
\kappa_s^2(z)  \, &= \, -4\, z \, \kappa_s^1(z) \, \dfrac{\mu_y +\mu_x \, \kappa_s^1(z)}{z^2-1 +2\, \mu_x \, z \, \kappa_s^0(z)} \label{expressionkappas2} \\
&= \, \dfrac{4 \, \mu_y^2 \, z^2}{\big( z^2-1 +2\, \mu_x \, z \, \kappa_s^0(z) \big)^2} \, 
\dfrac{\kappa_s^0(z) \, \big( z^2-1 +\mu_x \, z \, \kappa_s^0(z) \big)}{z^2-1 +2\, \mu_x \, z \, \kappa_s^0(z)} \notag \\
&= \, \dfrac{4 \, \mu_y^2 \, z^2}{\big( z^2-1 +2\, \mu_x \, z \, \kappa_s^0(z) \big)^2} \, 
\dfrac{\mu_x \, z}{z^2-1 +2\, \mu_x \, z \, \kappa_s^0(z)} \, .\notag
\end{align}
In other words, we have derived the relation
$$
\kappa_s^2(z)  \, = \, \dfrac{4 \, \mu_x \, \mu_y^2}{z^3} \, \dfrac{1}{(1-z^{-2})^2+4\, \mu_x^2} \, \dfrac{1}{\sqrt{(1-z^{-2})^2+4\, \mu_x^2}} \, ,
$$
and it only remains to use \eqref{eq:Legendre1} as well as the generating function of the Tchebychev polynomials of the second kind (see \cite{szego}):
\begin{equation}
\label{eq:Tchebychev}
\dfrac{1}{1-2 \, x \, t+t^2} \, = \, \sum_{n=0}^\infty \, U_n(x) \, t^n \, ,
\end{equation}
to derive
$$
\kappa_s^2(z) \, = \, \dfrac{4 \, \mu_x \, \mu_y^2}{z^3} \, \sum_{n=0}^\infty 
\left( \sum_{m=0}^n U_m(\alpha_x) \, P_{n-m}(\alpha_x) \right) \, z^{-2\, n} \, .
$$
Writing the Laurent series expansion of $\kappa_s^2$ under the form:
\begin{equation}
\label{Laurentkappas2}
\kappa_s^2(z) \, = \, \sum_{n=0}^\infty \dfrac{s_n^2}{z^{2\, n+1}} \, ,
\end{equation}
we have derived the relations $s_0^2=0$ and
\begin{equation}
\label{eq:sn2}
\forall \, n \ge 1 \, ,\quad s_n^2 \, = \, 4 \, \mu_x \, \mu_y^2 \,  \sum_{m=0}^{n-1} U_m(\alpha_x) \, P_{n-1-m}(\alpha_x) \, ,
\end{equation}
with $\alpha_x=1-2\, \mu_x^2$. We can for instance compute
$$
s_1^2 \, = \, 4 \, \mu_x \, \mu_y^2 \, ,\quad s_2^2 \, = \, 12 \, \mu_x \, \mu_y^2 \, (1-2 \, \mu_x^2) \, .
$$
It does not seem so easy to infer from \eqref{laplace} and from the relation
$$
U_m(\cos \theta) \, = \, \dfrac{\sin ((m+1) \, \theta)}{\sin \theta} \, ,
$$
the asymptotic behavior of $(s_n^2)_{n \in \mathbb{N}}$. However, we verify below on some numerical experiments that $(s_n^2)_{n \in \mathbb{N}}$ 
does not seem to be bounded. 
It actually seems to grow like $\sqrt{n}$, see Figure \ref{lapourifig}, 
up to an oscillating behavior similar to the one we have exhibited in \eqref{asymptotsn0} and \eqref{asymptotsn1}. This growth in $\sqrt{n}$ is consistent 
with the singularities of the form $(z/z_0-1)^{-3/2}$ displayed by $\kappa_s^2$. The rigorous justification of this asymptotic behavior is left to a future work.

\paragraph{Inductive determination of the expansions.} As we have already seen in the analysis of the one-dimensional case, the Laurent series 
of $\kappa_s^0$ can be also determined inductively by plugging the expansion
$$
\kappa_s^0(z) \, = \, \sum_{n=1}^\infty \dfrac{\sigma_n^0}{z^n} \, ,
$$
by then first identifying $\sigma_1^0$ and then deriving an induction formula which gives the expression of $\sigma_n^0$ in terms of $\sigma_1^0, \dots, 
\sigma_{n-1}^0$. This led to the formula \eqref{inductionsigman0}. This strategy can be extended to the `correctors' $\kappa_s^1$ and $\kappa_s^2$. 
Namely, the function $\kappa_s^1$ is given by the relation \eqref{eq:eq_kappa_s_1}. Since $z \, \kappa_s^0(z)$ has a finite limit at infinity, we already 
see on \eqref{eq:eq_kappa_s_1} that $\kappa_s^1$ tends to zero at infinity. It actually decays like $z^{-2}$. Hence its Laurent series expansion reads
\[
\kappa_s^1(z) \, = \, \sum_{n=1}^\infty \dfrac{\sigma_n^1}{z^n} \, ,
\]
and the coefficients $(\sigma_n^1)_{n \ge 1}$ are computed by plugging the latter expression in \eqref{eq:eq_kappa_s_1} (which itself has been 
obtained by differentiating \eqref{eq:rec2d} with respect to $\theta$). After a few simplifications, we obtain $\sigma_1^1=0$ and the recurrence 
relation
\[
\forall \, n \ge 1 \, ,\quad  \sigma_{n+1}^1 -\sigma_{n-1}^1+2 \, \mu_x \, \sum_{m=0}^n \sigma_m^1 \, \sigma_{n-m}^0 \, = \, -\mu_y \, \sigma_n^0 \, ,
\]
where we use the convention $\sigma_0^1=0$. We easily deduce that all odd coefficients $\sigma_{2 \, p+1}^1$ vanish (recall that the even coefficients 
$\sigma_{2\, m}^0$ vanish), and the even coefficients $\sigma_{2 \, p}^1$, which we also write $s_p^1$ to be consistent with \eqref{Laurentkappas1}, 
satisfy $s_0^1=0$ and
\begin{equation}
\label{eq:coef_sigma_1}
\forall \, n \ge 0 \, ,\quad s_{n+1}^1 \, = \, s_n^1 -2 \, \mu_x \, \sum_{m=0}^n s_m^1 s_{n-m}^0 - \mu_y \, s_n^0 \, .
\end{equation}
Recalling the initial values $s_0^0=\mu_x$ and $s_1^0=\mu_x\, (1-\mu_x^2)$, we recover for instance from \eqref{eq:coef_sigma_1} the two first values 
$s_1^1=-\mu_x\, \mu_y$ and $s_2^1=-\mu_x\, \mu_y\, (2-3\, \mu_x^2)$ as in the previous method.

We can follow the same lines to derive the Laurent series expansion of $\kappa_s^2$, by starting from the relation \eqref{expressionkappas2} and 
plugging (observing on \eqref{expressionkappas2} that $\kappa_s^2$ tends to zero at infinity):
$$
\kappa_s^2(z) \, = \, \sum_{n=1}^\infty \dfrac{\sigma_n^2}{z^n} \, ,
$$
which yields $\sigma_1^2=0$ and the induction relation
\[
\forall \, n \ge 1 \, ,\quad  \sigma_{n+1}^2 -\sigma_{n-1}^2+2 \, \mu_x \, \sum_{m=0}^n \sigma_m^2 \, \sigma_{n-m}^0 \, = \, 
-4\, \mu_y \, \sigma_n^1 -4\, \mu_x \, \sum_{m=0}^n \sigma_m^1 \, \sigma_{n-m}^1 \, .
\]
Consistently with \eqref{Laurentkappas2}, we find that all even coefficients $\sigma_{2\, p}^2$ vanish, and the odd coefficients $\sigma_{2\, p+1}^2 
=s_p^2$ satisfy $s_0^2=0$ as well as the induction relation
$$
\forall \, n \ge 0 \, ,\quad  s_{n+1}^2 \, = \, s_n^2 -2 \, \mu_x \, \sum_{m=1}^n s_m^2 \, s_{n-m}^0 
-4\, \mu_y \, s_{n+1}^1 -4\, \mu_x \, \sum_{m=1}^n s_m^1 \, s_{n+1-m}^1 \, .
$$
For instance, we can recover the first value $s_1^2=s_0^2-4\, \mu_y \, s_1^2=4\, \mu_x \, \mu_y^2$ as given by the previous method.
\bigskip

Independently of the method we choose to compute the coefficients in the Laurent series expansion of $\kappa_s^1$ and $\kappa_s^2$, we can 
implement them numerically up to any prescribed final time index $N_f$. We now go back to \eqref{eq:dtbc2d_2'}. Once we have performed the 
inverse Fourier $\mathcal{Z}$-transform, the approximate DTBC \eqref{eq:dtbc2d_2'} respectively read:
\begin{subequations}
\label{eq:dtbc2d_left}
\begin{align}
u_{0,k}^{n+2} \, &= \, -\sum_{0 \le m \le (n+1)/2} \, s_m^0 \, u_{1,k}^{n+1-2\, m} \, ,\label{DTBC2d-left-a} \\
u_{0,k}^{n+2} \, &= \, -\sum_{0 \le m \le (n+1)/2} \, s_m^0 \, u_{1,k}^{n+1-2\, m} 
-\sum_{1 \le m \le (n+2)/2} \, s_m^1 \, (u_{1,k+1}^{n+2-2\, m}-u_{1,k-1}^{n+2-2\, m}) \, ,\label{DTBC2d-left-b} \\
u_{0,k}^{n+2} \, &= \, -\sum_{0 \le m \le (n+1)/2} \, s_m^0 \, u_{1,k}^{n+1-2\, m} 
-\sum_{1 \le m \le (n+2)/2} \, s_m^1 \, (u_{1,k+1}^{n+2-2\, m}-u_{1,k-1}^{n+2-2\, m}) \notag \\
& \qquad \qquad \qquad 
-\sum_{1 \le m \le (n+1)/2} \, s_m^2 \,  \big( u_{1,k+1}^{n+1-2\, m}-2\, u_{1,k}^{n+1-2\, m}+u_{1,k-1}^{n+1-2\, m} \big) \, ,\label{DTBC2d-left-c}
\end{align}
\end{subequations}
depending on the tangential accuracy we aim at achieving on the boundary.

Of course, the asymptotic expansion that we have derived for $\kappa_s(z,\theta)$ also provides the approximate DTBC at the right boundary 
$j=J+1$:
\begin{subequations}
\label{eq:dtbc2d_right}
\begin{align}
u_{J+1,k}^{n+2} \, &= \, \sum_{0 \le m \le (n+1)/2} \, s_m^0 \, u_{J,k}^{n+1-2\, m} \, ,\label{DTBC2d-right-a} \\
u_{J+1,k}^{n+2} \, &= \, \sum_{0 \le m \le (n+1)/2} \, s_m^0 \, u_{J,k}^{n+1-2\, m} 
+\sum_{1 \le m \le (n+2)/2} \, s_m^1 \, (u_{J,k+1}^{n+2-2\, m}-u_{J,k-1}^{n+2-2\, m}) \, ,\label{DTBC2d-right-b} \\
u_{J+1,k}^{n+2} \, &= \, \sum_{0 \le m \le (n+1)/2} \, s_m^0 \, u_{J,k}^{n+1-2\, m} 
+\sum_{1 \le m \le (n+2)/2} \, s_m^1 \, (u_{J,k+1}^{n+2-2\, m}-u_{J,k-1}^{n+2-2\, m}) \notag \\
& \qquad \qquad \qquad 
+\sum_{1 \le m \le (n+1)/2} \, s_m^2 \,  \big( u_{J,k+1}^{n+1-2\, m}-2\, u_{J,k}^{n+1-2\, m}+u_{J,k-1}^{n+1-2\, m} \big) \, .\label{DTBC2d-right-c}
\end{align}
\end{subequations}
Let us observe that the boundary conditions \eqref{eq:dtbc2d_left}, \eqref{eq:dtbc2d_right} are nonlocal in time but they are local in space 
(unlike the original exact DTBC). The boundary conditions \eqref{eq:dtbc2d_left}, \eqref{eq:dtbc2d_right} only involve a three point stencil, as 
depicted on Figure \ref{fig:coin} (left picture) on the example of the point located closest to the upper right corner. We could have pushed the 
Taylor expansion with respect to $\theta$ further by including for instance the following term $\theta^3$. The problem is that there is no linear 
combination of $1$, $\cos \theta$ and $\sin \theta$ that behaves like $\theta^3$ at the origin. In other words, if one writes $\theta^3$ as a 
trigonometric polynomial $p(\theta)$ up to an $O(\theta^4)$ term, then the trigonometric polynomial $p$ has degree at least $2$, which means 
that the associated finite difference operator in $k$ will involve at least $4$ points (and more likely five if one wishes to keep some symmetry). 
Hence this extension causes some trouble when getting close to the four corners of the rectangle because the very last boundary points close 
to the corner would require a specific treatment. We therefore do not pursue this higher order extension for the time being and leave it to future 
investigations.
\bigskip

At this stage, we have made the numerical boundary conditions on the `left' and `right' boundaries, resp. $\{ j=0 \}$ and $\{ j=J+1 \}$, explicit. 
The analysis can be carried out almost word for word for the `bottom' and `top' boundaries, namely $\{ k=0 \}$ and $\{ k=K+1 \}$. (Recall that 
our final goal is to perform numerical simulations on the rectangle depicted on Figure \ref{fig:maillage}.) The only thing is to observe that the 
indices $j,k$ in \eqref{LF2d} can be exchanged by simply switching the roles of $\mu_x$ and $\mu_y$. From now on, we let $(t_n^0)_{n \ge 0}$, 
$(t_n^1)_{n \ge 0}$, $(t_n^2)_{n \ge 0}$ denote the sequences obtained by the same procedure as $(s_n^0)_{n \ge 0}$, $(s_n^1)_{n \ge 0}$, 
$(s_n^2)_{n \ge 0}$, but switching the roles of $\mu_x$ and $\mu_y$. For instance, the sequence $(t_n^0)_{n \ge 0}$ is defined by $t_0^0=\mu_y$, 
$t_1^0 :=\mu_y \, (1-\mu_y^2)$, and
$$
\forall \, n \ge 2 \, ,\quad t_n^0 \, = \, \dfrac{2\, n-1}{n+1} \, (1-2\, \mu_y^2) \, t_{n-1}^0 -\dfrac{n-2}{n+1} \, t_{n-2}^0 \, .
$$
Then the approximate DTBC that we consider on the bottom and top boundaries read
\begin{subequations}
\label{eq:dtbc2d_bottom}
\begin{align}
u_{j,0}^{n+2} \, &= \, -\sum_{0 \le m \le (n+1)/2} \, t_m^0 \, u_{j,1}^{n+1-2\, m} \, ,\label{DTBC2d-bottom-a} \\
u_{j,0}^{n+2} \, &= \, -\sum_{0 \le m \le (n+1)/2} \, t_m^0 \, u_{j,1}^{n+1-2\, m} 
-\sum_{1 \le m \le (n+2)/2} \, t_m^1 \, (u_{j+1,1}^{n+2-2\, m}-u_{j-1,1}^{n+2-2\, m}) \, ,\label{DTBC2d-bottom-b} \\
u_{j,0}^{n+2} \, &= \, -\sum_{0 \le m \le (n+1)/2} \, t_m^0 \, u_{j,1}^{n+1-2\, m} 
-\sum_{1 \le m \le (n+2)/2} \, t_m^1 \, (u_{j+1,1}^{n+2-2\, m}-u_{j-1,1}^{n+2-2\, m}) \notag \\
& \qquad \qquad \qquad 
-\sum_{1 \le m \le (n+1)/2} \, t_m^2 \,  \big( u_{j+1,1}^{n+1-2\, m}-2\, u_{j,1}^{n+1-2\, m}+u_{j-1,1}^{n+1-2\, m} \big) \, ,\label{DTBC2d-bottom-c}
\end{align}
\end{subequations}
and
\begin{subequations}
\label{eq:dtbc2d_top}
\begin{align}
u_{j,K+1}^{n+2} \, &= \, \sum_{0 \le m \le (n+1)/2} \, t_m^0 \, u_{j,K}^{n+1-2\, m} \, ,\label{DTBC2d-top-a} \\
u_{j,K+1}^{n+2} \, &= \, \sum_{0 \le m \le (n+1)/2} \, t_m^0 \, u_{j,K}^{n+1-2\, m} 
+\sum_{1 \le m \le (n+2)/2} \, t_m^1 \, (u_{j+1,K}^{n+2-2\, m}-u_{j-1,K}^{n+2-2\, m}) \, ,\label{DTBC2d-top-b} \\
u_{j,K+1}^{n+2} \, &= \, \sum_{0 \le m \le (n+1)/2} \, t_m^0 \, u_{j,K}^{n+1-2\, m} 
+\sum_{1 \le m \le (n+2)/2} \, t_m^1 \, (u_{j+1,K}^{n+2-2\, m}-u_{j-1,K}^{n+2-2\, m}) \notag \\
& \qquad \qquad \qquad 
+\sum_{1 \le m \le (n+1)/2} \, t_m^2 \,  \big( u_{j+1,K}^{n+1-2\, m}-2\, u_{j,K}^{n+1-2\, m}+u_{j-1,K}^{n+1-2\, m} \big) \, .\label{DTBC2d-top-c}
\end{align}
\end{subequations}

Before reporting on several numerical simulations using the boundary conditions \eqref{eq:dtbc2d_left}, \eqref{eq:dtbc2d_right}, 
\eqref{eq:dtbc2d_bottom}, \eqref{eq:dtbc2d_top}, we first study whether each subcase in these boundary conditions, for instance 
\eqref{DTBC2d-left-a}, \eqref{DTBC2d-left-b} and \eqref{DTBC2d-left-c} satisfy some basic stability requirements. This stability 
analysis for \emph{half-space problems} is carried out in the following paragraph.

\subsection{Stability analysis}

The stability analysis of initial boundary value problems for hyperbolic systems in regions with corners is a delicate matter, see for instance 
\cite{osher2,osher3,sarason-smoller,benoit} and references therein for a partial account of the theory (of which much, if not all, remains to 
be done). The discrete counterpart is not easier and has been left (almost) undone so far. At a theoretical level, we do not claim to study the 
interaction of either one of the numerical boundary conditions in \eqref{eq:dtbc2d_top} with either of the conditions in \eqref{eq:dtbc2d_right} 
at the upper right corner of the rectangle. Such an analytical study seems to be out of reach at the present time. However we shall report 
below on some interesting numerical observations which call for a lot of care in the coupling of numerical strategies on each side of the 
rectangle. Namely, we shall show on one example that coupling two specific numerical boundary conditions on the top and right boundaries, 
namely \eqref{DTBC2d-top-c} and \eqref{DTBC2d-right-c} yields strong instabilities (even though the two separate half-space problems seem 
to be stable, as our partial proof below indicates). This phenomenon is well-known in the PDE context, see an elementary example in 
\cite{osher4} or more elaborate examples in \cite{sarason-smoller,benoit}.

We now examine a more simple problem which consists in studying the stability of the half-space problem where the leap-frog 
scheme \eqref{LF2d} holds in the half space $\{ (j,k) \in \mathbb{Z}^2 \, , \, j \ge 1 \}$ and either of the numerical boundary conditions 
\eqref{DTBC2d-left-a}, \eqref{DTBC2d-left-b}, \eqref{DTBC2d-left-c} is imposed on the boundary $\{ j=0 \}$. Following what we have 
done in the one-dimensional case, the stability analysis is carried out here by means of the so-called normal mode analysis. Extending 
slightly what we have done in the one-dimensional case, the normal mode analysis amounts here to determining the solutions to the 
leap-frog scheme \eqref{LF2d} which are of the form:
$$
u_{j,k}^n \, = \, z^n \, \exp (i \, k \, \theta) \, v_j \, ,
$$
with $|z|>1$, $\theta \in \mathbb{R}$ and $(v_j)_{j \ge 0} \in \ell^2$, that also satisfy the homogeneous numerical boundary conditions, 
be they \eqref{DTBC2d-left-a}, \eqref{DTBC2d-left-b} or \eqref{DTBC2d-left-c}. The arguments are quite similar in each of the three 
cases (except that we shall not be able to carry them out completely in the last case). We deal below with \eqref{DTBC2d-left-a}, 
\eqref{DTBC2d-left-b} and \eqref{DTBC2d-left-c} in increasing order of complexity.

\paragraph{The `zero order' numerical boundary condition \eqref{DTBC2d-left-a}.} First of all, we plug the expression $u_{j,k}^n =z^n \, 
\exp (i \, k \, \theta) \, v_j$ in \eqref{LF2d} and find that $(v_j)_{j \ge 0}$ must satisfy the recurrence relation
$$
\forall \, j \ge 1 \, ,\quad \big( z^2-1 \big) \, v_j \, + \, \mu_x \, z \, \big( v_{j+1}-v_{j-1} \big) +2\, i \, \sin \theta \, z \, v_j \, = 0 \, .
$$
(Compare with \eqref{eq:LF2d_four_Z}.) Since $(v_j)_{j \ge 0}$ should decay to zero at infinity, we get
$$
\forall \, j \ge 0 \, ,\quad v_j \, = \, v_0 \, \kappa_s(z,\theta)^j \, ,\quad v_0 \in \mathbb{C} \, ,
$$
where we recall that $\kappa_s(z,\theta)$ is the only root of modulus $<1$ to \eqref{eq:rec2d} for $|z|>1$ and $\theta \in \mathbb{R}$ (and 
$1/\kappa_s(z,\theta) =-\kappa_u(z,\theta)$ where $\kappa_u(z,\theta)$ is the only root of modulus $>1$ to the same equation).

Going back to \eqref{DTBC2d-left-a}, or rather to its Fourier-$\mathcal{Z}$-transform counterpart \eqref{DTBC2d-a}, we now need to determine 
whether there holds
$$
1 \, = \, -\kappa_s^0(z) \, \kappa_s(z,\theta) \, ,
$$
or equivalently
$$
\kappa_u(z,\theta) \, = \, \kappa_s^0(z) \, \, ( \, = \, \kappa_s(z,0) \, )
$$
for some pair $(z,\theta)$ with $|z|>1$ and $\theta \in \mathbb{R}$. The latter relation can never hold because $|\kappa_u(z,\theta)|>1$ and 
$|\kappa_s^0(z)|<1$ so the Godunov-Ryabenkii condition holds as for the one-dimensional problem which we have considered in the previous 
Section. Here again, there does not exist any unstable eigenvalue for the half-space problem $\{ j \ge 1 \}$ with the numerical boundary condition 
\eqref{DTBC2d-left-a}. The same argument would apply for any of the other half space problems in $\{ j \le J \}$, $\{ k \ge 1 \}$ or $\{ k \le K \}$ 
with their corresponding `zero order' numerical boundary condition.

\paragraph{The `first order' numerical boundary condition \eqref{DTBC2d-right-b}.} We can use part of the previous analysis, namely the determination 
of the `stable' solutions to \eqref{LF2d} of the prescribed normal mode form. However we now consider the numerical boundary condition 
\eqref{DTBC2d-right-b} or rather its Fourier-$\mathcal{Z}$-transform counterpart \eqref{DTBC2d-b}, so we need to determine whether there holds
$$
1 \, = \, -\kappa_s^0(z) \, \kappa_s(z,\theta) -2\, i \, \sin \theta \, \kappa_s^1(z) \,\kappa_s(z,\theta) \, ,
$$
or equivalently
\begin{equation}
\label{lop2dordre1}
\kappa_u(z,\theta) \, = \, \kappa_s^0(z) +2\, i \, \sin \theta \, \kappa_s^1(z) \, ,
\end{equation}
for some pair $(z,\theta)$ with $|z|>1$ and $\theta \in \mathbb{R}$. Let us recall that $\kappa_s^0$ is a root to \eqref{eq:rec1d} and that $\kappa_s^1$ 
is given by the relation \eqref{eq:eq_kappa_s_1}.

Let us first consider the case $\theta =0 \, (\pi)$, that is, $\sin \theta =0$. Then $\kappa_u(z,\theta) =\kappa_u^0(z)$ and \eqref{lop2dordre1} can not hold 
since the left hand side has modulus $>1$ while the right hand side has modulus $<1$. We therefore assume from now on $\sin \theta \neq 0$. Let us 
assume for a moment that the equality \eqref{lop2dordre1} holds for some pair $(z,\theta)$ with $|z|>1$ and $\sin \theta \neq 0$. Then the right hand side 
of \eqref{lop2dordre1} must be a solution to the second degree polynomial equation \eqref{eq:rec2d}. Plugging the right hand side of \eqref{lop2dordre1} 
in \eqref{eq:rec2d} and expanding, we eventually get\footnote{Here we use the equation \eqref{eq:rec1d} satisfied by $\kappa_s^0$ as well as the relation 
\eqref{eq:eq_kappa_s_1} satisfied by $\kappa_s^1$.}:
$$
(2\, i \, \sin \theta)^2 \, z \, \kappa_s^1(z) \, \big( \mu_y +\mu_x \, \kappa_s^1(z) \big) \, = \, 0 \, .
$$
Simplifying by $\sin \theta$, we get either $\kappa_s^1(z)=0$ or $\mu_y +\mu_x \, \kappa_s^1(z)=0$. Let us exclude the first of these two options. 
For $|z|>1$, the stable root $\kappa_s^0(z)$ to \eqref{eq:rec1d} does not vanish. Hence we see on the relation \eqref{eq:eq_kappa_s_1} that 
$\kappa_s^1(z)$ does not vanish either. This means that if \eqref{lop2dordre1} holds for some pair $(z,\theta)$, then there necessarily holds 
$\mu_y +\mu_x \, \kappa_s^1(z)=0$. We use again \eqref{eq:eq_kappa_s_1} to compute:
\begin{align*}
\mu_y +\mu_x \, \kappa_s^1(z) \, &= \, \mu_y \, \dfrac{z^2-1+\mu_x \, z \, \kappa_s^0(z)}{z^2-1+2\, \mu_x \, z \, \kappa_s^0(z)} \\
&= \, -\dfrac{\mu_x \, \mu_y \, \kappa_u^0(z)}{z^2-1+2\, \mu_x \, z \, \kappa_s^0(z)} \, ,
\end{align*}
and the latter quantity can vanish only if $\mu_y=0$ (recall that we have assumed $\mu_x \neq 0$). At this point of the analysis, we have shown that the 
only possibility for \eqref{lop2dordre1} to hold is to have $\mu_y=0$, that is $c_y=0$. However, in that case, we have $\kappa_u(z,\theta) =\kappa_u^0(z)$ 
for all $\theta$, see \eqref{eq:rec2d}, and $\kappa_s^1(z)=0$ (see \eqref{eq:eq_kappa_s_1}). Hence \eqref{lop2dordre1} can never hold because for 
$\mu_y=0$, the right hand side has modulus $<1$ while the left hand side has modulus $>1$. In other words, the Godunov-Ryabenkiii condition holds 
and there does not exist any unstable eigenvalue for the half-space problem $\{ j \ge 1 \}$ with the numerical boundary condition \eqref{DTBC2d-left-b}. 
Of course, the same arguments would apply for any of the other half space problems.

\paragraph{The `second order' numerical boundary condition \eqref{DTBC2d-right-c}.} We can use the same approach as above. Verifying the 
Godunov-Ryabenkii condition amounts to determining whether there exists a pair $(z,\theta)$ with $|z|>1$ and $\theta \in \mathbb{R}$ such that:
\begin{equation}
\label{lop2dordre2}
\kappa_u(z,\theta) \, = \, \kappa_s^0(z) +2\, i \, \sin \theta \, \kappa_s^1(z) -4\, \sin^2 \dfrac{\theta}{2} \, \kappa_s^2(z) \, .
\end{equation}
We verify again that \eqref{lop2dordre2} can not hold if $\theta =0 \, (\pi)$ or if $\mu_y =0$. Hence we assume from now on $\sin \theta \neq 0$ and 
$\mu_y \neq 0$. Let us assume that the relation \eqref{lop2dordre2} holds, meaning that not only the right hand side of \eqref{lop2dordre2} is a root 
to \eqref{eq:rec2d} but also that it is the only root of modulus $>1$. We first plug the right hand side of \eqref{lop2dordre2} in \eqref{eq:rec2d} and 
collect the terms in increasing powers of $\sin \theta/2$ (we use again the equation \eqref{eq:rec1d} satisfied by $\kappa_s^0$, the relation 
\eqref{eq:eq_kappa_s_1} satisfied by $\kappa_s^1$ and we also use the relation \eqref{expressionkappas2} satisfied by $\kappa_s^2$). After 
simplifying (quite a bit...), we eventually get the polynomial equation
\begin{multline}
\label{lop2dordre2-eq}
i \, \sin \dfrac{\theta}{2} \, \big(  (z^2-1)^4 +8\, \mu_x^2 \, z^2 \, (z^2-1)^2 +16 \, \mu_x^2 \, (\mu_x^2-\mu_y^2) \, z^4 \big) \\
-4 \, \cos \dfrac{\theta}{2} \, \mu_y \, z \, (z^2-1) \, \big( (z^2-1)^2 +4 \, \mu_x^2 \, z^2 \big) \, = \, 0 \,  .
\end{multline}
Since $\sin \theta$ is nonzero, then $\sin \theta/2$ is also nonzero, and \eqref{lop2dordre2-eq} is consequently an eight degree polynomial 
equation in $z$. Changing $z$ into $i \, z$, \eqref{lop2dordre2-eq} becomes a polynomial equation with real coefficients (after simplifying by 
$i$). Moreover, the roots of \eqref{lop2dordre2-eq} are nonzero and are invariant by the transformation $z \to -1/z$. Let us also recall that we 
are interested in determining the roots of \eqref{lop2dordre2-eq} that satisfy $|z|>1$.

We have not been able to obtain a complete proof of the facts we claim below, hence we do not have a complete proof for the verification of the 
Godunov-Ryabenkii condition, but repeated numerical experiments indicate the following facts:
\begin{itemize}
 \item For $\theta \in (0,\theta_{\rm max})$, where the angle $\theta_{\rm max}$ is defined by
$$
\tan \dfrac{\theta_{\rm max}}{2} \, = \, \dfrac{2\, \mu_y \, (1-\mu_x^2)}{(1-\mu_x^2)^2 -\mu_x^2 \, \mu_y^2} \, ,
$$
 \eqref{lop2dordre2-eq} has a unique root $z_\theta$ verifying $|z_\theta|>1$. This root is a purely imaginary number (one other root of 
 \eqref{lop2dordre2-eq} is $-1/z_\theta$ since the roots of \eqref{lop2dordre2-eq} are invariant under the transformation $z \to -1/z$, and 
 the six other roots have modulus one).
 
 \item For $\theta \in [\theta_{\rm max},\pi)$, all roots to \eqref{lop2dordre2-eq} have modulus one (the value $\theta_{\rm max}$ is obtained by 
 considering the case where $-i$ is a (double) root to \eqref{lop2dordre2-eq}).
\end{itemize}

For \eqref{lop2dordre2-eq} to have a root $z$ that satisfies $|z|>1$, we must necessarily have $\theta \in (0,\theta_{\rm max})$, and in that case the 
remaining question is to determine whether the quantity
$$
\kappa_s^0(z_\theta) +2\, i \, \sin \theta \, \kappa_s^1(z_\theta) -4\, \sin^2 \dfrac{\theta}{2} \, \kappa_s^2(z_\theta) 
$$
is either the stable or the unstable root to \eqref{eq:rec2d}. On the numerical experiments we have conducted, we have verified the property:
\begin{equation}
\label{inegmodule}
\big| \kappa_s^0(z_\theta) +2\, i \, \sin \theta \, \kappa_s^1(z_\theta) -4\, \sin^2 \dfrac{\theta}{2} \, \kappa_s^2(z_\theta) \big| \, < \, 1 \, ,
\end{equation}
for $\theta \in (0,\theta_{\rm max})$, meaning that the right hand side of \eqref{lop2dordre2} may coincide with a root \eqref{eq:rec2d} for 
some specific values of $(z,\theta)$ but in that case it coincides with the stable root $\kappa_s(z,\theta)$ and not with the unstable root 
$\kappa_u(z,\theta)$ (which would complete the verification of the Godunov-Ryabenkii condition). The inequality \eqref{inegmodule} can 
be shown rigorously in the regime where $\theta$ is small for in that case, we have $z_\theta \to \infty$, see \eqref{lop2dordre2}, and all 
functions $\kappa_s^0$, $\kappa_s^1$, $\kappa_s^2$ tend to zero at infinity, so the modulus of
$$
\kappa_s^0(z_\theta) +2\, i \, \sin \theta \, \kappa_s^1(z_\theta) -4\, \sin^2 \dfrac{\theta}{2} \, \kappa_s^2(z_\theta) 
$$
tends to zero as $\theta$ tends to zero. We can also verify \eqref{inegmodule} in the limit regime $\theta \to \theta_{\rm max}$ for in that case 
we get $z_\theta \to -i$ and we compute
$$
\kappa_s^0(z_\theta) \to \dfrac{i \, \big( 1-\sqrt{1-\mu_x^2} \big)}{\mu_x} \, ,\quad 
\kappa_s^1(z_\theta) \to \dfrac{\mu_y \, \big( 1-\sqrt{1-\mu_x^2} \big)}{2 \, \mu_x \, \sqrt{1-\mu_x^2}} \, ,\quad 
\kappa_s^2(z_\theta) \to -\dfrac{i \, \mu_y^2 \, \mu_x}{2 \, (1-\mu_x^2)^{3/2}} \, .
$$
We can then verify numerically that the inequality
$$
\left| \dfrac{i \, \big( 1-\sqrt{1-\mu_x^2} \big)}{\mu_x} 
+2 \, i \, \sin \theta_{\rm max} \, \dfrac{\mu_y \, \big( 1-\sqrt{1-\mu_x^2} \big)}{2 \, \mu_x \, \sqrt{1-\mu_x^2}} 
+4 \, \sin^2 \dfrac{\theta_{\rm max}}{2} \, \dfrac{i \, \mu_y^2 \, \mu_x}{2 \, (1-\mu_x^2)^{3/2}} \right| \, < \, 1
$$
holds for all positive values of $\mu_x,\mu_y$ such that $\mu_x+\mu_y<1$ (here we use the classical formula $\sin \theta_{\rm max}=2\, t/(1+t^2)$ 
and $\sin^2 \theta_{\rm max}/2=t^2/(1+t^2)$ with $t:=\tan \theta_{\rm max}/2$). Hence it does seem that \eqref{lop2dordre2} can not hold for some 
pair $(z,\theta)$, which means that the Godunov-Ryabenkii condition is verified (though a complete analytical proof of this fact still remains open). 
We postpone this rigorous justification to a future work.

\subsection{Numerical experiments on a rectangle}

We go back to our original motivation which is the simulation of the leap-frog scheme \eqref{LF2d} on the rectangle $\{ 1 \le j \le J \, , \, 1 \le k \le K \}$ 
with some approximate DTBC on each side. In the numerical simulations reported below, the computations were run on the rectangle $(x,y) \in (-3,3) 
\times(-2,2)$ and on the time interval $t\in[0,8]$. The discretization parameters are $J=300$ and $K=200$ (from which one can deduce the values of 
the space steps $\delta x$, $\delta y$). The time step $\delta t$ is chosen in order to satisfy the CFL condition
\[
\mu_x +\mu_y \, = \, \dfrac{1}{2} \, .
\]
(Recall the $\ell^2$-stability condition \eqref{cfl2D} for \eqref{LF2d} on the whole space.) The initial datum that we consider for the transport equation 
\eqref{eq2} is a Gaussian, namely $u_0(x,y) = \exp (-5\, (x^2+y^2))$. The discrete initial condition $(u_{j,k}^0)$ is defined by setting
$$
\forall \, j=0,\dots,J+1 \, ,\quad \forall \, k=0,\dots,K+1 \, ,\quad u_{j,k}^0 \, := \, u_0(x_j,y_k) \, .
$$
The first time step value $(u_{j,k}^1)$ is defined by imposing the second order two-dimensional Lax-Wendroff scheme \cite{laxwendroff}, namely:
\begin{align*}
\forall \, j=1,\dots,J \, ,\quad &\forall \, k=1,\dots,K \, ,\\
u_{j,k}^1 \, := \, u_{j,k}^0 &-\dfrac{\mu_x}{2} \, (u_{j+1,k}^0-u_{j-1,k}^0) -\dfrac{\mu_y}{2} \, (u_{j,k+1}^0-u_{j,k-1}^0) \\
&+\dfrac{\mu_x^2}{2} \, (u_{j+1,k}^0-2\, u_{j,k}^0+u_{j-1,k}^0) +\dfrac{\mu_y^2}{2} \, (u_{j,k+1}^0-2\, u_{j,k}^0+u_{j,k-1}^0) \\
&+\dfrac{\mu_x\, \mu_y}{4} \, (u_{j+1,k+1}^0-u_{j+1,k-1}^0-u_{j-1,k+1}^0+u_{j-1,k-1}^0) \, .
\end{align*}
The boundary values for $u^1$ (all other relevant values of $j,k$) are set equal to zero for simplicity.

In the numerical simulations below, we try various coupling strategies of numerical boundary conditions on the `right' and `top' boundaries $\{j=J+1\}$ 
and $\{k=K+1\}$. We first make it clear on Figure \ref{fig:coin} that the two sides of the upper right corner are coupled through the boundary conditions 
if one chooses for instance to impose \eqref{DTBC2d-right-b} on $j=J+1$ and \eqref{DTBC2d-top-b} on $k=K+1$. The coupling is weak though since 
the boundary conditions \eqref{DTBC2d-right-b} and \eqref{DTBC2d-top-b} are explicit so the trace values $u_{J+1,K}^{n+2}$ and $u_{J,K+1}^{n+2}$ 
are coupled only through the preceding time steps.

\begin{figure}[htbp]
\begin{center}
\begin{tikzpicture}[scale=2,>=latex]
\draw [thin, dashed] (-3.5,0.5) grid [step=0.5] (-2,2);
\draw [thin, dashed] (2,0.5) grid [step=0.5] (3.5,2);
\node (centre) at (-2,2){$\bullet$};
\node (centre) at (3.5,2){$\bullet$};
\node (centre) at (-2,1.5){${\color{red}\otimes}$};
\node (centre) at (-2.5,1){${\color{blue}\otimes}$};
\node (centre) at (-2.5,1.5){${\color{blue}\otimes}$};
\node (centre) at (-2.5,2){${\color{blue}\otimes}$};
\node (centre) at (3,2){${\color{red}\otimes}$};
\node (centre) at (2.5,1.5){${\color{blue}\otimes}$};
\node (centre) at (3,1.5){${\color{blue}\otimes}$};
\node (centre) at (3.5,1.5){${\color{blue}\otimes}$};
\draw (3.35,0.35) node[right]{$x_{J+1}$};
\draw (-2.15,0.35) node[right]{$x_{J+1}$};
\draw (2.85,0.35) node[right]{$x_J$};
\draw (-2.65,0.35) node[right]{$x_J$};
\draw (-3.8,2.1) node[right]{$y_{K+1}$};
\draw (1.7,2.1) node[right]{$y_{K+1}$};
\draw (-3.8,1.6) node[right]{$y_K$};
\draw (1.7,1.6) node[right]{$y_K$};
\draw [thin, dashed] (2,0.5) -- (2,2);
\draw (-3.5,2) -- (-2,2);
\draw (2,2) -- (3.5,2);
\draw (-2,0.5) -- (-2,2);
\draw (3.5,0.5) -- (3.5,2);
\end{tikzpicture}
\caption{The stencil of the numerical boundary conditions \eqref{DTBC2d-right-b} and \eqref{DTBC2d-top-b} near the upper right corner. Left: the 
points (in blue crosses) involved in the computation of $u_{J+1,K}^{n+2}$ (in red). Right: the points (in blue crosses) involved in the computation 
of $u_{J,K+1}^{n+2}$ (in red).}
\label{fig:coin}
\end{center}
\end{figure}
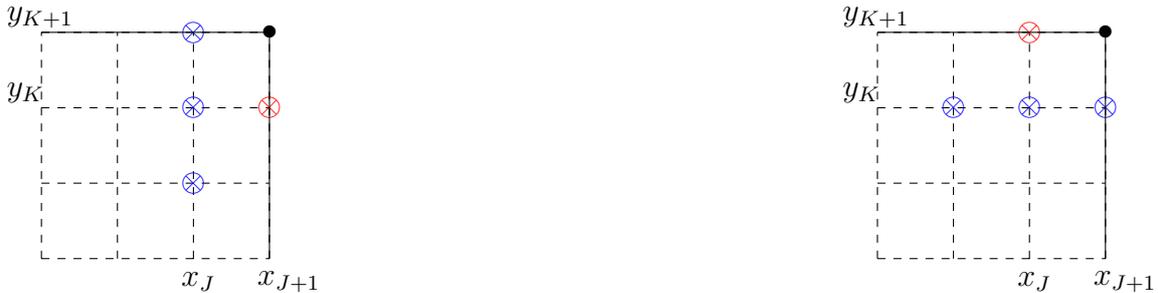

\begin{figure}[htbp]
\centering
\begin{tikzpicture}
\pgfplotscolorbardrawstandalone[ 
   colormap/jet,
   colorbar horizontal,
   point meta min=-10,
   point meta max=0,
   colorbar style={
      width=5cm,
      xtick={-10,-8,...,0}}]
\end{tikzpicture}
\caption{Legend of the 2D representation of the evolution of $\log_{10} |u_{j,k}^n|$.}
\label{fig:legend2D}
\end{figure}
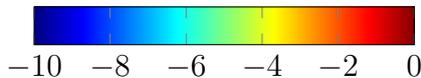

\begin{figure}[htbp]
\centering
\def\monechelle{1}
\begin{tikzpicture}[scale=\monechelle,axes/.style={scale=\monechelle}]
  \centering
  \fill[fill={rgb,255:red,128; green,0; blue,0}]     (8.0,0.0) rectangle (8.5,0.3); 
  \fill[fill={rgb,255:red,255; green,83; blue,25}]   (6.4,0.0) rectangle (6.9,0.3);
  \fill[fill={rgb,255:red,242; green,255; blue,64}]  (4.8,0.0) rectangle (5.3,0.3);
  \fill[fill={rgb,255:red,131; green,255; blue,226}] (3.2,0.0) rectangle (3.7,0.3);
  \fill[fill={rgb,255:red,153; green,191; blue,255}] (1.6,0.0) rectangle (2.1,0.3);
  \fill[fill={rgb,255:red,204; green,204; blue,233}] (0.0,0.0) rectangle (0.5,0.3);
  \draw foreach[count=\n from -7] \x in {0.5,2.1,...,10}{(\x,0.18) node[right,axes]{$10^{\n}$}};  
  \draw (-.1,-.1) rectangle (9.5,0.43);
\end{tikzpicture}
\caption{Legend of the 3D representation of the evolution of $\log_{10} |u_{j,k}^n|$.}
\label{fig:legend}
\end{figure}
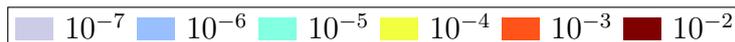

In order to highlight the wave reflections on each side of the rectangle, we show below on some two and three dimensional figures the logarithm of the 
(absolute value of the) solution. The logarithm of the solution is represented on two-dimensional figures at various time steps. The full time evolution is 
represented on a three dimensional figure. Logarithmic representations can be far more enlightening if one wishes to reveal the magnitude of the reflected 
waves. The color legend for all 2D logarithmic representations reported below is given in Figure \ref{fig:legend2D} and the color legend for all 3D logarithmic 
representations reported below is given in Figure \ref{fig:legend}.

\def\original{3cm}
\def\originale{60.73mm}
\def\oldtabcolsep{\tabcolsep}
\def\oldtarraystretch{\arraystretch}
\setlength{\tabcolsep}{0pt}
\renewcommand{\arraystretch}{0}
\def\newsize{3.0cm}
\def\mascale{.8}

\begin{figure}[htbp]
\begin{minipage}[htbp]{.5\linewidth}
  \begin{tabular}[htbp]{cc}
    \logplottext{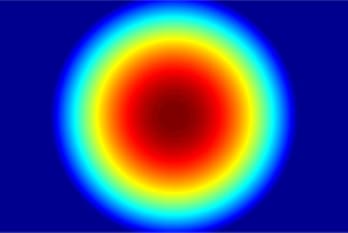}{\newsize}{\original}{\mascale}{t=0}   & 
    \logplottext{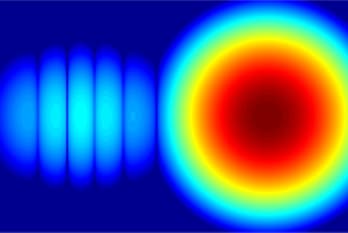}{\newsize}{\original}{\mascale}{t=1.6} \\
    \logplottext{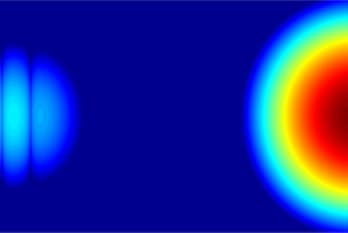}{\newsize}{\original}{\mascale}{t=3.2} & 
    \logplottext{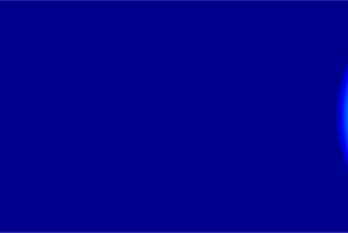}{\newsize}{\original}{\mascale}{t=4.8} \\
    \logplottext{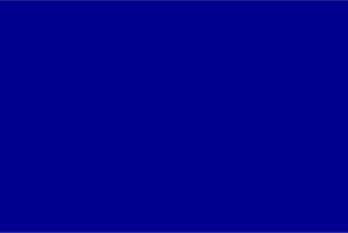}{\newsize}{\original}{\mascale}{t=6.4} & 
    \logplottext{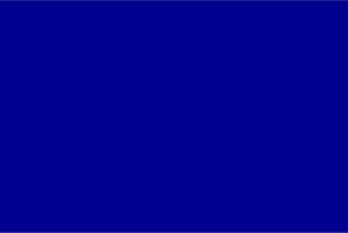}{\newsize}{\original}{\mascale}{t=8}   \\
  \end{tabular}
\end{minipage}
\begin{minipage}[htbp]{.5\linewidth}
  \def\newsizee{.8\textwidth}
  \myTDrepr{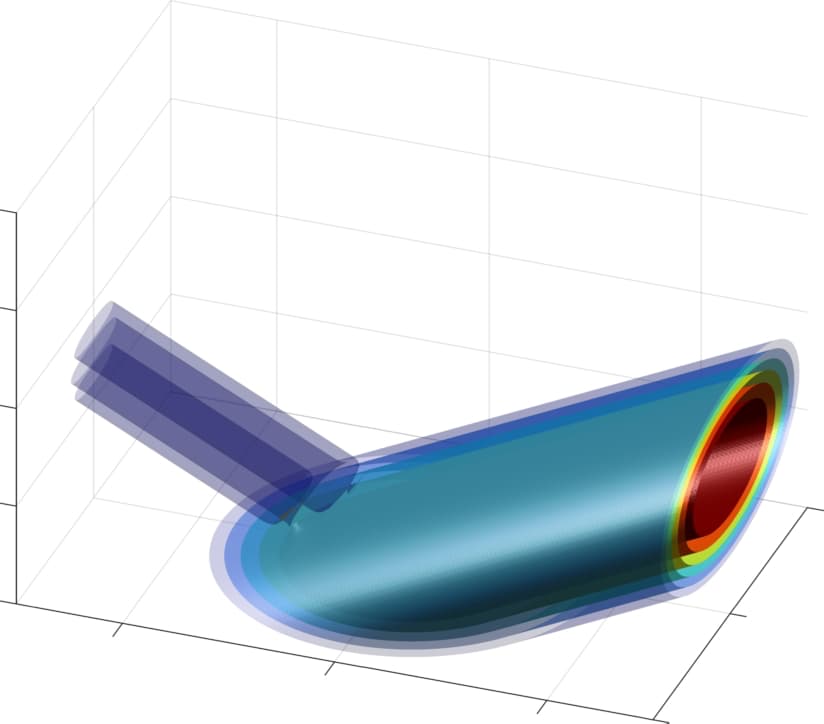}{\newsizee}{\originale}
\end{minipage}
\caption{Evolution of $\log_{10} |u_{j,k}^n|$ (left: six successive times, right: full time evolution) for $\mathbf{c}=(1,0)$ and DTBC of order $0$.}
\label{fig:test2D-1}
\end{figure}

\begin{figure}[htbp]
\begin{minipage}[htbp]{.5\linewidth}
  \begin{tabular}[htbp]{cc}
    \logplottext{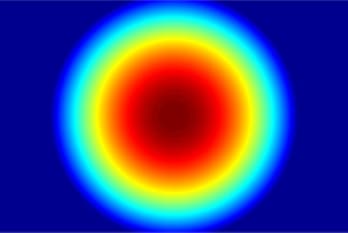}{\newsize}{\original}{\mascale}{t=0}   &
    \logplottext{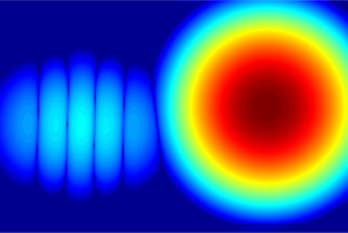}{\newsize}{\original}{\mascale}{t=1.6} \\    
    \logplottext{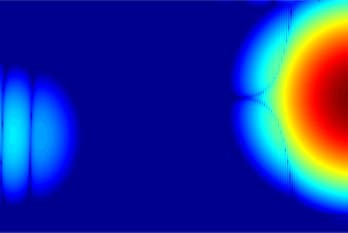}{\newsize}{\original}{\mascale}{t=3.2} &
    \logplottext{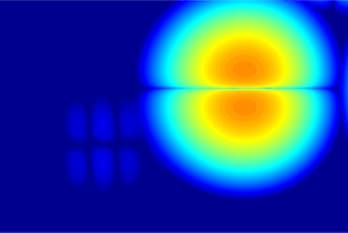}{\newsize}{\original}{\mascale}{t=4.8} \\    
    \logplottext{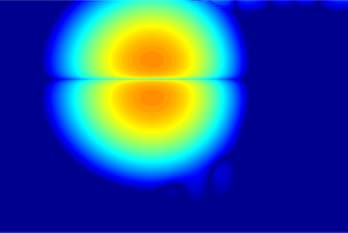}{\newsize}{\original}{\mascale}{t=6.4} &
    \logplottext{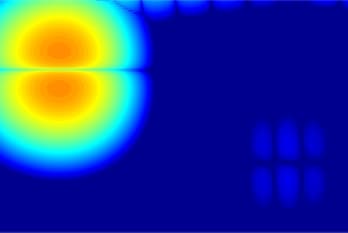}{\newsize}{\original}{\mascale}{t=8}   \\    
  \end{tabular}
\end{minipage}
\begin{minipage}[htbp]{.5\linewidth}
  \def\newsizee{.8\textwidth}
  \myTDrepr{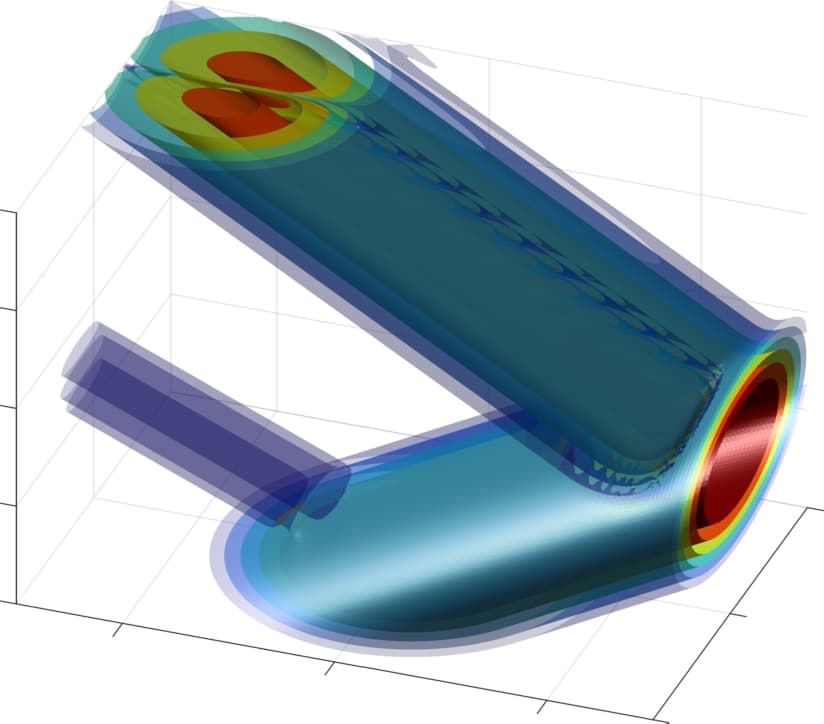}{\newsizee}{\originale}
\end{minipage}
\caption{Evolution of $\log_{10} |u_{j,k}^n|$ (left: six successive times, right: full time evolution) for $\mathbf{c}=(1,0.1)$ and DTBC of order $0$.}
\label{fig:test2D-2}
\end{figure}

\begin{figure}[htbp]
\begin{minipage}[htbp]{.5\linewidth}
  \begin{tabular}[htbp]{cc}
    \logplottext{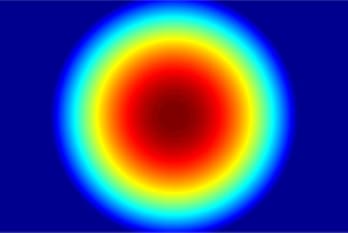}{\newsize}{\original}{\mascale}{t=0}   &
    \logplottext{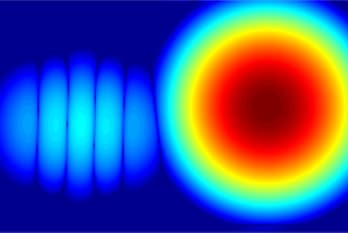}{\newsize}{\original}{\mascale}{t=1.6} \\    
    \logplottext{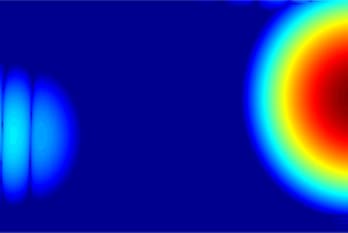}{\newsize}{\original}{\mascale}{t=3.2} &
    \logplottext{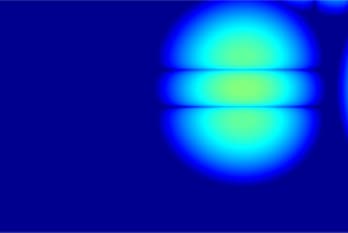}{\newsize}{\original}{\mascale}{t=4.8} \\    
    \logplottext{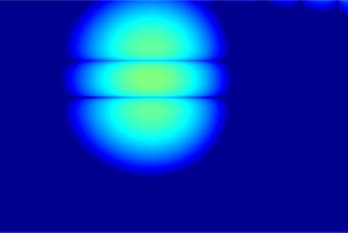}{\newsize}{\original}{\mascale}{t=6.4} &
    \logplottext{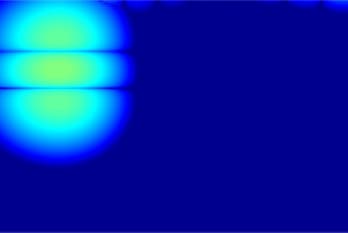}{\newsize}{\original}{\mascale}{t=8}   \\    
  \end{tabular}
\end{minipage}
\begin{minipage}[htbp]{.5\linewidth}
  \def\newsizee{.8\textwidth}
  \myTDrepr{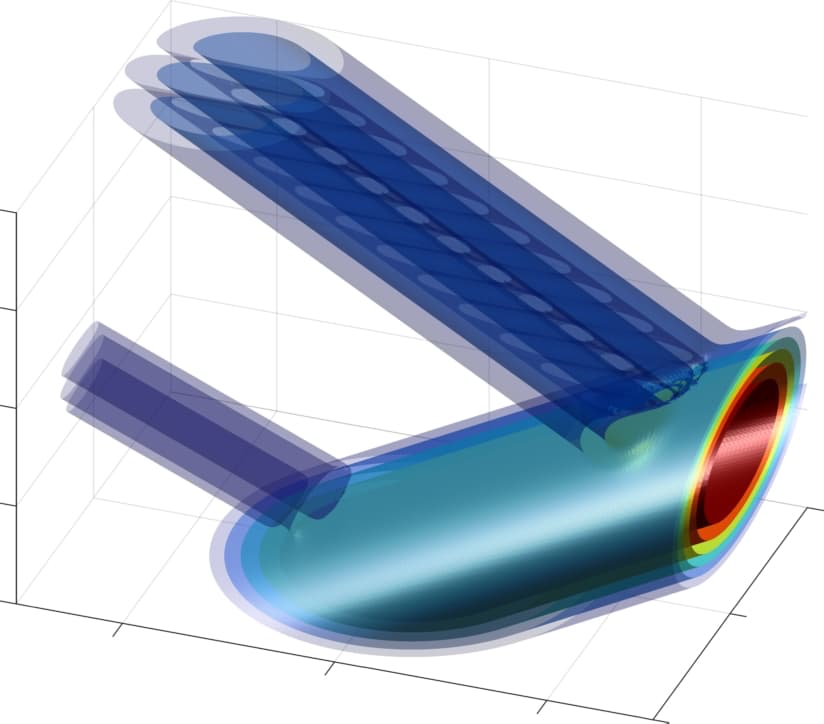}{\newsizee}{\originale}
\end{minipage}
\caption{Evolution of $\log_{10} |u_{j,k}^n|$ (left: six successive times, right: full time evolution) for $\mathbf{c}=(1,0.1)$ and DTBC of order $1$.}
\label{fig:test2D-3}
\end{figure}

\begin{figure}[htbp]
\begin{minipage}[htbp]{.5\linewidth}
  \begin{tabular}[htbp]{cc}
    \logplottext{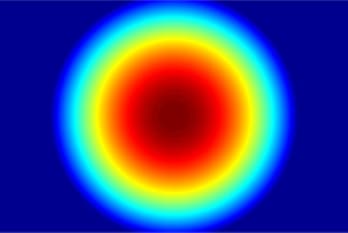}{\newsize}{\original}{\mascale}{t=0}   &
    \logplottext{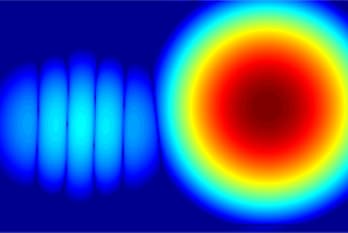}{\newsize}{\original}{\mascale}{t=1.6} \\    
    \logplottext{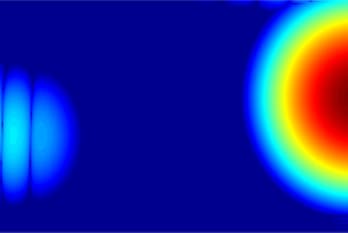}{\newsize}{\original}{\mascale}{t=3.2} &
    \logplottext{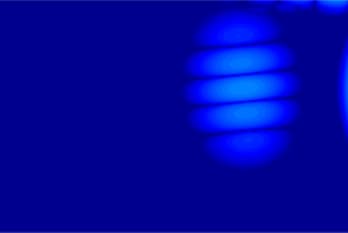}{\newsize}{\original}{\mascale}{t=4.8} \\    
    \logplottext{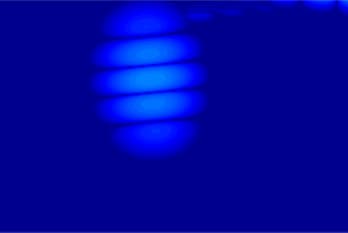}{\newsize}{\original}{\mascale}{t=6.4} &
    \logplottext{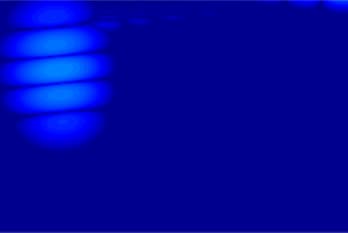}{\newsize}{\original}{\mascale}{t=8}   \\    
  \end{tabular}
\end{minipage}
\begin{minipage}[htbp]{.5\linewidth}
  \def\newsizee{.8\textwidth}
  \myTDrepr{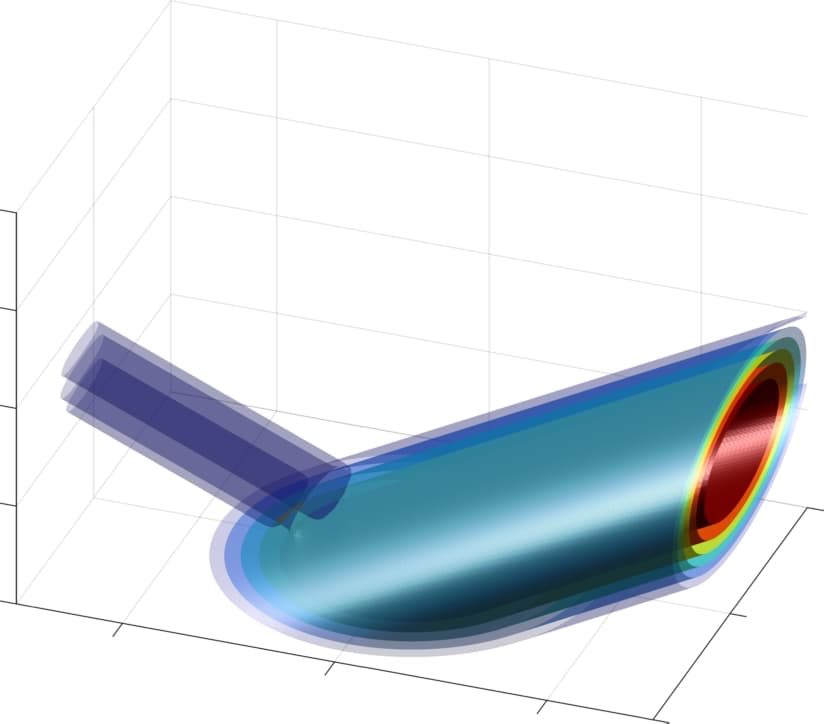}{\newsizee}{\originale}
\end{minipage}
\caption{Evolution of $\log_{10} |u_{j,k}^n|$ (left: six successive times, right: full time evolution) for $\mathbf{c}=(1,0.1)$ and DTBC of order $2$ 
in the $x-$direction and order $1$ in the $y-$direction.}
\label{fig:test2D-4}
\end{figure}

We first report on some cases where the initial Gaussian condition is transported towards one edge of the rectangle and does not give rise 
(in the considered time interval) to multiple wave reflections. The first case is depicted in Figure \ref{fig:test2D-1} where the chosen velocity 
is $\mathbf{c}=(1,0)$ (and therefore all boundary conditions \eqref{DTBC2d-right-a}, \eqref{DTBC2d-right-b}, \eqref{DTBC2d-right-c} coincide 
since $\mu_y=0$). We call them DTBC of order $0$ since \eqref{DTBC2d-right-a} amounts to not retaining any tangential dependence. The 
results are of course as accurate as in one space dimension (because the scheme \eqref{LF2d} is one-dimensional in that case and the 
tangential index $k$ only enters as a parameter). We now turn on the tangential dependence by considering the case $\mathbf{c}=(1,0.1)$. 
Both $\mu_x$ and $\mu_y$ are positive then. On the right and left faces of the rectangle, we either implement:

\begin{itemize}
  \item The approximate DTBC \eqref{DTBC2d-right-a}, \eqref{DTBC2d-left-a}. These are called DTBC of order $0$; accordingly our choice of approximate 
  DTBC for the top and bottom faces of the rectangle are \eqref{DTBC2d-top-a}, \eqref{DTBC2d-bottom-a}. The results are shown in Figure \ref{fig:test2D-2}. 
  The reflected wave at the (right) outflow boundary has magnitude $10^{-3}$.
  
  \item The approximate DTBC \eqref{DTBC2d-right-b}, \eqref{DTBC2d-left-b}. These are called DTBC of order $1$; accordingly our choice of approximate 
  DTBC for the top and bottom faces of the rectangle are \eqref{DTBC2d-top-b}, \eqref{DTBC2d-bottom-b}. The results are shown in Figure \ref{fig:test2D-3}. 
  The reflected wave at the (right) outflow boundary now has magnitude $10^{-5}$, which shows significant improvement.
  
  \item The approximate DTBC \eqref{DTBC2d-right-c}, \eqref{DTBC2d-left-c}. These are called DTBC of order $2$; there is a subtlety here because 
  our choice of approximate DTBC for the top and bottom faces of the rectangle are still \eqref{DTBC2d-top-b}, \eqref{DTBC2d-bottom-b}, and not 
  \eqref{DTBC2d-top-c}, \eqref{DTBC2d-bottom-c} as one might expect. This is further explained below. The results are shown in Figure \ref{fig:test2D-4}. 
  The reflected wave at the (right) outflow boundary now has magnitude $10^{-8}$ and is therefore not visible. This shows again significant improvement.
\end{itemize}

The overall conclusion that can be drawn from these numerical results is that, up to stability issues (which can be significant from an analytical point of 
view), adding more tangential terms in the Taylor expansion \eqref{eq:dtbc2d_2} seems to lead to more accurate results, which was to be expected of 
course since we consider smooth solutions. The numerical simulations quantify the gain of adding each new term in \eqref{DTBC2d-b} and 
\eqref{DTBC2d-c} with respect to \eqref{DTBC2d-a}.

\begin{figure}[htbp]
\begin{minipage}[htbp]{.5\linewidth}
  \begin{tabular}[htbp]{cc}
    \logplottext{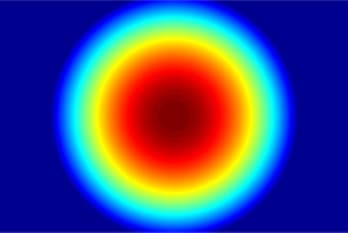}{\newsize}{\original}{\mascale}{t=0}   &
    \logplottext{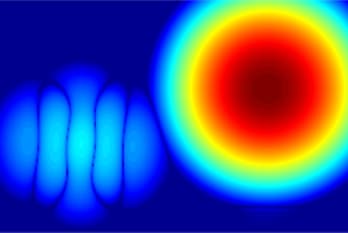}{\newsize}{\original}{\mascale}{t=1.6} \\    
    \logplottext{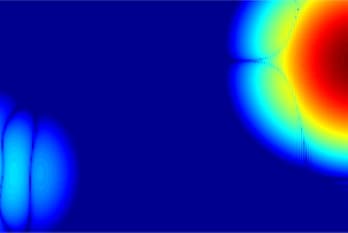}{\newsize}{\original}{\mascale}{t=3.2} &
    \logplottext{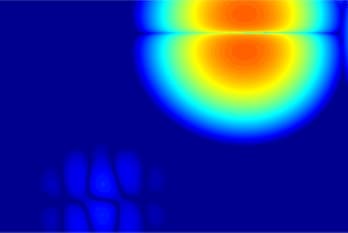}{\newsize}{\original}{\mascale}{t=4.8} \\    
    \logplottext{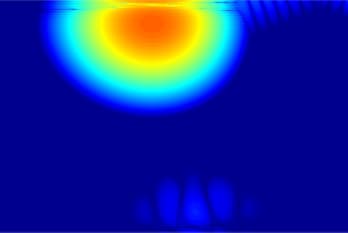}{\newsize}{\original}{\mascale}{t=6.4} &
    \logplottext{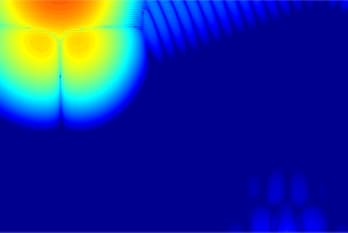}{\newsize}{\original}{\mascale}{t=8}   \\    
  \end{tabular}
\end{minipage}
\begin{minipage}[htbp]{.5\linewidth}
\def\newsizee{.8\textwidth}
\myTDrepr{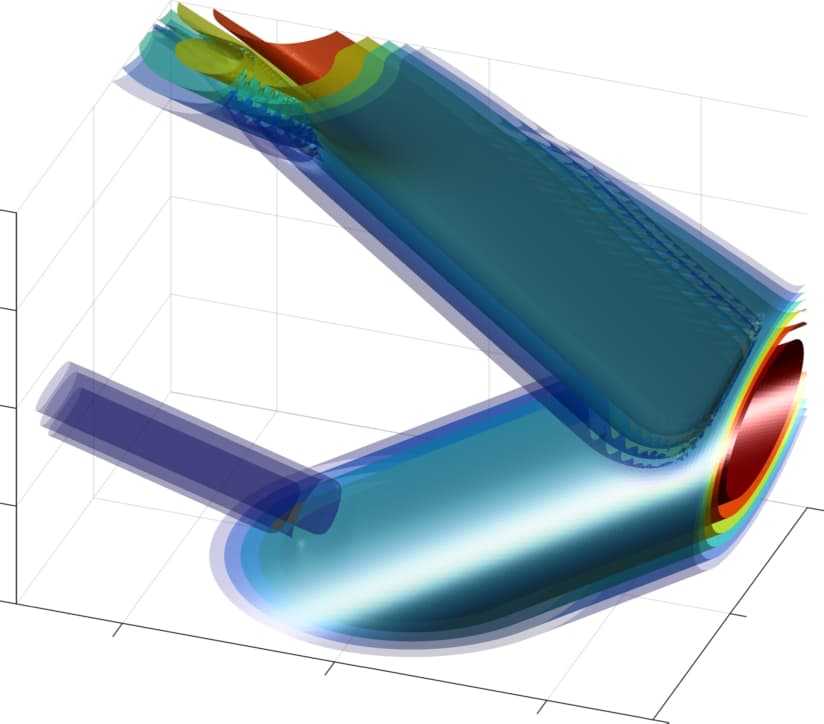}{\newsizee}{\originale}
\end{minipage}
\caption{Evolution of $\log_{10} |u_{j,k}^n|$ (left: six successive times, right: full time evolution) for $\mathbf{c}=(1,0.3)$ and DTBC of order $0$.}
\label{fig:test2D-5}
\end{figure}

\begin{figure}[htbp]
\begin{minipage}[htbp]{.5\linewidth}
  \begin{tabular}[htbp]{cc}
    \logplottext{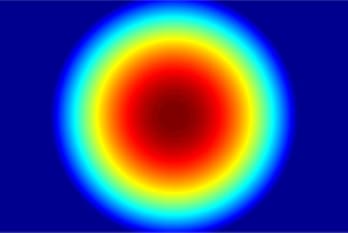}{\newsize}{\original}{\mascale}{t=0}   &
    \logplottext{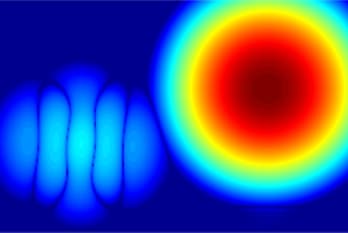}{\newsize}{\original}{\mascale}{t=1.6} \\    
    \logplottext{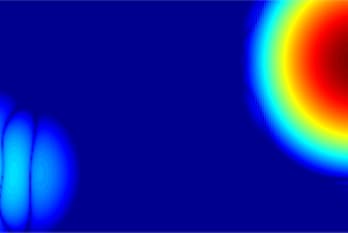}{\newsize}{\original}{\mascale}{t=3.2} &
    \logplottext{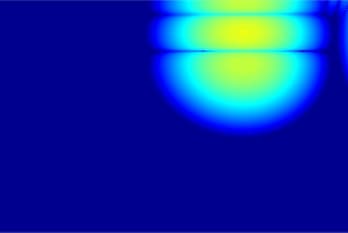}{\newsize}{\original}{\mascale}{t=4.8} \\    
    \logplottext{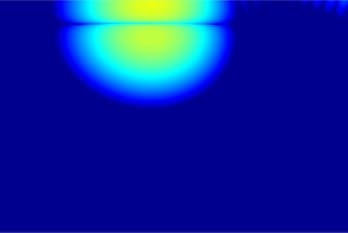}{\newsize}{\original}{\mascale}{t=6.4} &
    \logplottext{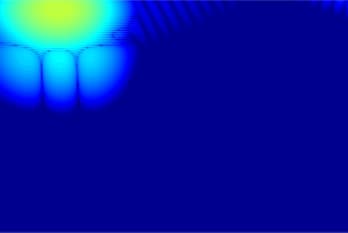}{\newsize}{\original}{\mascale}{t=8}   \\    
  \end{tabular}
\end{minipage}
\begin{minipage}[htbp]{.5\linewidth}
  \def\newsizee{.8\textwidth}
  \myTDrepr{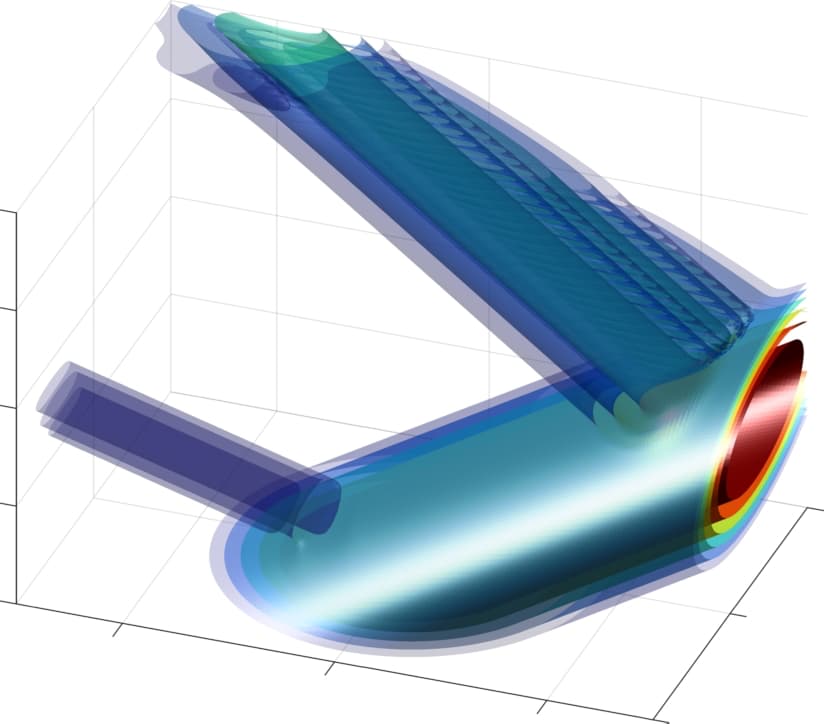}{\newsizee}{\originale}
\end{minipage}
\caption{Evolution of $\log_{10} |u_{j,k}^n|$ (left: six successive times, right: full time evolution) for $\mathbf{c}=(1,0.3)$ and DTBC of order $1$.}
\label{fig:test2D-6}
\end{figure}

\begin{figure}[htbp]
\begin{minipage}[htbp]{.5\linewidth}
  \begin{tabular}[htbp]{cc}
    \logplottext{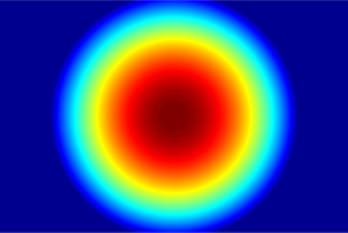}{\newsize}{\original}{\mascale}{t=0}   &
    \logplottext{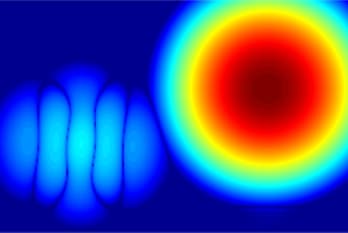}{\newsize}{\original}{\mascale}{t=1.6} \\    
    \logplottext{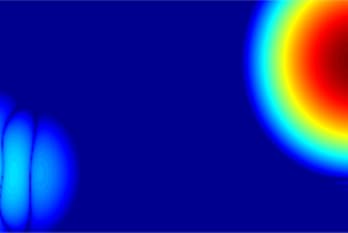}{\newsize}{\original}{\mascale}{t=3.2} &
    \logplottext{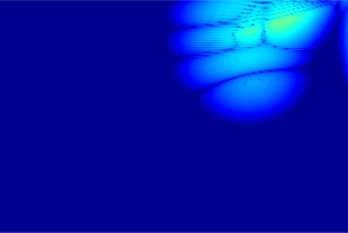}{\newsize}{\original}{\mascale}{t=4.8} \\    
    \logplottext{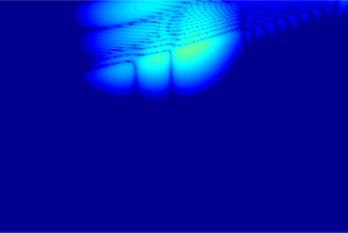}{\newsize}{\original}{\mascale}{t=6.4} &
    \logplottext{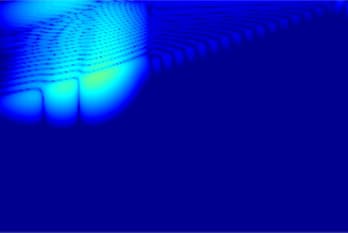}{\newsize}{\original}{\mascale}{t=8}   \\    
  \end{tabular}
\end{minipage}
\begin{minipage}[htbp]{.5\linewidth}
  \def\newsizee{.8\textwidth}
  \myTDrepr{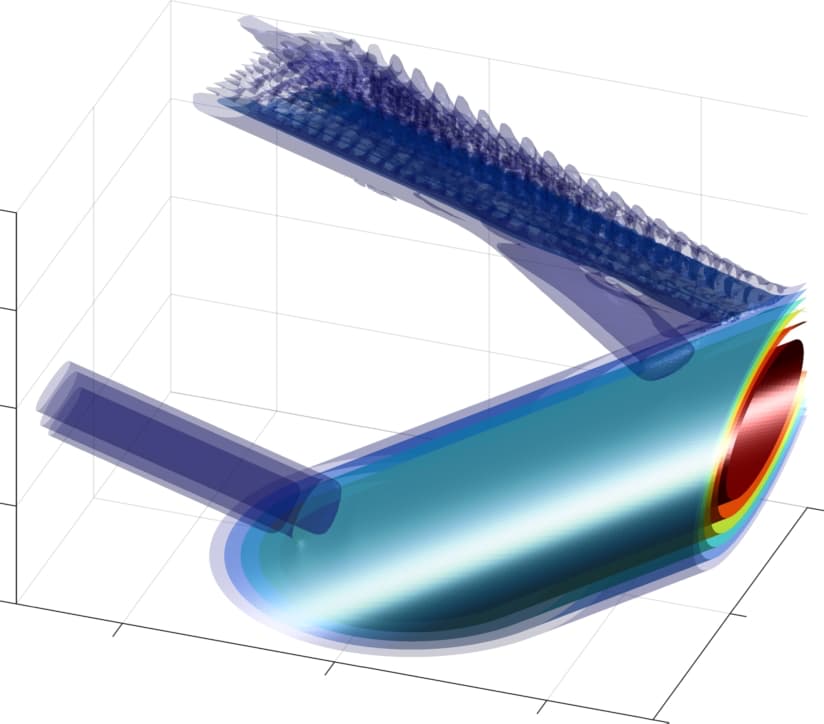}{\newsizee}{\originale}
\end{minipage}
\caption{Evolution of $\log_{10} |u_{j,k}^n|$ (left: six successive times, right: full time evolution) for $\mathbf{c}=(1,0.1)$ and DTBC of order $2$ 
in the $x-$direction and order $1$ in the $y-$direction.}
\label{fig:test2D-7}
\end{figure}

We now report on more intricate cases that exhibit multiple wave reflections. Namely, we now choose the velocity $\mathbf{c}=(1,0.3)$. With this 
velocity, the initial Gaussian function should leave the rectangle through the right edge of the boundary. Because the numerical boundary conditions 
\eqref{eq:dtbc2d_right} are not exactly transparent, wave reflection occurs. For the leap-frog scheme, the reflection takes place according to the 
billiard law. A (small amplitude highly oscillating) wave then moves towards the top face of the rectangle and a second wave reflection occurs. For 
this second reflection, the small tangential frequency assumption that motivated the expansion \eqref{eq:dtbc2d_2} is not any longer valid. Hence 
it is likely that there will be little absorption through the top boundary. The numerical simulations confirm that expectation. We follow the same plan 
as for the case $\mathbf{c}=(1,0.1)$ and report on some simulations with increasing tangential approximation. We thus implement:

\begin{itemize}
 \item The zero order DTBC \eqref{DTBC2d-right-a}, \eqref{DTBC2d-left-a}, \eqref{DTBC2d-top-a}, \eqref{DTBC2d-bottom-a} on the four sides of the rectangle. 
 This case is shown in Figure \ref{fig:test2D-5}; after the first reflection on the right face, the wave has magnitude $10^{-3}$. The magnitude remains of order 
 $10^{-3}$ after the second reflection on the top boundary.
 
 \item The first order DTBC \eqref{DTBC2d-right-b}, \eqref{DTBC2d-left-b}, \eqref{DTBC2d-top-b}, \eqref{DTBC2d-bottom-b} on the four sides of the rectangle. 
 This case is shown in Figure \ref{fig:test2D-6}; after the first reflection on the right face, the wave has magnitude $10^{-5}$. The magnitude remains of order 
 $10^{-5}$ after the second reflection on the top boundary. The overall accuracy is improved.
 
 \item The second order DTBC \eqref{DTBC2d-right-c}, \eqref{DTBC2d-left-c} on the right and left boundaries with the first order DTBC \eqref{DTBC2d-top-b}, 
 \eqref{DTBC2d-bottom-b} on the top and bottom boundaries. This case is shown in Figure \ref{fig:test2D-7}; after the first reflection on the right face, the wave 
 has magnitude $10^{-6}$. The magnitude remains of order $10^{-6}$ after the second reflection on the top boundary. The overall accuracy is improved again.
\end{itemize}

\begin{figure}[htbp]
\begin{minipage}[htbp]{.5\linewidth}
  \begin{tabular}[htbp]{cc}
    \logplottext{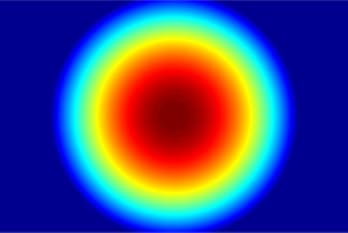}{\newsize}{\original}{\mascale}{t=0}   &
    \logplottext{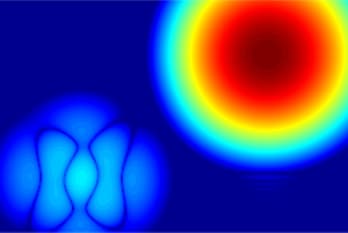}{\newsize}{\original}{\mascale}{t=1.6} \\    
    \logplottext{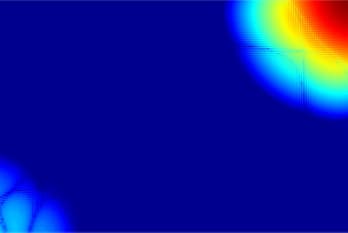}{\newsize}{\original}{\mascale}{t=3.2} &
    \logplottext{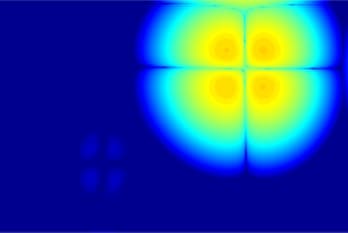}{\newsize}{\original}{\mascale}{t=4.8} \\    
    \logplottext{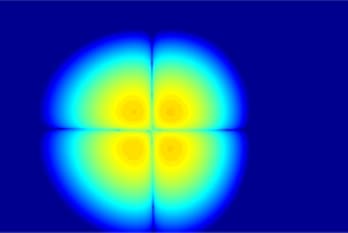}{\newsize}{\original}{\mascale}{t=6.4} &
    \logplottext{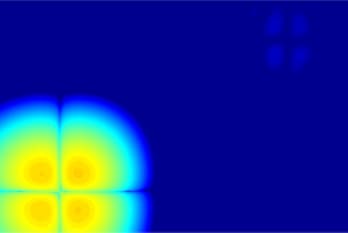}{\newsize}{\original}{\mascale}{t=8}   \\    
  \end{tabular}
\end{minipage}
\begin{minipage}[htbp]{.5\linewidth}
  \def\newsizee{.8\textwidth}
  \myTDrepr{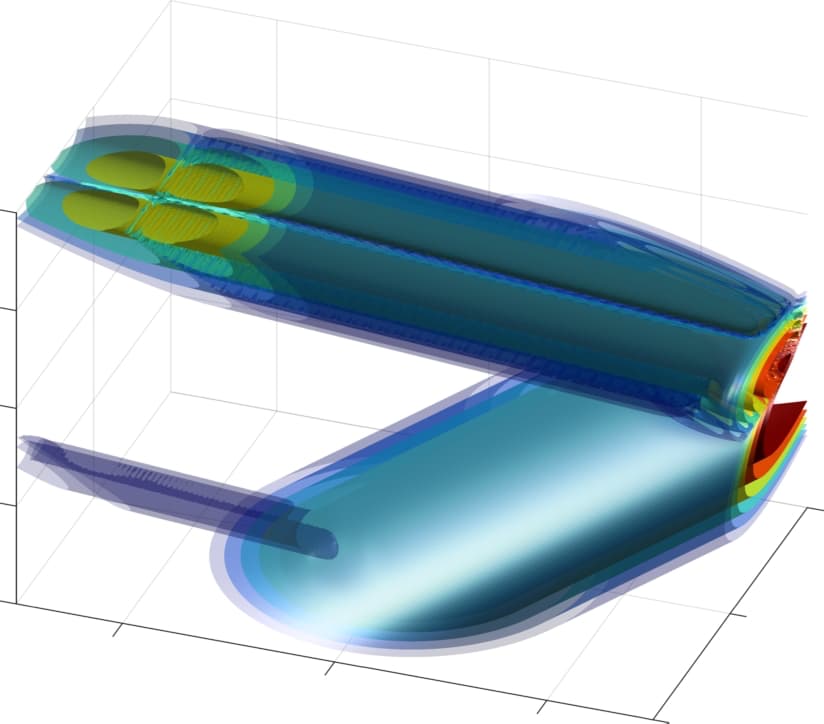}{\newsizee}{\originale}
\end{minipage}
\caption{Evolution of $\log_{10} |u_{j,k}^n|$ (left: six successive times, right: full time evolution) for $\mathbf{c}=(1,2/3)$ and DTBC of order $0$.}
\label{fig:test2D-8}
\end{figure}

\begin{figure}[htbp]
\begin{minipage}[htbp]{.5\linewidth}
  \begin{tabular}[htbp]{cc}
    \logplottext{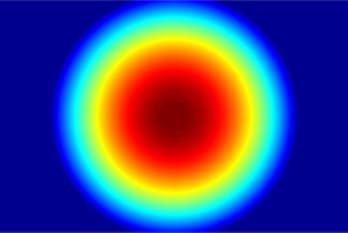}{\newsize}{\original}{\mascale}{t=0}   &
    \logplottext{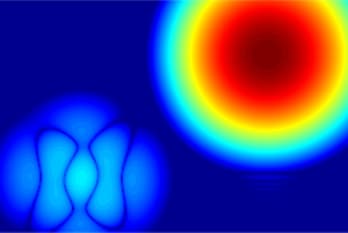}{\newsize}{\original}{\mascale}{t=1.6} \\    
    \logplottext{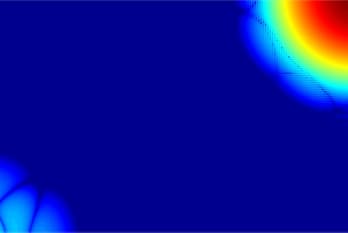}{\newsize}{\original}{\mascale}{t=3.2} &
    \logplottext{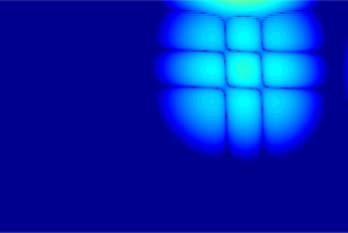}{\newsize}{\original}{\mascale}{t=4.8} \\    
    \logplottext{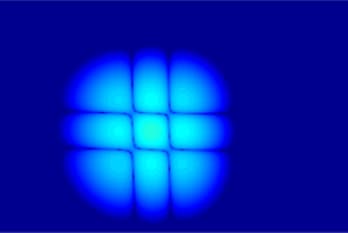}{\newsize}{\original}{\mascale}{t=6.4} &
    \logplottext{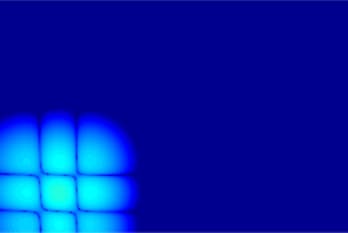}{\newsize}{\original}{\mascale}{t=8}   \\    
  \end{tabular}
\end{minipage}
\begin{minipage}[htbp]{.5\linewidth}
  \def\newsizee{.8\textwidth}
  \myTDrepr{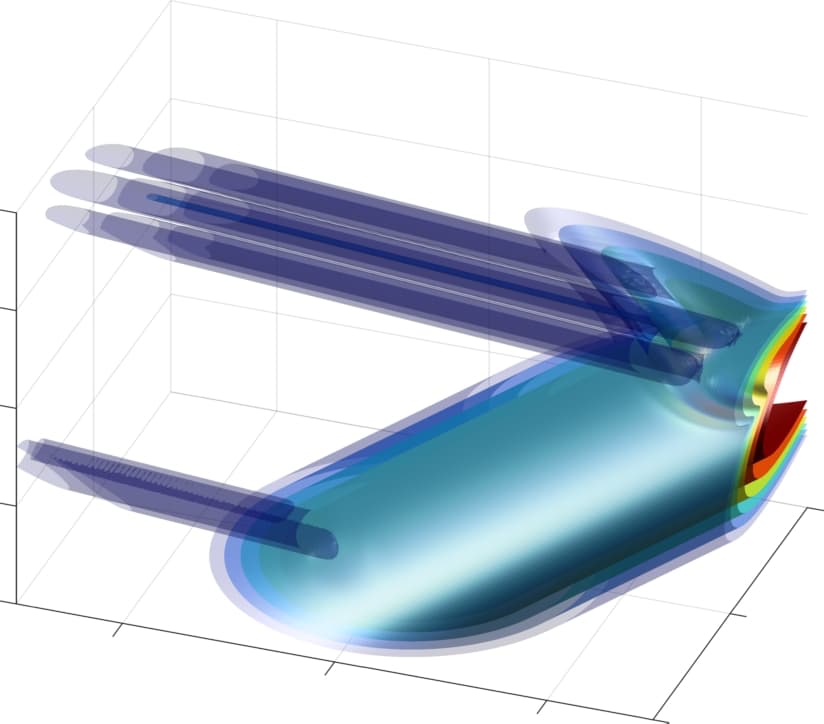}{\newsizee}{\originale}
\end{minipage}
\caption{Evolution of $\log_{10} |u_{j,k}^n|$ (left: six successive times, right: full time evolution) for $\mathbf{c}=(1,2/3)$ and DTBC of order $1$.}
\label{fig:test2D-9}
\end{figure}

\begin{figure}[htbp]
\begin{minipage}[htbp]{.5\linewidth}
  \begin{tabular}[htbp]{cc}
    \logplottext{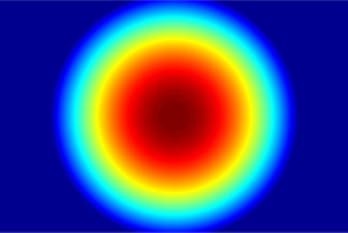}{\newsize}{\original}{\mascale}{t=0}   &
    \logplottext{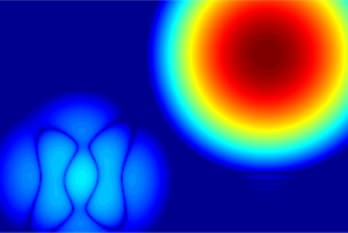}{\newsize}{\original}{\mascale}{t=1.6} \\    
    \logplottext{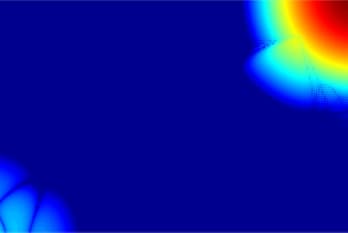}{\newsize}{\original}{\mascale}{t=3.2} &
    \logplottext{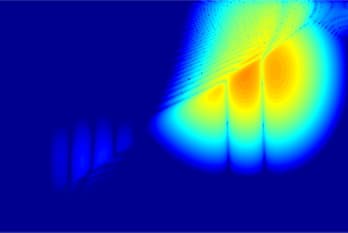}{\newsize}{\original}{\mascale}{t=4.8} \\    
    \logplottext{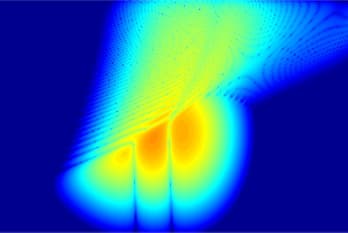}{\newsize}{\original}{\mascale}{t=6.4} &
    \logplottext{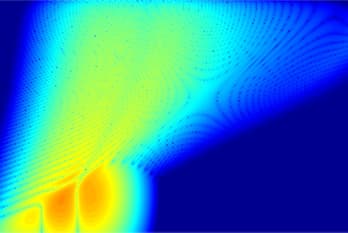}{\newsize}{\original}{\mascale}{t=8}   \\    
  \end{tabular}
\end{minipage}
\begin{minipage}[htbp]{.5\linewidth}
  \def\newsizee{.8\textwidth}
  \myTDrepr{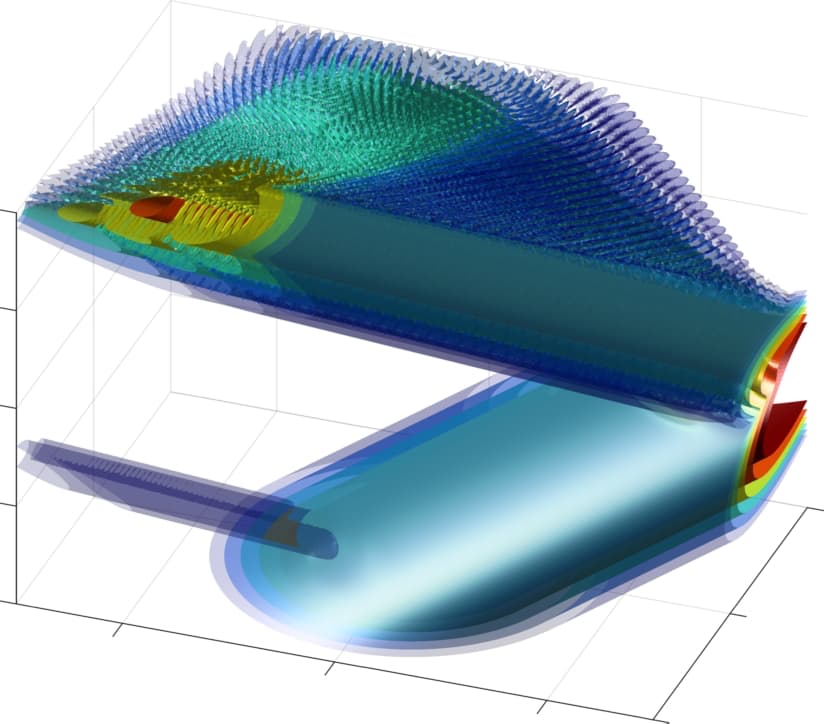}{\newsizee}{\originale}
\end{minipage}
\caption{Evolution of $\log_{10} |u_{j,k}^n|$ (left: six successive times, right: full time evolution) for $\mathbf{c}=(1,2/3)$ and DTBC of order $2$ 
in the $x-$direction and order $1$ in the $y-$direction.}
\label{fig:test2D-10}
\end{figure}

Figures \ref{fig:test2D-8}, \ref{fig:test2D-9} and \ref{fig:test2D-10} represent similar simulations as Figures \ref{fig:test2D-5}, \ref{fig:test2D-6} and \ref{fig:test2D-7} 
with the velocity $\mathbf{c}=(1,2/3)$. In that case, the initial condition is transported towards the upper right corner of the computational domain. There is again 
an improvement when passing from zero order DTBC to first order DTBC (compare the scales between Figures \ref{fig:test2D-8} and \ref{fig:test2D-9}). Figure 
\ref{fig:test2D-10} shows that coupling approximate DTBC with different tangential approximation orders on the two sides of a corner may affect the overall 
accuracy. Namely, for the velocity $\mathbf{c}=(1,2/3)$, using DTBC of order $1$ on the four sides of the rectangle yields better results than coupling DTBC 
of order $2$ on two opposite sides with DTBC of order $1$ on the two remaining sides.

\begin{figure}[htbp]
\begin{minipage}[htbp]{.45\linewidth}
  \input{evol_l2_O3.tex}
\end{minipage}
\begin{minipage}[htbp]{.55\linewidth}
          \begin{tikzpicture}
            \begin{axis}[
              xlabel=$x$,ylabel=$t$,
              enlargelimits=false,
              axis on top, width=.81\textwidth,
              colorbar,
              colorbar style={point meta min=-20,point meta max=14.4747,    ytick pos=right,
                tick label style={font=\footnotesize},
              },
              colorbar/width=3mm,    
              ]
              \addplot graphics [
              xmin=-3,xmax=3,
              ymin=-2,ymax=2,
              ] {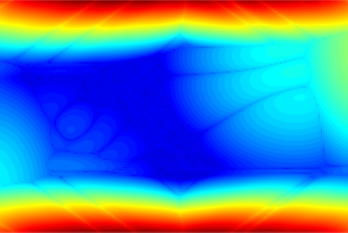};
            \end{axis}
          \end{tikzpicture}
\end{minipage}
\caption{The instability for $\mathbf{c}=(1,0.3)$ and DTBC of order $2$ in both $x$ and $y$ directions. Left: time evolution of the logarithm of 
the $\ell^2$-norm of the solution. Right: the unstable profile at $t=4$ (logarithmic scale).}
\label{fig:test2D-11}
\end{figure}

Eventually, we show in Figure \ref{fig:test2D-11} the evolution of the $\ell^2$ norm of $(u_{j,k}^n)$ as a function of the time $t^n$ in the case 
$\mathbf{c}=(1,0.3)$ when the second order DTBC \eqref{DTBC2d-right-c}, \eqref{DTBC2d-left-c}, \eqref{DTBC2d-top-c}, \eqref{DTBC2d-bottom-c} 
are implemented on the four sides of the rectangle. The computation exhibits two stages. In the first stage, the initial condition (mostly) exits through 
the right boundary. As long as the solution remains small on the boundary, the effect of the numerical boundary conditions is rather invisible and the 
$\ell^2$ norm remains approximately constant. In a second stage, the trace of the solution is no longer small and an exponential instability takes 
place. The logarithmic scale in Figure \ref{fig:test2D-11} clearly shows the exponential behavior in time. The numerical solution displays a profile 
that is concentrated along the top and right boundaries (right of Figure \ref{fig:test2D-11}). It thus seems clear that the coupling between 
\eqref{DTBC2d-right-c} and \eqref{DTBC2d-top-c} gives rise to an unstable eigenvalue that is associated with a `surface wave' that has exponential 
decay with respect to the normal directions to both the right and top boundaries. We have unfortunately not been able to obtain an analytical proof 
of this fact, but our conclusion is that despite good stability properties of \eqref{DTBC2d-right-c} and \eqref{DTBC2d-top-c} when considered separately 
on each side of the rectangle, the coupling at the corner may yield strong instabilities.

\renewcommand\figurescale{0.4}

At last, we have also tested the efficiency of the sum of exponential approximation compared to the standard convolution procedure in the case of DTBC 
of order $1$ \eqref{DTBC2d-right-b}, \eqref{DTBC2d-left-b}, \eqref{DTBC2d-top-b}, \eqref{DTBC2d-bottom-b}. Namely, we have compared the computational 
effort in that case with the sum of exponential approximation of both sequences $(s_n^0)$ and $(s_n^1)$. We have chosen $(M,N)=(50,20)$ for the degrees 
of the Pad\'e approximant. The results are presented in Table \ref{tab:tab01}. Though there is little gain for the sum of exponential approximation in one space 
dimension, there is a factor $2$ gain in terms of computational effort in two space dimensions for an overall accuracy that is comparable.

\setlength{\tabcolsep}{6pt}
\renewcommand{\arraystretch}{1}

\begin{table}[htbp]
  \centering
  \begin{tabular}{|c|c|c|} \hline
    & Minimal cputime & Mean cputime \\ \hline
    DTBC & $21.9850$ sec & $22.2692$ sec \\ \hline 
SumExp & $9.7760 $ sec & $10.0946$ sec \\  \hline
  \end{tabular}
  \caption{Comparison of cputime between the implementation of DTBC and the implementation of the sum of exponential approximation in two space dimensions.}
  \label{tab:tab01}
  \end{table}

For the reader's convenience, all numerical codes are made available on the web page \url{http://nabuco.math.cnrs.fr}.
  
\bibliographystyle{alpha}
\bibliography{Transparent}
\end{document}